\pdfoutput=1

\documentclass[11pt]{article}
\RequirePackage[OT1]{fontenc}
\RequirePackage{amsthm,amsmath}
\RequirePackage[numbers]{natbib}
\RequirePackage[colorlinks,citecolor=blue,urlcolor=blue]{hyperref}

\usepackage{graphicx,subfigure}
\usepackage{amssymb,amsfonts}
\usepackage{enumitem}
\usepackage{algorithm}
\usepackage{algorithmic}

\newtheorem{theorem}{Theorem}

\newtheorem{proposition}{Proposition}

\def \bal#1\eal{\begin{align}#1\end{align}}
\def \bals#1\eals{\begin{align*}#1\end{align*}}
\def\P{\mathbb{P}}

\def\C{\mathcal{C}}
\def\real{\mathbb{R}}

\def\XPa{X_{\textnormal{PA}(i)}}
\def\XV{X_{V_i}}
\def\Pe{\P_{\exp}}
\def\CV{\C_{V_i}}

\def\Cmin{\rho_{\min}}

\def\by{{ Y}}
\def\bys{\widebar{ Y}_s}
\def\byi{\widebar{ Y}_1}

\def\regV{\lambda}
\def\regPa{\mu}

\newlength{\widebarargwidth}
\newlength{\widebarargheight}
\newlength{\widebarargdepth}
\DeclareRobustCommand{\widebar}[1]{%
  \settowidth{\widebarargwidth}{\ensuremath{#1}}%
  \settoheight{\widebarargheight}{\ensuremath{#1}}%
  \settodepth{\widebarargdepth}{\ensuremath{#1}}%
  \addtolength{\widebarargwidth}{-0.3\widebarargheight}%
  \addtolength{\widebarargwidth}{-0.3\widebarargdepth}%
  \makebox[0pt][l]{\hspace{0.15\widebarargheight}%
    \hspace{0.3\widebarargdepth}%
    \addtolength{\widebarargheight}{0.3ex}%
    \rule[\widebarargheight]{0.95\widebarargwidth}{0.1ex}}%
  {#1}}

\newcounter{NLCTR}
\newenvironment{mylist}
{\begin{list}{\Roman{NLCTR}.}
	{
	\vskip-0.1in
	\usecounter{NLCTR}
	\setlength{\parsep}{0.05in} \setlength{\itemsep}{\parskip}
	\setlength{\leftmargin}{0.0in} \setlength{\rightmargin}{0.08in}
	\setlength{\listparindent}{0in}
	\setlength{\labelwidth}{-0.045in}
	\setlength{\labelsep}{0.05in} \setlength{\itemindent}{0in}}}
{\end{list}}

\usepackage{setspace}

\begin{document}

\title{A General Framework for Mixed Graphical Models}
\date{}
\author{
    Eunho Yang$^{1}$, Pradeep Ravikumar$^{2}$, Genevera I. Allen$^{3,4}$,\\ Yulia Baker$^{3}$, Ying-Wooi Wan$^{4}$ and Zhandong Liu$^{4}$\\
    {\small $^{1}$IBM T.J. Watson Research Center,}\\
    {\small $^{2}$Department of Computer Science, University of Texas Austin,}\\
    {\small $^{3}$ Department of  Statistics, Rice University,} \\
    {\small $^{4}$ Department of Pediatrics-Neurology, Baylor College of Medicine}
}
\maketitle
%
%
%
%
%
%
%
%
%

\begin{abstract}
``Mixed Data'' comprising a large number of heterogeneous variables
  (e.g.  count, binary, continuous, skewed continuous, among other
  data types) are prevalent in varied areas such as genomics and
  proteomics, imaging genetics, national security, social networking,
  and Internet advertising.  There have been limited efforts at
  statistically modeling such mixed data jointly, in part because of
  the lack of computationally amenable multivariate distributions that
  can capture direct dependencies between such mixed variables of
  different types.  In this paper, we address this by introducing a
  novel class of Block Directed Markov Random Fields (BDMRFs).  Using
  the basic building block of node-conditional univariate exponential
  families from Yang et al. (2012), we introduce a class of mixed
  conditional random field distributions, that are then chained
  according to a block-directed acyclic graph to form our class of
  Block Directed Markov Random Fields (BDMRFs). The Markov
  independence graph structure underlying a BDMRF thus has both
  directed and undirected edges. We introduce conditions under which
  these distributions exist and are normalizable, study several
  instances of our models, and propose scalable penalized conditional
  likelihood estimators with statistical guarantees for recovering the
  underlying network structure. Simulations as well as an application
  to learning mixed genomic networks from next generation sequencing
  expression data and mutation data demonstrate the versatility of our
  methods.  
\end{abstract}

%
%
%

\section{Introduction}\label{Sec:Intro}

\subsection{Motivation}
Acquiring and storing data is steadily becoming cheaper in this Big
Data era. A natural consequence of this is more varied and complex
data sets, that consist of variables of {\em mixed types}; we refer to
mixed types as variables measured on the same set of samples that can
each belong to differing domains such as 
binary, categorical, ordinal, counts, continuous, and/or skewed
continuous, among others.
Consider, for instance, three popular Big Data applications, all of
which are comprised of such mixed data:  
\begin{itemize}
\item \emph{High-throughput biomedical data.}  With new biomedical
technologies, scientists can now measure full genomic, proteomic, and
metabolomic scans as well as biomedical imaging on a single subject or
biological sample.   
\item \emph{National security data.}  Technologies exist to collect varied
information such as call-logs, 
geographic coordinates, text messages, tweets, demographics,
purchasing history, and Internet browsing history, among others.  
\item \emph{Internet-scale marketing data.}  Internet companies in an effort to
optimize advertising revenues collect information such as Internet browsing
history, social media postings, friends in social networks, status updates,
tweets, purchasing history, ad-click history, and online video viewing
history, among others.  
\end{itemize}
In each of these examples, variables of many different types are
collected on the same samples, and these variables are clearly
dependent. Consider the motivating example of high-throughput
``omics'' data in further detail.  There has been a recent
proliferation of genomics technologies that can measure nearly every
molecular aspect of a given sample.  These technologies, however,
produce mixed types of variables: mutations and aberrations such as
SNPs and copy number variations are  binary or categorical, functional
genomics comprising gene expression and miRNA expression as measured
by RNA-sequencing are count-valued, and epigenetics as measured by
methylation arrays are continuous. Clearly, all of these genomic
biomarkers are closely inter-related as they belong to the same
complex biological system. Multivariate distributions such as
graphical models applied to one type of data, typically gene
expression, are popularly used, for tasks ranging from data
visualization, finding important biomarkers, and estimating regulatory
relationships.  Yet, to understand the complete molecular basis of
diseases, requires us to understand relationships not only {\em
  within} a specific type of biomarkers, but also {\em between}
different types of biomarkers. Thus, developing a class of
multivariate distributions that can directly model dependencies
between genes based on gene expression levels (counts) as well as the
mutations (binary) that influence gene expression levels is needed, to
holistically model the genomic system. Few such multivariate
distributions exist, however, that can directly model dependencies
between mixed types 
of variables. In this paper, our goal is to address this lacuna, and
define a parametric family of multivariate distributions, that can
model rich dependence structures over mixed variables. We leverage
theory of exponential family distributions to do so, and term our
resulting novel class of statistical models as mixed Block Directed Markov
Random Fields (BDMRFs). 


A popular class of statistical studies of such mixed, multivariate data sidestep multivariate densities altogether. Instead, they relate a set of multivariate response variables of one type to multivariate covariate variables of another type, using multiple regression models or multi-task learning models \citep{argyriou2007multi}.  These are especially popular in 
expression quantitative trait loci (eQTL) analyses, which seek to link changes in functional gene expression levels to specific genomic mutations \citep{kim2009statistical}. Recent approaches~\citep{yang2009heterogeneous} further allow these multiple regression models to associate covariates with mixed types of responses. More general regression and predictive models such as Classification and Regression Trees have also been proposed for such mixed types of covariates~\citep{esl_2nd}. Other approaches implicitly account for variables of mixed types in many machine learning procedures using suitable distance or entropy-based measures~\citep{esl_2nd,hsu2007hierarchical}.  There have also been non-parametric extensions of probabilistic graphical models using copulas \cite{dobra_2011,liu2012high} or rank-based estimators \cite{xue2012regularized,liu_nonpara_2009}, which could potentially be used for mixed data; non-parametric methods, however, may suffer from a loss of statistical efficiency when compared to parametric families, especially under very high-dimensional sampling regimes.  Others have proposed to build network models based on random forests \cite{fellinghauer2013stable} which are able to handle mixed types of variables, but these do not correspond to a multivariate density.

Among parametric statistical modeling approaches for such mixed data, the most popular, especially in survey statistics and spatial statistics \cite{cressie1993statistics}, are hierarchical models that permit dependencies through latent variables. For example, \citet{sammel1997latent} propose a latent variable model for mixed continuous and count variables, while \citet{rue2009approximate} propose latent Gaussian models that permit dependencies through a latent Gaussian MRF.  While these methods provide statistical models for mixed data, they model dependencies between observed variables via a latent layer that is not observed; estimating these models with strong statistical guarantees is thus typically computationally expensive and possibly intractable.


In this paper, we seek to specify parametric multivariate
distributions over mixed types of variables, that directly model
dependencies among these variables, without recourse to latent
variables, and that are computationally tractable with statistical
guarantees for high-dimensional data. Due in part to its
importance, there has been some recent set of proposals towards such
direct parametric statistical models, building on some seminal earlier
work by \citet{lauritzen1989graphical}. We review these in the next
sub-section after first providing some background on Markov Random
Fields (MRFs). As we will show in the next sub-section, these are
however largely targeted to the case where there are variables of two
types: discrete and continuous. Some other recent proposals, including
some of our prior work, consider more general mixed multivariate
distributions, but as we will show these are not sufficiently
expressive to allow for rich dependencies between disparate data
types. Our proposed class of Block Directed Markov Random Fields
(BDMRFs) serve as a vast generalization of these proposals, and indeed
as the teleological porting of the work by
\citet{lauritzen1989graphical} to the completely heterogeneous mixed
data setting. 

\subsection{Background and Related Work}
\label{Sec:Background}

\subsubsection{Markov Random Fields} Suppose $X = (X_1,\hdots,X_p)$ is
a random vector, with each variable $X_i$ taking values in a set
$\mathcal{X}$. Let $G_X = (V_X, E_X)$ be an undirected graph over $p$
nodes corresponding to the $p$ variables $\{X_s\}_{s=1}^{p}$. The
graphical model over $X$ corresponding to $G_X$ is \emph{a set of
distributions} that satisfy \emph{Markov independence assumptions}
with respect to the graph $G_X$~\citep{Lauritzen}. By the
Hammersley-Clifford theorem~\citep{Clifford90}, any such distribution
that is strictly positive over its domain also factors according to
the graph in the following way. Let $\C_X$ be a set of cliques
(fully-connected subgraphs) of the graph $G_X$, and let
$\{ \phi_{c}(X_c) \}_{c \in \C_X}$ be a set of clique-wise sufficient
statistics. Then, any strictly positive distribution of $X$ within the
graphical model family represented by the graph $G_X$ takes the form: 
\begin{eqnarray}
\label{EqnGenMRFPot}
	\P[X] & \propto & \exp \biggr \{\sum_{c \in \C_X}
\theta_c \phi_c(X_c) \biggr \},
\end{eqnarray}
where $\{\theta_c\}$ are weights over the sufficient statistics. Popular instances of this model include Ising models \citep{WaiJor08} for discrete-valued qualitative variables, and Gaussian MRFs~\citep{Speed86} for continuous-valued quantitative variables. Ising models specify joint distributions over a set of binary variables each with domain $\mathcal{X} = \{0,1\}$, with the form
\begin{align*}
	\P[X] \propto \exp\bigg\{\sum_{(s,t) \in E_X}\theta_{st} \, X_s \, X_{t} \bigg\},
\end{align*}
where we have ignored singleton terms for simplicity. Gaussian MRFs on the other hand specify joint distributions over a set of continuous real-valued variables each with domain $\mathcal{X} = \real$, with the form
\begin{align}\label{EqnGMRF}
	&\P[X] \propto \exp\Bigg\{ \sum_{s \in V_X} \frac{\theta_s}{\sigma_s} X_s + \sum_{(s,t) \in E_X}\frac{\theta_{st}}{\sigma_s\sigma_t} \, X_s \, X_t - \sum_{s \in V_X}\frac{X_s^2}{2\sigma_s^2} \Bigg\}.
\end{align}

\subsubsection{Conditional Gaussian Models}
We now review conditional Gaussian models, which were the first proposed class of mixed graphical models, introduced in \citep{lauritzen1989graphical}, and further studied in  \citep{frydenberg1989decomposition,lauritzen1992propagation,lauritzen1989mixed,Lauritzen}. Let $Y := (Y_1,\hdots,Y_p) \in \real^p$ be a continuous \emph{response} random vector, and let $X := (X_1,\hdots,X_q) \in \{1,\hdots,k\}^q$ be a discrete \emph{covariate} random vector. Taken together, $(X,Y)$ is then a \emph{mixed} random vector with both continuous and discrete components. For such a mixed random vector, \citep{lauritzen1989graphical} proposed the following joint distribution:
\begin{align}
	\P[X,Y] = \exp \left\{ g(X) + h(X)^T Y - \frac{Y^T K(X) Y}{2} \right\},
	\label{EqnCGModel}
\end{align}
parameterized by $g(X) \in \real$, $h(X) \in \real^p$, and $K(X) \in \real^{p \times p}$, such that $K(X) \succeq 0$. They termed this model a \emph{conditional Gaussian (CG)} model, since the conditional distribution of $\P[Y|X]$ given the joint distribution in \eqref{EqnCGModel} is given by a multivariate Gaussian distribution:
\begin{align*}
	\P[Y|X] \equiv \mathcal{N}(K(X)^{-1} h(X),\, K(X)^{-1}).
\end{align*}

Consider the set of vertices $V_X$ and $V_Y$ corresponding to random variables in $X$ and $Y$ respectively, and consider the set joint vertices $V = V_X \cup V_Y$. Would Markov assumptions with respect to a graph $G = (V,E)$ over all the vertices entail restrictions on $g(X), h(X), K(X)$ in the joint distribution in \eqref{EqnCGModel}? In the following theorem from \citet{lauritzen1989graphical}, which we reproduce for completeness, they provide an answer to this question:
\begin{theorem}
	If the CG distribution in \eqref{EqnCGModel} is Markov with respect to a graph $G$, then $g(X), h(X), K(X)$ can be written as 
	\begin{align*}
	g(X) = \sum_{d \subseteq [q]}g_{d}(X), \ 
	h(X) = \sum_{d \subseteq [q]}h_{d}(X), \
	K(X) = \sum_{d \subseteq [q]}K_{d}(X),
	\end{align*}
	where we use $[q]$ to denote $\{1,\hdots,q\}$ and 
	\begin{itemize}
	\item $g_{d}(X),h_d(X), K_{d}(X)$ depend on $X$ only through the subvector $X_{d}$, 
	\item for any subset $d \subseteq [q]$ which is not complete with respect to $G$, the corresponding components are zero, so that $g_{d}(X) = 0, h_{d}(X) = 0, K_{d}(X) = 0$.
	\end{itemize}
\end{theorem}

Recently, there have been several proposals for estimating the graph
structure of these CG models in high-dimensional
settings. \citet{lee2012learning} consider a specialization of CG
models involving only pairwise interactions between any two variables
and propose sparse node-wise estimators for graph
selection.  \citet{cheng2013high} further consider three-way
interactions between two binary variables and one continuous variable
and also propose sparse node-wise estimators to select the network
structure.  


\subsubsection{Graphical Models beyond Ising and Gaussian MRFs, via Exponential Families}

A key caveat with the CG models as defined in \eqref{EqnCGModel}
and \citep{lauritzen1989graphical}, however, is that it is
specifically constructed for mixed discrete and thin-tailed continuous
data. What if the responses and/or the covariates were count-valued,
or skewed-continuous, or belonged to some other non-categorical
non-thin-tailed-continuous data type? This is the question we address
in this paper. Towards this, we first briefly review here a recent
line of work~\citep{YRAL12,YRAL13CRF,YBRAL14} (which extends earlier
work by \citet{besag1974spatial}) which specified undirected graphical
model distributions where the variables all belong to one data-type,
but which could be any among a wide class of data-types. Their
development was as follows. Consider the general class of univariate
exponential family distributions (which include many popular
distributions such as Bernoulli, Gaussian, Poisson, negative binomial,
and exponential, among others): 
\begin{align}\label{EqnExpFamDist}
	\P[Z;\theta] = \exp \Big\{\theta \, B(Z) + M(Z) - D(\theta)\Big\},
\end{align}
with sufficient statistics $B(Z)$, base measure $M(Z)$, and log-normalization constant $D(\theta)$. Suppose that $p$-dimensional random vector $X$ has node-conditional distributions specified by an exponential family,
\begin{align*}
	\log \P[X_s | X_{V_X \backslash s}] = E(X_{V_X \backslash s})\, B(X_s) + M(X_s) - \bar{D}(X_{V_X \backslash s}) ,
\end{align*}
where the function $E(X_{V_X \backslash s})$ is any function that depends on the rest of all random variables except $X_s$. Further suppose that the corresponding joint distribution $\Pe[X]$ factors according to the set of cliques $\C_X$ of a graph $G_X$. \citet{YRAL12} then showed that such a joint distribution consistent with the above node-conditional distributions exists, and moreover necessarily has the form
\begin{align}\label{EqnExpMRF}
	\log \P_{\exp}[X] = \sum_{C \in \C_{X}} \theta_{C} \prod_{s \in C} B(X_s) + \sum_{s \in V_X} M(X_{s}) - A (\theta)
\end{align}
where the function $A(\theta)$ is so-called log-partition function, that provides the log-normalization constant for the multivariate distribution.
 
\citet{YRAL13CRF} then developed the conditional extension of \eqref{EqnExpMRF} above. Consider a response random vector $Y := (Y_1,\hdots,Y_p) \in \mathcal{Y}^p$ and a covariate random vector $X := (X_1,\hdots,X_q) \in \mathcal{X}^q$. Suppose that the node conditional distributions of $Y_s$ for all $s \in V_Y$ follow exponential families:
\begin{align*}
	\log \P[Y_s | Y_{V_Y \backslash s}, X ] =  E(Y_{V_Y \backslash s},X)\, B(Y_s) + M(Y_s) - \bar{D}(Y_{V_Y \backslash s},X) ,
\end{align*}
and that the corresponding joint $\P_{\exp}[Y|X]$ factors according to set of cliques $\C_Y$ among random variables in $Y$. Then, \citep{YRAL13CRF} showed that such a joint distribution consistent with the node-conditional distributions does exist, and necessarily has the form:
\begin{align}\label{EqnExpCRF}
		\log \P_{\exp}[Y|X] = \sum_{C \in \C_{Y}} \theta_{C}(X) \prod_{s \in C} B(Y_s) + \sum_{s \in V_Y} M(Y_{s}) - A_{Y|X}\big(\theta(X)\big),
\end{align}
where $\theta_C(X)$ is any function that only depends on the random vector $X$. 

\subsubsection{Mixed MRFs}
\label{SecMixedMRFs}
While \eqref{EqnExpMRF},\eqref{EqnExpCRF} specify multivariate distributions for variables of varied data-types, they are nonetheless specified for the setting where all the variables belong to the same type. Accordingly, there have been some recent extensions~\citep{YBRAL14,Chen2014} of the above for the more general setting of interest in this paper, where each variable belongs to a potentially different type. Their construction was as follows, and can be seen to be an extension not only of the class of exponential family MRFs in \eqref{EqnExpMRF}, but also of the class of conditional Gaussian models in \eqref{EqnCGModel}. Suppose that the node conditional distributions of each variable $X_s$ for $s \in V_X$ now belongs to potentially differing univariate exponential families: 
\begin{align*}
	\log \P[X_s | X_{V_X \backslash s}] = E_s(X_{V_X \backslash s})\, B_s(X_s) + M_s(X_s) - \bar{D}_s(X_{V_X \backslash s}) ,
\end{align*}
while as before, we require the corresponding joint distribution
$\P_{\exp}[X]$ to factor according to the set of cliques $\C_X$ of a
graph $G_X$. In a precursor to this paper~\citet{YBRAL14}, we showed that such a joint distribution consistent with the node-conditional distributions does exist, and moreover necessarily has the form
\begin{align}\label{EqnExpMixedMRF}
	\log \P_{\exp}[X] = \sum_{C \in \C_{X}} \theta_{C} \prod_{s \in C} B_s(X_s) + \sum_{s \in V_X} M_s(X_{s}) - A (\theta)
\end{align}
where the sufficient statistics $B_s(\cdot)$ and base measure $M_s(\cdot)$ in this case can be different across random variables. While these provide multivariate distributions over heterogeneous variables, as we will show in the main section of the paper, this class of distributions sometimes have stringent normalizability restrictions on the set of parameters $\{\theta_C\}$. In this paper, we thus develop a far more extensive generalization that leverages block-directed graphs.

\subsection{Summary \& Organization}

This paper is organized as follows.  In Section~\ref{Sec:ElemChain},
we first introduce the class of Elementary Block Directed Markov
Random Fields (EBDMRFs), which are the simplest subclass of our class
of graphical models. Here, we assume the heterogeneous set of
variables can be grouped into two groups: $X$ and $Y$ (each of which
could have heterogeneous variables in turn). Our class of EBDMRFs are
then specified via a simple application of the chain rule as  $\P[X,Y]
= \Pe[Y | X] \, \Pe[X]$, where $\Pe[X]$ is set to an exponential family
mixed MRF as in \eqref{EqnExpMixedMRF} from \citep{YBRAL14}, and
$\Pe[Y|X]$ is a novel class of what we call exponential family mixed
CRFs, that extend our prior work on exponential family CRFs in
\eqref{EqnExpCRF} from \cite{YRAL13CRF}, and  exponential family mixed
MRFs in \eqref{EqnExpMixedMRF} from \citep{YBRAL14}. This can be seen
as a generalization of the seminal mixed graphical models work of
\citep{lauritzen1989graphical}. We then discuss the properties of this
class of distribution, including conditions or restrictions on the
parameter space under which the distribution is normalizability. As we
show, our formulation of EBDMRFs have substantially weaker normalizability
restrictions when compared with our preliminary work on exponential
family mixed MRFs~\cite{YBRAL14}.  

In Section~\ref{Sec:RecChain}, we then extend this construction further by recursively applying the chain rule respecting a directed acyclic graph over blocks of variables, resulting in a class of graphical models we call Block Directed Markov Random Fields (BDMRFs). The overall underlying graph of this class of graphical models thus has both directed edges between blocks of variables and undirected edges within blocks of variables. Our construction yields a very general and flexible class of mixed graphical models that directly parameterizes dependencies over mixed variables.
We study the problem of parameter estimation and graph structure learning for our class of BDMRF models in Section~\ref{Sec:Learning}, providing statistical guarantees on the recovery of our models even under high-dimensional regimes. Finally, in Section~\ref{Sec:SimData}
and Section~\ref{Sec:RealData}, we demonstrate the flexibility and applicability of our models via both an empirical study of our estimators through a series of simulations, as well as an application to high-throughput cancer genomics data.

\section{Elementary Block Directed Markov Random Fields (EBDMRFs)}\label{Sec:ElemChain}

In this section, we will introduce a simpler subclass of our eventual
class of graphical models which we term Elementary Block Directed
Markov Random Fields (EBDMRFs). Before doing so, we first develop a
key building block: a novel class of conditional distributions.

\subsection{Exponential Family Mixed CRFs}
We now consider the modeling of the conditional distribution of a heterogeneous random response vector $Y := (Y_1,\hdots, Y_p) \in \mathcal{Y}_1 \times \hdots \times \mathcal{Y}_p$, conditioned on a heterogeneous random covariate vector $X := (X_1,\hdots,X_q) \in 
\mathcal{X}_1 \times \hdots \times \mathcal{X}_q$. Suppose that we
have a graph $G_Y = (V_Y,E_Y)$, with nodes $V_Y$ associated with
variables in $Y$. Denote the set of cliques of this graph by $C_Y$,
and the set of \emph{response} neighbors in $V_Y$ for any response
node $s \in V_Y$ by $N_Y(s)$. Suppose further that we also have a set
of nodes $V_X$ associated with the covariate variables in $X$, and
that for any response node $s \in V_{Y}$, we have a set of \emph{covariate
neighbors} in $V_X$ denoted by $N_{YX}(s)$.

%

Suppose that the response variables $Y$ are locally Markov with respect to their specified neighbors, so that 
\bal
\P[Y_s| Y_{V_Y \backslash s},X] = \P[Y_s|Y_{N_{Y}(s)},X_{N_{YX}(s)}].
\label{CondMixedCRFMarkov}
\eal
Moreover, suppose that this conditional distribution $Y_s$ conditioned on the rest of $Y_{V_Y\backslash s}$ and $X$ is given by an arbitrary univariate exponential family:
\begin{align}\label{EqnExpMixedCRF_Cond}
	\log \P[Y_s | Y_{V_Y \backslash s},X] =  E_s\big(Y_{V_Y \backslash s},X\big)\, B_s(Y_s) + M_s(Y_s) - \bar{D}_s\big(Y_{V_Y \backslash s},X\big),
\end{align}
with sufficient statistic $B_{s}(\cdot)$ and base measure $M_{s}(\cdot)$. Note that there is no assumption on the form of functions $E_s\big(Y_{V_Y \backslash s},X\big)$. The following theorem then specifies the algebraic form of the conditional distribution $\Pe[Y|X]$.

\begin{theorem}\label{ThmExpMixedCRF}
Consider a $p$-dimensional random vector $Y = (Y_1,\hdots,Y_p)$ denoting the set of responses, and let $X = (X_1,\hdots,X_q)$ be a $q$-dimensional covariate vector. Then, the node-wise conditional distributions satisfying the Markov condition in \eqref{CondMixedCRFMarkov} as well as the exponential family condition in \eqref{EqnExpMixedCRF_Cond}, are indeed consistent with a graphical model joint distribution, that factors according to $G_Y$, 
and has the form:
\begin{align}
	\log \Pe[Y|X] = \sum_{C \in \C_{Y}} \theta_{C} \big( X_{N_{YX}(C)} \big) \prod_{s \in C} B_{s}(Y_s) + \sum_{s \in V_Y} M_{s}(Y_{s}) - A_{Y|X}\big(\theta(X)\big),
	\label{EqnExpMixedCRF}
\end{align}
where $N_{YX}(C) = {\cap_{s \in C}}N_{YX}(s)$ and $A_{Y|X}\big(\theta(X)\big)$ is the log-normalization function, which is the function on the set of parameters $\{\theta_{C} \big( X_{N_{YX}(C)} \big)\}_{C\in \C_{Y}}$.
\end{theorem}

The proof of the theorem follows along similar lines to that of the
Hammersley Clifford Theorem, and is detailed in
Appendix~\ref{SecProofMixedCRF}. In
Appendix~\ref{Sec:NormalizabilityMixedCRFs}, we also discuss
constraints on the covariate parameter functions
$\{\theta_{C}(\cdot)\}_{C \in \C_Y}$ that ensure the distribution is
normalizable.  Note that to ensure the Global Markov Property holds,
these can be arbitrarily specified as long as 
they are functions solely of the covariate neighborhoods.

We term this class of conditional distributions as exponential family
mixed CRFs. Notice that this framework provides a multivariate density
over the random response vector $Y$, but not a joint density over both
$X$ and $Y$ as we ultimately desire.

\begin{figure}[tb]
	\centering
	\includegraphics[width=0.4\textwidth]{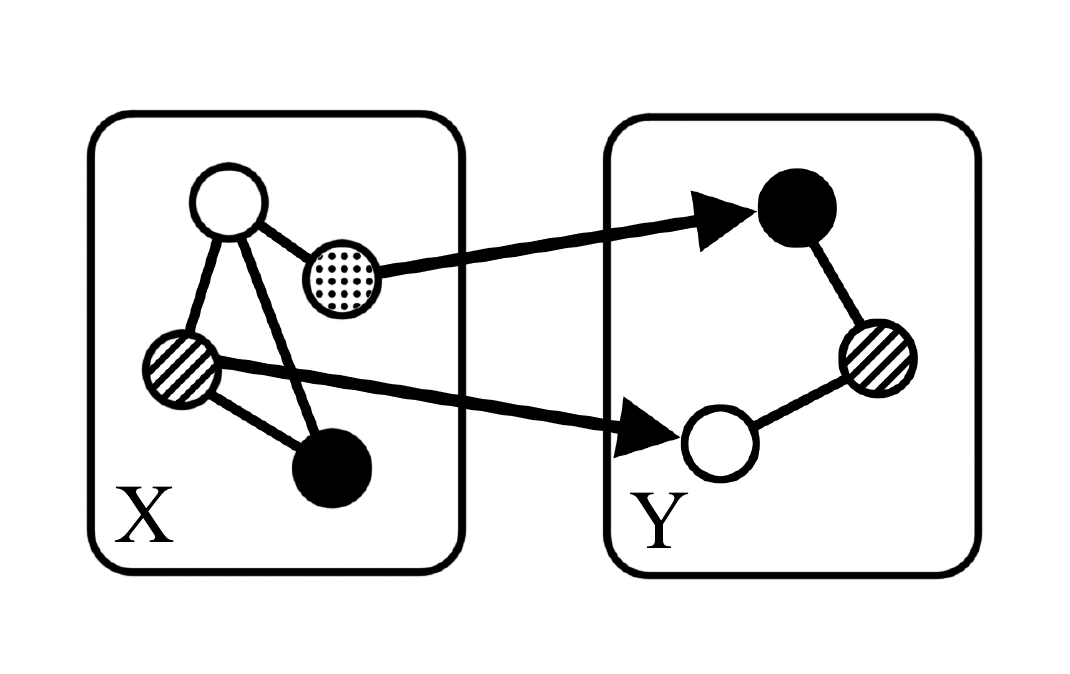}
	\caption{Elementary Block Directed MRF}
	\label{Fig:ElemHet}
\end{figure}

\subsection{Model Specification: EBDMRFs}
We assume that the heterogeneous set of variables can be partitioned into two groups $X := (X_1,\hdots,X_q) \in \mathcal{X}_1 \otimes \hdots \mathcal{X}_q$ and  $Y := (Y_1,\hdots,Y_p) \in 
\mathcal{Y}_1 \otimes \hdots \mathcal{Y}_p$; note that each group could in turn be heterogeneous. Such a delineation of the overall set of variables into two groups is natural in many settings. For instance, the variables in $X$ could be the set of covariates, while the variables in $Y$ could be the set of responses, or $X$ could be cause variables, while $Y$ could be effect variables, and so on. Given the ordering of variables $(X,Y)$, suppose it is of interest to specify dependencies among $Y$ conditioned on $X$, and then of the marginal dependencies among $X$. 

Towards this, suppose that we have an undirected graph $G_Y =
(V_Y,E_Y)$, with nodes $V_Y$ associated with variables in $Y$, and set
of cliques $C_Y$, and in addition, an undirected graph $G_X =
(V_X,E_X)$, with nodes $V_X$ associated with variables in $X$, and set
of cliques $C_X$. Suppose in addition, we have directed edges $E_{XY}$
from nodes in $V_X$ to $V_Y$.  Thus, the overall graph structure has
both undirected edges $E_{X}$ and $E_Y$ among nodes solely in $X$ and
$Y$ respectively, as well as directed edges $E_{XY}$, from nodes in
$X$ to $Y$, as shown in Fig. \ref{Fig:ElemHet}. For any response node
$s \in V_Y$, we will denote the set of \emph{response-specific}
neighbors in $G_Y$ by $N_Y(s)$, and we will again denote \emph{covariate-specific} neighbors in $V_X$, by $N_{YX}(s)$.

Armed with this notation, we propose the following natural joint distribution over $(X,Y)$:
\[ \P[X,Y] = \Pe[Y|X] \, \Pe[X],\]
where the two pieces, $\Pe[Y|X]$ and $\Pe[X]$, are specified as follows. The conditional distribution $\Pe[Y|X]$ is specified by an exponential family \emph{mixed} CRF as in \eqref{EqnExpMixedCRF},
\[
\log \Pe[Y|X] = \sum_{C \in \C_{Y}} \theta_{C}\big( X_{N_{YX}(C)}  \big) \prod_{s \in C} B_{s}(Y_s) + \sum_{s \in V_Y} M_{s}(Y_{s}) - A_{Y|X}\big(\theta(X)\big),
\]
while the marginal distribution $\Pe[X]$ is specified by an exponential family \emph{mixed} MRF as in \eqref{EqnExpMixedMRF}, 
\[
\log \Pe[X] = \sum_{C \in \C_{X}} \theta_{C} \prod_{t \in C} B_{t}(X_t) + \sum_{t \in V_X} M_{t}(X_{t}) - A_{X}(\theta).
\]
Thus, the overall joint distribution, which we call an Elementary Block Directed Markov Random Field (EBDMRF) is given as:
\begin{align}
	\log \P[X,Y] & = \sum_{C \in \C_{Y}} \theta_{C}\big(X_{N_{YX}(C)} \big) \prod_{s \in C} B_{s}(Y_s) + \sum_{s \in V_Y} M_{s}(Y_{s}) + \sum_{C \in \C_{X}} \theta_{C} \prod_{t \in C} B_{t}(X_t) \nonumber\\
 	& \qquad + \sum_{t \in V_X} M_{t}(X_{t}) - A_{Y|X}\big(\theta(X)\big) - A_{X}(\theta).
	\label{EqnElemBDMRF}
\end{align}

We provide additional intuition on our class of EBDMRF distributions
in the next few sections.  First, we compare the conditional
independence assumptions entailed by the mixed 
graph of the EBDMRF with the Global Markov assumptions entailed by an
undirected graph in Section~\ref{SecCompEBDMRFGlobalMarkov}. Then in
Section~\ref{SecCompEBDMRFMixedMRF}, we compare the form of the EBDMRF
class of distributions with that of the mixed MRF distributions
introduced earlier in Section~\ref{SecMixedMRFs}. Next, we analyze
the domain of the parameters of the distribution by considering
parameter restrictions to ensure normalizability in 
Section~\ref{SecEBDMRFNormalizability}. Finally, in
Section~\ref{SecEBDMRFExamples}, we provide several examples of our
EBDMRF distributions, compare these to Mixed MRF distribution
counterparts, and place them in the context of the
larger graphical models literature.

\subsection{Global Markov Structure}
\label{SecCompEBDMRFGlobalMarkov}

The EBDMRF distribution is specified by a mixed graph with both
directed edges from $V_X$ to $V_Y$ and undirected edges within $V_X$
and $V_Y$. A natural question that arises then is 
what are the conditional independence assumptions specified by the
edges in the mixed graph? 
Consider the undirected edges in $G_X$: by construction, these
represent the Markov conditional independence assumptions in the
marginal distribution $\Pe[X]$. The remaining undirected and directed
edges incident on the response nodes $V_Y$ in turn represent Markov
assumptions in the conditional distribution $\Pe[Y|X]$ as described in
\eqref{CondMixedCRFMarkov} in the previous section. It can thus be
seen that the set of conditional independence assumptions entailed by
the mixed graph differ from those of a purely undirected graph obtained
from the mixed graph by dropping the orientations of the edges from
nodes in $V_Y$ to nodes in $V_X$. 

But under what additional restrictions would a EBDMRF entail Markov
independence assumptions with respect to an undirected graph over the
joint set of vertices $V = V_X \cup V_Y$? This was precisely the
question asked in the classical mixed graphical model work of
\citep{lauritzen1989graphical}. Specifically, we are interested in
outlining what restrictions on $\theta(\cdot)$ in the joint
distribution in \eqref{EqnElemBDMRF} would entail global
Markov assumptions with respect to a graph $G = (V,E)$ over all the
vertices $V = V_X \cup V_Y$. In the following theorem, we provide an
answer to this question: 
\begin{theorem}\label{ThmMarkovProp}
Consider an EBDMRF distribution of the form \eqref{EqnElemBDMRF}, with graph structure specified by undirected edges $E_{X}$ among nodes in $V_X$, and $E_{Y}$ among nodes in $V_Y$, as well as directed edges $E_{XY}$ from nodes in $V_X$ to $V_Y$. Then, if this distribution is globally Markov with respect to a graph $G = (V,E)$ with nodes $V = V_X \cup V_Y$ and edges $E \subseteq V \times V$, then 
\begin{itemize}
	\item[(a)] $E_{X},E_{Y},E_{XY} \subseteq E$.
	\item[(b)] For all response cliques $C \in \C_Y$, 
	\begin{align*}
\theta_{C}(X) = \sum_{d \subseteq N_{YX}(C)} \theta_{C,d}(X),
\end{align*}
where $\theta_{C,d}(X)$ depends on $X$ only through the subvector $X_{\{d\} \cup C}$, and that for any subset $d \subseteq [q]$ such that $\{d\} \cup C$ is not complete with respect to $G$,  $\theta_{C,d}(X) = 0$.
\end{itemize}
\end{theorem}
The theorem thus entails that the covariate parameters $\theta$ and
$\theta(\cdot)$ in \eqref{EqnElemBDMRF} factor with respect to the
overall graph $G$. In other words to ensure the Global Markov
structure holds, the covariate parameters can be
arbitrarily specified as long as they are functions solely of the
covariate-neighborhoods.


\subsection{Comparison to Mixed MRFs}
\label{SecCompEBDMRFMixedMRF}

The previous section investigated the implications as far as global
Markov properties for our EBDMRF distribution with both directed and undirected
edges as compared to the 
undirected edges of the mixed MRF distribution. 
Here, we compare and contrast the \emph{factorized form}
of these distributions, in part because the key reason for
the popularity of graphical model distributions is the factored form
of their joint distributions. 

Suppose we have two sets of variable, $X$ and $Y$, and consider the
case of pairwise interactions.  Then, this special case of the mixed
MRF distributions in Section~\ref{SecMixedMRFs} can be written as: 
\begin{align}\label{EqnMixedMRFPairwise}
	\log \Pe[X,Y] = & \sum_{s \in V_Y} \theta_s B_s(Y_s)  + \sum_{t \in V_X} \theta_{t} B_t(X_{t}) + \sum_{(s,s') \in E_Y }\theta_{ss'} \, B_s(Y_s) \, B_{s'}(Y_{s'}) \nonumber\\
	&  + \sum_{(t,t') \in E_X }\theta_{t t'} \, B_t(X_{t}) \, B_{t'}(X_{t'}) + \sum_{(s,t) \in E_{XY}} \theta_{st} \, B_s(Y_s) \,B_t(X_{t})  \nonumber\\
	&  + \sum_{s \in V_Y} M_s(Y_s) + \sum_{t \in V_X} M_t(X_{t})  -A(\theta)   
\end{align}

For this same pairwise case, the EBDMRF class of distributions in \eqref{EqnElemBDMRF} can be written as:
\begin{align}\label{EqnEBDMRFPairwise}
	\log \P[X,Y] = & \sum_{s \in V_{Y}} \theta_{s}(X_{N_{YX}(s)}) B_{s}(Y_s) + \sum_{(s,s') \in E_Y} \theta_{ss'}(X_{N_{YX}(s,s')}) \, B_s(Y_s) \, B_{s'}(Y_{s'}) \nonumber\\
	& +  \sum_{t \in V_X} \theta_{t} B_t(X_{t}) + \sum_{(t,t') \in E_X }\theta_{t t'} \, B_t(X_{t}) \, B_{t'}(X_{t'})  \nonumber\\
 	& + \sum_{s \in V_Y} M_s(Y_s) + \sum_{t \in V_X} M_t(X_{t}) - A_{Y|X}\big(\theta(X)\big) - A_X(\theta).
\end{align}

Noting that covariate functions can be set arbitrarily, it can be seen that the two distributions have almost similar forms when we set the covariate functions $\{\theta_{s}(X_{N_{YX}(s)}), \theta_{ss'}(X_{N_{YX}(s,s')})\}$ as 
\begin{align}\label{EqnLinParamFunc}
	\theta_{s}(X_{N_{YX}(s)}) :=  \theta_s + \sum_{t \in N_{YX}(s)} \theta_{st} \,B_t(X_{t})  , \quad \text{and} \quad 
	\theta_{ss'}(X_{N_{YX}(s,s')}) := \theta_{ss'}.
\end{align}
In this case, the EBDMRF distribution in \eqref{EqnEBDMRFPairwise} and
the mixed MRF distribution in \eqref{EqnMixedMRFPairwise} differ
primarily due to the non-linear term $A_{Y|X}(\theta(X))$. Notice that
in \eqref{EqnMixedMRFPairwise}, $A(\theta)$ is only dependent on parameters and hence the
form \eqref{EqnMixedMRFPairwise} and indeed more generally all mixed MRFs  belong to the
class of multivariate exponential family distributions.  On the other
hand, $A_{Y|X}(\theta(X))$ in \eqref{EqnEBDMRFPairwise} depends on $X$ and hence even
when the covariate functions are simple linear forms as in \eqref{EqnLinParamFunc},
these are not exponential family distributions.  

Similarly, the conditional distribution $\Pe[X|Y]$ of mixed MRF distribution in \eqref{EqnMixedMRFPairwise} can be easily derived as:
\begin{align*}
	\log \P[X|Y] \propto & \sum_{t \in V_X} \Big(\theta_t + \sum_{s \,|\, (s,t) \in E_{XY}}  \theta_{st} \,B_s(Y_{s}) \Big) B_t(X_t)  \nonumber\\
	& + \sum_{(t,t') \in E_X }\theta_{tt'} \, B_t(X_t) \, B_{t'}(X_{t'}) + \sum_{t \in V_X} M_t(X_t) , 
\end{align*}
while the conditional distribution $\P[X|Y]$  for the case of the EBDMRF, when setting the covariate functions $\{\theta_{s}(X_{N_{YX}(s)}), \theta_{ss'}(X_{N_{YX}(s,s')})\}$ as in \eqref{EqnLinParamFunc}, can be written as
\begin{align*}
	& \log \P[X|Y] \propto  \sum_{t \in V_X} \Big( \theta_t + \sum_{s \,|\, (s,t) \in E_{XY}} \theta_{st} \,B_s(Y_{s})  \Big) B_t(X_t) \nonumber\\
	& + \sum_{(t,t') \in E_X }\theta_{tt'} \, B_t(X_t) \, B_{t'}(X_{t'}) + \sum_{t \in V_X} M_t(X_t) - A_{Y|X}(\theta(X)), 
\end{align*}
which can again be seen to differ primarily due to the term $A_{Y|X}(\theta(X))$.

\noindent
It is also instructive to consider the differing forms of the conditional distributions $\P[Y|X]$ in general. For pairwise mixed MRFs, this can be written as

\begin{align}\label{EqnMixedMRFPairwiseCond}
	\log \Pe[Y|X] & \propto \sum_{s \in V_Y} \Big( \theta_s + \sum_{t \,|\, (s,t) \in E_{XY}} \theta_{st} \,B_t(X_{t})  \Big) B_s(Y_s)   \nonumber\\
	& + \sum_{(s,s') \in E_Y }\theta_{ss'} \, B_s(Y_s) \, B_{s'}(Y_{s'}) + \sum_{s \in V_Y} M_s(Y_s),
\end{align}

while that for the pairwise EBDMRF can be written as
\begin{align}\label{EqnEBDMRFPairwiseCond}
	\log \P[Y|X]  & \propto   \sum_{s \in V_{Y}} \theta_{s}(X_{N_{YX}(s)}) B_{s}(Y_s)  \nonumber\\
	& + \sum_{(s,s') \in E_Y} \theta_{ss'}(X_{N_{YX}(s,s')}) \, B_s(Y_s) \, B_{s'}(Y_{s'}) + \sum_{s \in V_Y} M_s(Y_s).
\end{align}
Thus, the two distributions differ only because the covariate
functions $\{\theta_{s}(X_{N_{YX}(s)})\}$,
$\{\theta_{ss'}(X_{N_{YX}(s,s')})\}$ in the EBDMRF distribution can be
set arbitrarily. However, when they are set precisely equal to the
expressions in \eqref{EqnLinParamFunc}, these two distributions can be
seen to identical.  Overall, the pairwise EBDMRF and Mixed MRF
distributions are strikingly similar in the case where the EBDMRF
distribution has linear covariate functions, differing only by the
normalization terms $A_{Y|X}(\theta(X) )$.  This seemingly minor
difference, however, has important consequences in terms of
normalizability discussed next.

\subsection{Normalizability}\label{SecEBDMRFNormalizability}

An important advantage of our EBDMRF class of distributions is that
the parameter restrictions for normalizability of \eqref{EqnElemBDMRF}
can be characterized simply. Recall that the problem of
normalizability refers to the set of restrictions on the parameter
space that is required to ensure the joint density integrates to one;
this entails ensuring that the log-partition functions are finitely
integrable.  

\begin{theorem}\label{ThmEBDMRFNormalizability}
	For any given set of parameters, the joint distribution in \eqref{EqnElemBDMRF} exists and is well-defined, so long as the corresponding marginal MRF $\Pe[X]$, as well as the conditional distribution $\Pe[Y|X]$ are well-defined.
\end{theorem}
Thus the normalizability conditions for the joint distribution in
\eqref{EqnElemBDMRF} reduces to those for the marginal MRF $\Pe[X]$
and for the mixed CRF $\Pe[Y|X]$.

In the previous section, we saw that the pairwise EBDMRF and the
pairwise mixed MRF distributions have very similar form when the
covariate functions $\{\theta_{s}(X_{N_{YX}(s)}),
\theta_{ss'}(X_{N_{YX}(s,s')})\}$ are set as in
\eqref{EqnLinParamFunc}. It will thus be instructive to compare the
\emph{normalizability restrictions} imposed on the mixed MRF
parameters with those on the EBDMRF parameters in this special case. 

\vskip0.1in
\noindent
Let us introduce shorthand notations for the following expressions:
\begin{align*}
	& F(X; \theta) =  \sum_{t \in V_X} \theta_{t} B_t(X_{t}) + \sum_{(t,t') \in E_X }\theta_{t t'} \, B_t(X_{t}) \, B_{t'}(X_{t'}) + \sum_{t \in V_X} M_t(X_{t}) \, , \nonumber\\
	& F(X,Y; \theta) = \sum_{s \in V_Y} \theta_s B_s(Y_s) + \sum_{(s,s') \in E_Y }\theta_{ss'} \, B_s(Y_s) \, B_{s'}(Y_{s'}) \nonumber\\
	& + \sum_{(s,t) \in E_{XY}} \theta_{st} \, B_s(Y_s) \,B_t(X_{t}) + \sum_{s \in V_Y} M_s(Y_s) .
\end{align*} 
Then, the log-partition function $A(\theta)$ in the mixed MRF distribution~\eqref{EqnMixedMRFPairwise} can be written as
\begin{align}\label{EqnMixedMRF-A}
	A(\theta) & := \log \left( \sum_{X,Y} \exp \Big\{ F(X;\theta)  + F(X,Y;\theta) \Big\} \right)  \nonumber\\
	 & =  \log \left( \sum_{X} \Big[ \exp \{ F(X;\theta) \} \sum_Y \exp \{ F(X,Y;\theta)\} \Big]\right) .
\end{align} 

On the other hand, the pairwise EBDMRF distribution in \eqref{EqnEBDMRFPairwise} has the following two ``normalization'' terms:
\begin{align*}
	A_{Y|X}\big(\theta(X)\big) & :=  \log  \left(\sum_Y \exp\{ F(X,Y;\theta)\} \right)\, ,  \nonumber\\
	A_X(\theta) & := \log  \left(\sum_X \exp \{ F(X;\theta) \}\right).
\end{align*} 
Thus, the overall ``normalization'' term in \eqref{EqnEBDMRFPairwise} can be written as:
\begin{align}\label{EqnRecMRF-A}
	A_{Y|X}\big(\theta(X)\big) + A_X(\theta)  = & \log  \left(\sum_Y \exp\{ F(X,Y;\theta)\} \right) + \log  \left(\sum_X \exp \{ F(X;\theta) \}\right) \nonumber\\
	 = & \log \left\{ \left(\sum_X \exp \{ F(X;\theta) \} \right) \left( \sum_Y \exp\{ F(X,Y;\theta)\} \right) \right\}.
\end{align}
In contrast to those of the mixed MRF, these terms are not
normalization \emph{constants} as the term
$A_{Y|X}\big(\theta(X)\big)$ is a function of $X$. Comparing
\eqref{EqnMixedMRF-A} and \eqref{EqnRecMRF-A}, it can be seen that the
two expressions become identical in the special case where
$\theta_{st} = 0$ for all $(s,t) \in E_{XY}$, which would entail that
$X$ and $Y$ are independent, and that $F(X,Y;\theta)$ is only a
function of $Y$. 

But more generally, how would the two different normalization terms affect the normalizability of the two classes of distributions? Would one necessarily be more restrictive than the other? Interestingly, the following theorem indicates that the EBDMRF distribution imposes strictly weaker conditions for normalizability.
\begin{theorem}\label{Prop:MRF-RecMRF}
Consider the pairwise mixed MRF distribution~\eqref{EqnMixedMRFPairwise}. Then, if its log-partition function $A(\theta)$ is finite, then the EBDMRF distribution form in \eqref{EqnEBDMRFPairwise}, with the covariate functions $\{\theta_{s}(X_{N_{YX}(s)}), \theta_{ss'}(X_{N_{YX}(s,s')})\}$ set to the linear forms in \eqref{EqnLinParamFunc}, is normalizable as well.
\end{theorem}  
We provide a proof of this assertion in Appendix~\ref{Proof:Prop:MRF-RecMRF}.
Namely if the log-partition function $A(\theta)$ in \eqref{EqnMixedMRF-A} is
finite, then both $A_{Y|X}(\theta(X))$ and $A_{X}(\theta)$ in \eqref{EqnRecMRF-A}
must also be finite; thus, if pairwise mixed MRFs are normalizable,
then so are EBDMRFs.  The inverse of statement in
Theorem~\ref{Prop:MRF-RecMRF} \emph{does 
  not} hold in general, which can be demonstrated by several
counter-examples discussed in the next section.

\subsection{Examples}
\label{SecEBDMRFExamples}

We provide several examples of EBDMRFs to better illustrate the
properties and implications of our models, and relate them to the
broader graphical models literature.

\paragraph{Gaussian-Ising EBDMRFs}

Gaussian-Ising mixed graphical models have been well studied in the
literature
\citep{lauritzen1989graphical,cheng2013high,lee2012learning,YBRAL14}.
Here, we show that each of these studied classes 
are a special case of our class of EBDMRF models, and which in addition provide other formulations of these mixed models. Suppose we have a set of binary valued random variables $X$, each with domain
$\mathcal{X} \in \{-1, 1\}$, and a set of continuous valued random
variables $Y$, each with domain, $\mathcal{Y} \in \real$.  We can
specify the node-conditional distributions associated with $X$ as
Bernoulli with sufficient statistics and base measure given by
$B_t(X_{t})=X_{t}$, $M_t(X_{t})=0$ for all $t \in V_X$; similarly, we
can specify the node-conditionals associated with $Y$ as Gaussian with
sufficient statistics and base measure given by
$B_s(Y_{s})=\frac{Y_{s}}{\sigma_s}$, and
$M_s(Y_{s})=-\frac{Y_{s}^{2}}{2\sigma^{2}_s}$, respectively for all $s 
\in V_{Y}$. Given these two sets of binary and real valued random
variables, there 
are then three primary ways of specifying joint EBDMRF distributions: the mixed
MRF $\Pe[X,Y]$ from~\citet{YBRAL14,Chen2014} , the EBDMRF specified by $\Pe[Y
| X] \, \Pe[X]$, 
and the EBDMRF specified by $\Pe[X | Y] \, \Pe[Y]$; note that these
various formulations do not coincide and give distinct ways of
modeling a joint density.  When we
consider only pairwise interactions and linear covariate functions as
in \eqref{EqnLinParamFunc}, these models 
all take a similar form, only differing in the log-normalization
terms:
\begin{align}
	\log \P[X,Y] &=  \sum_{s \in V_Y} \frac{\theta_s}{\sigma_s} Y_s
+  \sum_{t \in V_X} \theta_{t} X_{t}
+  \sum_{(s,s') \in E_Y } \frac{\theta_{ss'}}{\sigma_s \sigma_{s'}}
Y_s Y_{s'} \nonumber \\ 
	&  +  \sum_{(t,t') \in E_X } \theta_{tt'} \, X_{t} \, X_{t'}
+ \sum_{(s,t) \in E_{XY}} \frac{\theta_{st}}{\sigma_s} \, Y_{s}
\, X_{t} - \sum_{s \in 
  V_Y} \frac{Y_{s}^{2}}{2\sigma^{2}_s} \nonumber  \\
& - A(\theta)   \label{EqnGauss-Ising-mixedMRF} \\
& 
\ \  \textrm{or} \ \  -  A_{Y|X}\big(\theta(X)\big) -
A_{X}(\theta)    \label{EqnGauss-Ising-EBDMRF} \\
& \ \  \textrm{or} \ \ -  A_{X|Y}\big(\theta(Y)\big) - A_{Y}(\theta) .
  \label{EqnIsing-Gauss-EBDMRF}  
\end{align}
As discussed in Section~\ref{SecCompEBDMRFMixedMRF}, the only
differences between these 
three models are the differing normalization terms, which are determined by the
directionality of the edges between nodes of different types. 
If we define $\Theta$ as the matrix, $[\Theta]_{ss'} = \left\{ \begin{array}{l l}
                        -\frac{1}{\sigma_s^2} & \quad \text{if $s=s'$}\\[0.35em]
                        \frac{ \theta_{ss'}}{ \sigma_s \sigma_{s'}} & \quad \text{otherwise.}
                    \end{array} \right.$,  then as discussed in
\citet{YBRAL14} the
Gaussian-Ising mixed MRF is normalizable when $\Theta \prec 0$.  Then
following from Theorem~\ref{Prop:MRF-RecMRF}, both forms of the Gaussian-Ising 
EBDMRF are also 
normalizable when $\Theta \prec 0$.  Upon inspection, we can see that this
restriction cannot be further relaxed.  We further examine these three formulations
of pairwise MRF models for a set of binary and continuous random
variables through numerical examples in
Section~\ref{Sec:SimData}.

Notice that the form of \eqref{EqnGauss-Ising-mixedMRF} is precisely that of a
special case 
of the conditional Gaussian (CG) models first proposed by
\citet{lauritzen1989graphical} and 
reviewed previously.  
This class of models can also be extended to the case of higher-order
interactions.  For example, the
higher-order interactions in
\citet{cheng2013high} for the CG 
model are another special case of mixed MRF models.  Also, this
formulation of mixed MRF model can easily be extended further to 
consider categorical random variables as 
\citet{lauritzen1989graphical} and more recently 
 \citet{lee2012learning} considered for the CG models.
 
Besides the class of mixed MRFs, our EBDMRF models in \eqref{EqnGauss-Ising-EBDMRF} and \eqref{EqnIsing-Gauss-EBDMRF}  are also normalizable joint distributions, and have many potential applications. Returning to our genomics motivating example,
Gaussian-Ising EBDMRFs may be particularly useful for joint network
modeling of binary mutation variables (SNPs) and continuous gene
expression variables (microarrays).  Biologically, SNPs are fixed
point mutations that influence the dynamic and tissue specific gene
expression.  Thus, we can take the directionality of our EBDMRFs
following from known biological processes: $\Pe[Y | X]
\Pe[X]$ where binary random variables associated with SNPs
influence and hence form directed edges with continuous random variables
associated with genes.

\paragraph{Gaussian-Poisson EBDMRFs}

Existing classes of mixed MRF distributions do not permit dependencies between types of variables for Gaussian-Poisson graphical models
\citep{YBRAL14,Chen2014}, which is a major limitation.  Here, we show that in contrast our EBDMRF formulation permits a full
dependence structure with Gaussian-Poison graphical models, due to its weaker 
normalizability conditions. Consider a set of count-valued random variables $X$, each with domain
$\mathcal{X}=\{0,1,2, \hdots\}$, and a set of continuous real-valued
random variables $Y$, each with domain $\mathcal{Y}=\real$, with
corresponding node-conditional distributions specified by the Poisson
and Gaussian distributions respectively.  

Let us consider the simple case of pairwise models with linear
covariate functions from \eqref{EqnLinParamFunc}.  Then to specify
normalizable EBDMRF distributions,
Theorem~\ref{ThmEBDMRFNormalizability} says that we need 
only show that the corresponding CRF and MRF distributions are
normalizable.  First, consider the EBDMRF defined by $\Pe[ Y | X] \, \Pe[X]$ and consider the conditional Gaussian CRF
distribution, $\Pe[Y|X]$, defined as 
\begin{align}\label{EqnCondGaussianPoisson}
	\log \Pe[Y|X] \propto & \sum_{s \in V_Y} \Big( \theta_s + \sum_{t \text{ s.t } (s,t) \in E_{XY}} \theta_{st} \,X_{t} \Big) \frac{Y_s}{\sigma_s} + \sum_{(s,s') \in E_Y } \frac{\theta_{ss'}}{\sigma_s \sigma_{s'}} \, Y_s \, Y_{s'} - \sum_{s \in V_Y} \frac{Y_{s}^{2}}{2\sigma^{2}_s}. \, 
\end{align}
This conditional Gaussian distribution is well defined for any
value of $X$ as long as $\Theta \prec 0$ where $\Theta$ is the matrix
defined in the previous section.  Hence, as long as $\theta_{st} \leq
0$ for all $s,t \in V_{X}$ and the Poisson MRF is normalizable
\citep{YRAL13PGM}, then the EBDMRF distribution given by $\Pe[X,Y] = \Pe [Y|X]\,
\Pe[X]$ is normalizable and has the following form:
\begin{align}
	 \log \P[X,Y] & = \sum_{s \in V_Y} \Big( \theta_s + \sum_{t
           \text{ s.t } (s,t) \in E_{XY}} \theta_{st} \,X_{t} \Big)
         \frac{Y_s}{\sigma_s} + \sum_{(s,s') \in E_Y } \theta_{ss'} \,
         \frac{Y_s \, Y_{s'}}{\sigma_s \sigma_{s'}}  - \sum_{s \in
           V_Y} \frac{Y_{s}^{2}}{2\sigma^{2}_s} \nonumber\\ 
	& + \sum_{X \in V_X} \theta_t X_t + \sum_{(t,t') \in E_X}
         \theta_{tt'} \, X_t \, X_{t'} - \sum_{t \in V_X} \log(X_t!) -
         A_{Y|X}\big(\theta(X)\big) - A_X(\theta)\, 
	\label{EqnGaussianPoissonBDMRF}
\end{align}
where $A_{Y|X}\big(\theta(X)\big)$ is the log normalization constant of conditional Gaussian CRF and $A_X(\theta)$ is that of Poisson MRF.
 
Similarly, consider the EBDMRF over count and continuous valued
variables given by $\Pe[ X | Y] \, \Pe[Y]$.  Again, we have that the
Gaussian MRF is normalizable if $\Theta \prec 0$; the Poisson CRF is
also normalizable if $\theta_{st} \leq 0$ for all $s,t \in V_{X}$ as
discussed in \citet{YRAL13CRF}.  
Thus, {\em both} forms of our Gaussian-Poison EBDMRF permit non-trivial
dependencies between count-valued and continuous variables.

This interesting consequence should be contrasted with that for the mixed MRF~\eqref{EqnMixedMRFPairwise}, which \emph{does not} allow for interaction terms between the count-valued and continuous variables,
since otherwise the log partition function \eqref{EqnMixedMRF-A}
cannot be bounded \cite{YBRAL14,Chen2014}. In other words, the only way for a Gaussian-Poisson
mixed MRF distribution to exist would be a product of independent
distributions over the count-valued random vector and the continuous-valued vector.  Thus, our EBDMRF construction has important
implications permitting non-trivial dependencies in certain classes of mixed graphical models that were previously unachievable.

\paragraph{Other Examples of pairwise EBDMRF models}

As our EBDMRF framework yields a flexible class of models, there are
many possible other forms that these can take.  Here, we outline
classes of normalizable homogeneous pairwise EBDMRFs for easy
reference.  Note that we state normalizability conditions for these
classes of models without proof as these can easily be derived
from Theorem~\ref{ThmEBDMRFNormalizability} and the conditions
outlined in \cite{YBRAL14}.   
\begin{mylist}
\item {\em Poisson-Ising EBDMRFs}. If we let $X$ be a count-valued
  random vector and $Y$ a binary random vector, then we can specify
  the appropriate node-conditional distributions of $X$ as Poisson and
  of $Y$ as Bernoulli.  This gives us three ways of modeling
  Poisson-Ising mixed graphical models: $\Pe[X, Y]$ via mixed MRFs
  and  $\Pe[Y | X]\, \Pe[X]$ or $\Pe[X | Y]\,
  \Pe[Y]$ via our EBDMRFs.  All three of these mixed graphical
  model distributions are normalizable only if $\theta_{tt'} \leq 0$
  for all $t,t' \in V_{X}$.
\item {\em Exponential-Ising EBDMRFs}. Similar to the above classes of
  models, now let $X$ be a positive real-valued random vector and
  specify its node-conditional distributions via the exponential
  distribution, and let $Y$ be a binary random vector with Bernoulli node-conditional distributions as before. Then, the normalizability conditions for the construction of $\Pe[Y | X]\, \Pe[X]$ will be simply $\theta_{t} < 0$ and $\theta_{tt'} \leq 0$
  for all $t,t' \in V_{X}$. The constructions of $\Pe[X | Y]\, \Pe[Y]$ and mixed MRF $\Pe[X,Y]$ require the additional condition that $\theta_{st} \leq 0$
  for all $(s,t) \in E_{XY}$.     
\item {\em Gaussian-Exponential EBDMRFs}.  Let $X$ be a positive real-valued
  random vector and $Y$ a real-valued random vector; then we can specify
  the appropriate node-conditional distributions of $X$ as exponential and
  of $Y$ as Gaussian.  Similar to the Gaussian-Poisson case, the mixed
  MRF distribution {\em does not} permit dependencies between $X$ and $Y$
  \citep{YBRAL14,Chen2014}.  
  But again, our EBDMRF
  distribution $\Pe[Y | X] \, \Pe[X]$ exists,
  permits non-trivial dependencies between nodes in $X$ and $Y$, and is
  normalizable under very similar conditions as the Gaussian-Poisson
  EBDMRF case.  
  
\item {\em Exponential-Poisson EBDMRFs}. Interestingly, we can specify
  all three classes of mixed graphical model distributions for $X$, a
  count-valued random vector with node-conditionals specified as
  Poisson, and for $Y$, a positive real-valued random vector with
  node-conditionals specified as exponential. Here, the normalizability conditions for the construction of $\Pe[X | Y]\, \Pe[Y]$ will be $\theta_{tt'} \leq 0$ for all $t,t' \in V_{X}$ and $\theta_{s} < 0$, $\theta_{ss'} \leq  0$ for all $s,s' \in V_{Y}$. 
  The constructions of $\Pe[Y | X]\, \Pe[X]$ and mixed MRF $\Pe[X,Y]$ additionally require the condition that $\theta_{st} \leq 0$ for all $(s,t) \in E_{XY}$.

\end{mylist}

Note also, that other univariate exponential family distributions can
be used to specify these homogeneous pairwise EBDMRFs.  A particularly
interesting class of these could be the variants of the Poisson
distribution proposed by \citet{YRAL13PGM} to build Poisson graphical models
that permit both positive and negative conditional dependencies.
Within our EBDMRFs, these could be used to expand the possible
formulations of mixed Poisson graphical models that are not restricted
to negative conditional dependence relationships.

Additionally, we have only studied homogeneous pairwise models, but
heterogeneous pairwise EBDMRFs may be of interest in many
applications.  For example, suppose we have count-valued nodes, $X$,
associated with Poisson node-conditional distributions,
and let $Y = (Y_{1}, Y_{2})$ be a set of mixed nodes with $Y_{1}$
binary-valued associated with Bernoulli node-conditionals and $Y_{2}$
continuous associated with Gaussian node-conditionals.  Then, we could use 
heterogeneous EBDMRFs to model dependencies between these three sets
of variables resulting in a form of Gaussian-Ising-Poisson mixed
graphical model: $\P[Y_1,Y_2,X] = \Pe[Y_1,Y_2 | X ] \, \Pe[X]$ where
$\Pe[Y_1,Y_2 | X ]$ is the mixed Gaussian-Ising CRF as defined in
\eqref{EqnExpMixedCRF}.  This distribution is
normalizable if $\Theta \prec 0$ and $\theta_{tt'} \leq 0$ for all $t,t' \in V_{X}$, which follows from the above discussion.

Finally, the examples we have highlighted were focused on pairwise
MRFs with linear covariate functions.  Our class of models, however,
is more flexible and permits non-linear relationships in the covariate
functions $\theta$.  
Note that these non-linear covariate functions
are not permitted in the mixed MRF construction, giving our EBDMRFs yet another important advantage for flexibly modeling mixed multivariate
data.

\section{Block Directed Markov Random Fields (BDMRFs)}\label{Sec:RecChain}

The class of EBDMRF distributions in the earlier section was specified
in terms of a marginal MRF and a conditional CRF distribution
given a binary partition of the set of all variables. We will now
consider a generalization of this construction. 

\subsection{Model Specification.}

\vskip0.1in
\noindent
Let $X$ denote the set of random variables which will form the nodes
or vertices, $V$, of our network model. Suppose that $V$ can be
partitioned 
into an \emph{ordered} set of disjoint exhaustive sets
$V_{1},\hdots,V_{m}$, so that $V_{i} \cap V_{j} = \emptyset, \,
\forall i \neq j$ and $\cup_{j=1}^{m}V_{j} = V$. Consider a mixed
graph $G = (V,E)$ with both directed and undirected edges. Suppose
that the undirected edges are purely between vertices within a single
subset, so that for any undirected edge $(s,t) \in E$, we have that
$s,t \in V_i$, for some $i \in [m]$. Suppose further that any directed
edge points from a vertex in a set with a lower index to that in a set
with a higher index, so that for any directed edge $(s,t) \in E$, we
have that $s \in V_i, t \in V_j$, with $i < j$. 
In the existing literature, note that such a mixed graph
with ordered blocks, directed edges between nodes in different blocks,
and undirected edges within 
blocks, has been referred to as a block-recursive or chain-graph
\citep{Lauritzen}.

This mixed graph induces a directed acyclic graph (DAG) over
the subsets $\{V_{i}\}_{i=1}^{m}$, since there can be no directed
cycles due to the ordering constraint on the directed
edges. Conversely, any DAG over the subsets $\{V_{i}\}_{i=1}^{m}$ in
turn induces a partial ordering over these subsets which can be used
as a ordering constraint on the directed edges in the mixed graph. 
Our mixed graph is thus specified by a partial ordering of the
blocks, and correspondingly, the vertices. The previous 
section on elementary chain models over random vectors $(X,Y)$ can be
understood as using an elementary mixed graph with the partition
$V_{1} = V_{X}$ and $V_{2} = V_{Y}$ of the set of vertices $V = V_{X}
\cup V_{Y}$, and with the directed edges only leading from nodes in
$X$ to nodes in $Y$. 

Now, we seek to define the general class of mixed graphical models
associated with this blocked construction tat we term Block Directed
Markov Random Fields (BDMRFs). First, however, we need to set up some
graph-theoretic notation. For any $i \in [m]$, 
indexing the $m$ subsets $\{V_{i}\}_{i=1}^{m}$, we define the set of
``parent'' subsets 
\[ \text{PA}(i) = \cup_{j=1}^{m} \big\{V_{j}: \exists \, \text{directed} (s,t) \in E, \, s \in V_{j}, t \in V_{i}\big\}.\] 

We will also overload notation, and use $\text{PA}(t)$ to denote the set of parent \emph{nodes} of any node $t \in V$:
\[ \text{PA}(t) = \big\{s: \exists \, \text{directed} (s,t) \in E \big\}.\] 
For any subset of nodes $C \subset V$, we also let $\text{PA}(C) = {\cap_{s \in C}} \text{PA}(s)$. For $i \in [m]$, let $E_{i} := E \cap (V_i \otimes V_i)$ 
denote the set of undirected edges within block $V_i$, and let $\C_{V_i}$ denote the set of cliques (fully connected subgraphs) with respect to just the edges in $E_i$.

Armed with all this graph-theoretic notation, we can then define the following general class of Block Directed Markov Random Fields (BDMRFs):
\[ \P[X] = \prod_{i=1}^{m} \P_{\exp}[\XV|\XPa], \]
where $\Pe[\XV|\XPa]$ is specified by an exponential family \emph{mixed} CRF of the form detailed in \eqref{EqnExpMixedCRF}:
\begin{align}
	\log \P_{\exp}[\XV|\XPa] = \sum_{C \in \C_{V_i}} \theta_{C}\left(X_{\text{PA}(C)}\right) \prod_{s \in C} B_{s}(X_s) + \sum_{s \in V_i} M_{s}(X_{s}) - A_{i}\left(\theta\left(\XPa\right)\right),
\label{EqnBlockMixedCRF}
\end{align}

Substituting in these expressions for exponential mixed CRFs in the overall joint distribution of BDMRFs, we arrive at the following form:
\begin{align}
\log \P[X] = \sum_{i=1}^{m} \left\{ \sum_{C \in \CV} \theta_{C}\left(X_{\text{PA}(C)}\right) \prod_{s \in C} B_{s}(X_s)   + \sum_{s \in V_i} M_{s}(X_{s}) - A_{i}\left(\theta\left(\XPa\right)\right) \right\}.
\label{EqnBDMRF}
\end{align}

Because of the recursive conditioning which results in directed edges
between nodes in different blocks, our BDMRFs are distinct from the
mixed MRFs of \citet{YBRAL14} which consist of only undirected edges.  Note
also that like the EBDMRF counterparts and unlike mixed MRFs, our
BDMRFs do not necessarily correspond to an exponential family
distribution because of the log-normalization constants (see
Section~\ref{SecCompEBDMRFMixedMRF}).  Arguably,
these could be considered as the teleological 
endpoint of the CG chain models first proposed by
\citet{lauritzen1989graphical}.

\subsection{Global Markov Structure \& Normalizability.}
As in the classical CG chain model as well as our EBDMRFs discussed
previously, we can derive the form the covariate parameters,
$\theta(\XPa)$, must take to ensure global Markov assumptions hold.
Recall the notation that for any $i \in 
[m]$, $E_{i} := E \cap (V_i \otimes V_i)$ denotes the set of
undirected edges within block $V_i$, and $\C_{V_i}$ denotes the set of
cliques with respect to just the edges in $E_i$. Then it can be seen
that the BDMRF distribution in \eqref{EqnBDMRF} is specified by
covariate functions $\{\theta_{C}(X)\}$ for $C \in \C_{V_i}, i \in
[m]$. The question then remains: Under what additional restrictions on
these parameters 
$\{\theta_{C}(X)\}$ would the BDMRF joint distribution in
\eqref{EqnBDMRF} entail Markov independence assumptions with respect
to an \emph{undirected graph}  $G = (V,E)$ over all the vertices $V$?
In the following theorem, we provide an answer to this question: 
\begin{theorem}\label{ThmRecurMarkovProp}
Consider a BDMRF distribution of the form \eqref{EqnBDMRF}, with an
underlying mixed graph $G = (V,E)$. Now, suppose this distribution is
globally Markov with respect to a graph $H = (V,E')$, with undirected
edges $E' \subseteq V \times V$. Then,  
\begin{itemize}
\item[(a)] 
The undirected skeleton, $\text{skeleton}(E)$, of the mixed graph $G =
(V,E)$, consisting of all its edges sans directions, satisfies:
$\text{skeleton}(E) \subset E'$. 
\item[(b)]
For all $i \in [m]$, and within-block cliques $C \in \C_{V_i}$, 
\begin{align*}
\theta_{C}(X_{\text{PA}(C)}) = \sum_{d \subseteq \text{PA}(C)} \theta_{C,d}(X),
\end{align*}
where $\theta_{C,d}(X)$ depends on $X$ only through the subvector
$X_{\{d\} \cup C}$, and that for any subset $d \subseteq [p]$ such
that $\{d\} \cup C$ is not complete with respect to $H$,
$\theta_{C,d}(X) = 0$. 
\end{itemize}
\end{theorem}
The theorem thus entails that the covariate parameters $\theta(\cdot)$
in \eqref{EqnBDMRF} \emph{factor} with respect to the overall
undirected graph $H$.  Thus as with our EBDMRFs, to ensure a graph
consistent with the global Markov structure, covariate parameters
can be arbitrarily specified as long as they are functions solely of
the node-neighborhoods.

Additionally, we can easily generalize the conditions for
normalizability  discussed in Section~\ref{Sec:ElemChain} from our
EBDMRFs to our new BDMRF distributions:
\begin{theorem}\label{ThmBDMRFNormalizability}
For any given set of parameters, the joint distribution in
\eqref{EqnBDMRF} exists and is well-defined, so long as the
corresponding conditional mixed CRFs $\P_{\exp}[\XV|\XPa]$ in
\eqref{EqnBlockMixedCRF} are well-defined. 
\end{theorem}
Hence, the distributional form in \eqref{EqnBDMRF} is normalizable and
valid as long as each individual mixed CRF $\P_{\exp}[\XV|\XPa]$
\eqref{EqnBlockMixedCRF} is normalizable. The conditions under which
individual mixed CRF is normalizable are detailed in
Appendix~\ref{Sec:NormalizabilityMixedCRFs} in Propositions
\ref{Prop_Finite} and \ref{Prop_Infinite}.  As with our EBDMRFs, these
BDMRF normalizability conditions are weaker than those for mixed MRFs, as
can be easily seen in the analogous extension of
Theorem~\ref{Prop:MRF-RecMRF}.

\subsection{Examples.}
\label{Sec:BDMRFExamples}

Our class of BDMRF distributions are a very flexible class that allows
one to specify mixed graphical models in numerous ways.  We illustrate
this flexibility through a simple example: Gaussian-Ising-Poisson
graphical models for mixed multivariate distributions over continuous,
count, and binary-valued random variables.

\paragraph{Gaussian-Ising-Poisson Mixed Graphs.}

Consider a partition of variables into three homogeneous blocks, $X,
Y, Z$, such that that variables $\{X_t\}$ in block $X$ are binary,
with domain $\mathcal{X}_t = \{-1,1\}$, variables $\{Y_s\}$ in block $Y$ are
continuous, with domain $\mathcal{Y}_s = \mathbb{R}$, and variables
$\{Z_u\}$ in block $Z$ are count-valued, with domain $\mathcal{Z}_u =
\{0,1,\hdots\}$.  We can specify the corresponding node-conditional
distributions as follows: Bernoulli for block $X$ with 
sufficient statistic and base measure $B_t(X_{t})=X_{t}$ and
$M_t(X_{t})=0$ respectively; Gaussian with known variance $\sigma^2$
for block $Y$ with
$B_s(Y_{s})=\frac{Y_{s}}{\sigma_s}$ and $M_s(Y_{s})=-
\frac{Y_s^2}{2\sigma_s^2}$; and finally, Poisson for block $Z$ with
$B_u(Z_{u})=Z_{u}$ and  $M_u(Z_u)=-\log (Z_u!)$. 

\begin{figure}[b]
	\centering
	\includegraphics[width=0.4\textwidth]{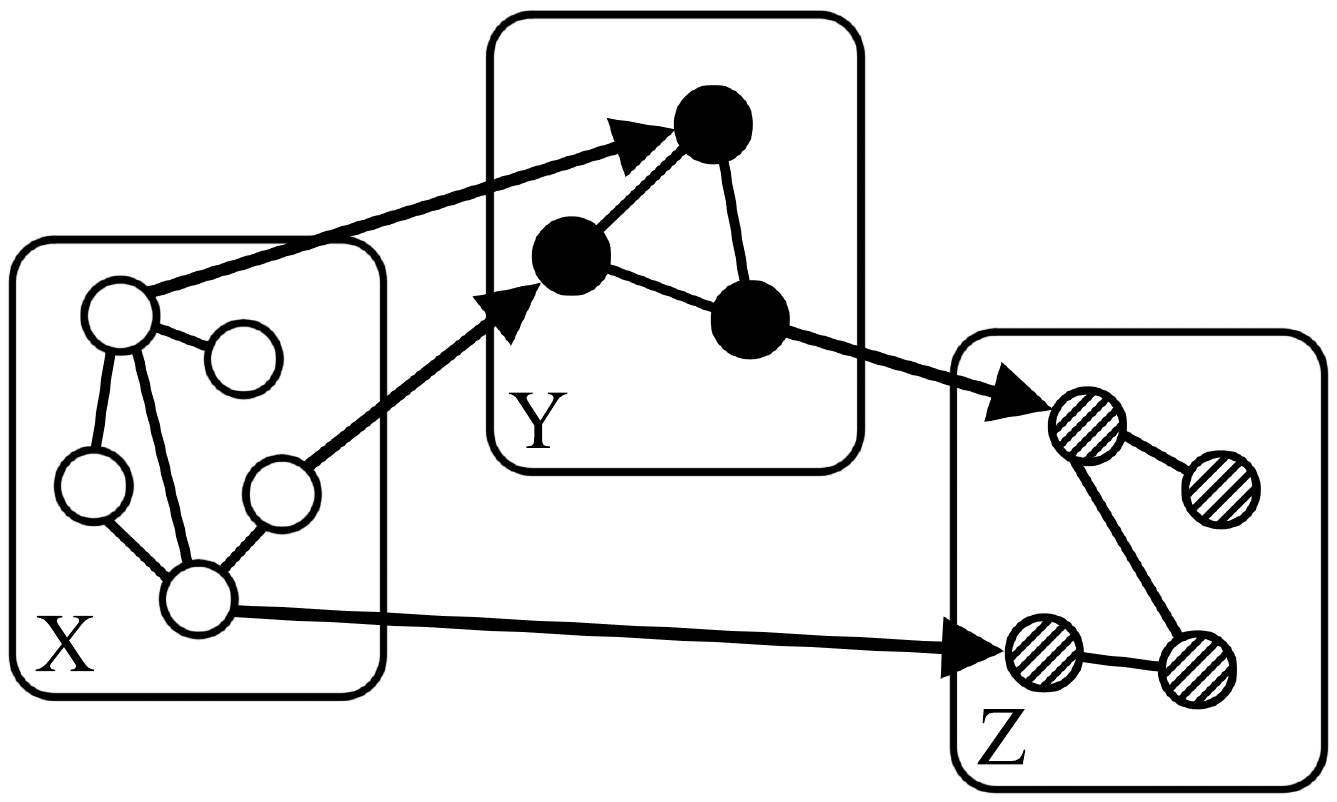}
	\caption{Block Directed MRFs: The joint distribution of $(X,Y,Z)$ can be written as $\Pe[X]\Pe[Y|X]\Pe[Z|X,Y]$}
	\label{Fig:ElemHom3way}
\end{figure}

Then, one possible way of specifying a joint distribution, $\P(X,Y,Z)$
over these variables is given by the following BDMRF: 
$\Pe[X,Y,Z] = \Pe[X] \allowbreak \Pe[Y|X] \, \Pe[Z|X,Y]$.  For pairwise graphical
models, this distribution is denoted by directed edges extending from
nodes in $X$ to nodes in $Z$ and $Y$
and from nodes in $Y$ to nodes in $Z$ as illustrated in Figure
\ref{Fig:ElemHom3way}.   We can write the form of this joint
BDMRF distribution as follows: 
\begin{align*}
	& \log \P[X,Y,Z] = \log \Pe[X] + \log \Pe[Y|X] + \log \Pe[Z|X,Y] 
\end{align*}
where $\Pe[X]$ is an Ising MRF, $\Pe[Y|X]$ is a 
Gaussian CRF as described earlier, and $\Pe[Z|X,Y]$ follows a Poisson 
CRF.  When we specify linear covariate functions as in \eqref{EqnLinParamFunc}, these
take the following forms: 
\begin{align*}
	& \log \Pe[X] = \sum_{t \in V_X} \theta_t \, X_t + \sum_{(t,t') \in E_X} \theta_{tt'} \, X_t \, X_{t'} -A_X(\theta) \, , \nonumber\\
	& \log \Pe[Y|X] = \sum_{s \in V_Y} \Big( \theta_s + \sum_{t \text{ s.t } (s,t) \in E_{XY}} \theta_{st} \,X_{t} \Big) \frac{Y_s}{\sigma_s} + \sum_{(s,s') \in E_Y } \frac{\theta_{ss'}}{\sigma_s \sigma_{s'}} \, Y_s \, Y_{s'} \nonumber\\
	& \hspace{2cm} - \sum_{s \in V_Y} \frac{Y_{s}^{2}}{2\sigma^{2}_s} - A_{Y|X}\big(\theta(X)\big) \, ,\nonumber\\ 
	& \log \Pe[Z|X,Y] = \sum_{u \in V_Z} \Big( \theta_u + \sum_{t \text{ s.t } (t,u) \in E_{XZ}} \theta_{tu} \,X_{t} + \sum_{s \text{ s.t } (s,u) \in E_{YZ}} \theta_{su} \,\frac{Y_{s}}{\sigma_s} \Big) Z_u  \nonumber\\
	& \hspace{2cm} + \sum_{(u,u') \in E_Z } \theta_{uu'} \, Z_u \, Z_{u'} - \sum_{u \in V_Z} \log (Z_u!) - A_{Z|X,Y}\big(\theta(X,Y)\big) \, .  
\end{align*}
Notice that our BDMRF model is normalizable if each of the Ising MRF,
Gaussian CRF and Poisson CRFs are normalizable.  We discussed these
conditions in Section~\ref{SecEBDMRFExamples}; namely, $\Theta \prec
0$ for the Gaussian CRF as defined previously and 
$\theta_{uu'}  \leq 0$ for all $u,u' \in V_Z$ for the Poisson CRF.
This particular BDMRF formulation then exists and permits non-trivial
dependencies between count-valued, continuous and binary-valued random
vectors.

Even in this simple case of three homogeneous blocks of variables, we
have given just one possible formulation of a BDMRF model and there exists a
combinatorial number of ways to specify a valid BDMRF 
model in this case as the Ising CRF, Gaussian CRF, and Poisson CRF
distributions exist under fairly mild conditions.  Thus, even for
fixed variable types, our BDMRFs allow for a most flexibly and rich
class of dependence structures between data of mixed types.
Furthermore, the recursive conditioning via our mixed CRFs
substantially weakens the normalizability conditions so that
non-trivial dependencies are permitted in many more cases such as our
above example where the mixed MRF distribution does not permit
dependencies between variables of different types.

\paragraph{Applicability of BDMRFs}

To specify our BDMRF model, one must know both the block partitions of
variables as well as the ordering of these blocks which determine the
directionality of edges between nodes in different blocks.  While
these assumptions may seem restrictive, they naturally arise in many
settings.  For example, sets of variables measured over time naturally
form blocks corresponding to each time point and the blocks can
clearly be ordered according to a chain graph.  Analogously, this
natural ordering arises with applications to sequential treatments and
natural language processing among many others.

While it may be less obvious, natural directionality
between blocks of variables can also arise in settings that are not
associated with time or sequences.  Consider the example of
different types of high-throughput genomics data measured on the same
set of samples as previously discussed in
Section~\ref{Sec:Intro}. Here, we have
mutations and aberrations (via SNPs and copy number 
aberrations) which are binary, gene expression (via RNA-sequencing)
which are counts,
and epigenetic markers (via methylation arrays) which are continuous.
Thus, each set of genetic variables forms a natural block of binary
sequence genetic markers, $X$, count-valued functional genetic
markers, $Y$, and continuous epigenetic markers, $Z$.  Biologically,
however, scientists have established how these genomic variables
interact with each other.  For example, mutations and
aberrations are fixed alterations in the DNA sequence that influence
gene expression.  An epigenetic marker or methylated region determines
whether a sequence of 
DNA is ready for transcription, and thus also influences gene
expression.  Then based on scientific knowledge, we expect both
mutation and methylation markers to point to gene expression
markers, thus determining the directionality between these blocks of
variables.  Within these types of genomic
biomarkers, however, we expect Markov conditional dependencies
corresponding to undirected edges, which precisely correspond to the
mixed graph construction underlying our BDMRFs.

\section{Learning BDMRFs}\label{Sec:Learning}

\begin{algorithm}[tb]
	\caption{Graph Selection Algorithm for BDMRFs}
	\label{alg:blockrec}
	\begin{algorithmic}
		\STATE {\bfseries Input:} $n$ i.i.d. samples
                $\big\{X^{(j)}\big\}_{j=1}^n$ and the set of partial
                DAG orderings $\{V_1, V_2, \hdots, V_m\}$.
		\FOR{$i=1$ {\bfseries to} $m$}
		\STATE Learn conditional mixed CRF
                $\P_{\exp}[\XV|\XPa]$ via neighborhood selection.
		\ENDFOR
		\STATE {\bfseries Output:} $\P[X] = \prod_{i=1}^{m} \P_{\exp}[\XV|\XPa]$.
	\end{algorithmic}
\end{algorithm}

\begin{figure}[tb]
	\centering
	\subfigure[]{
		\hspace{-0.5cm}\includegraphics[width=0.32\textwidth]{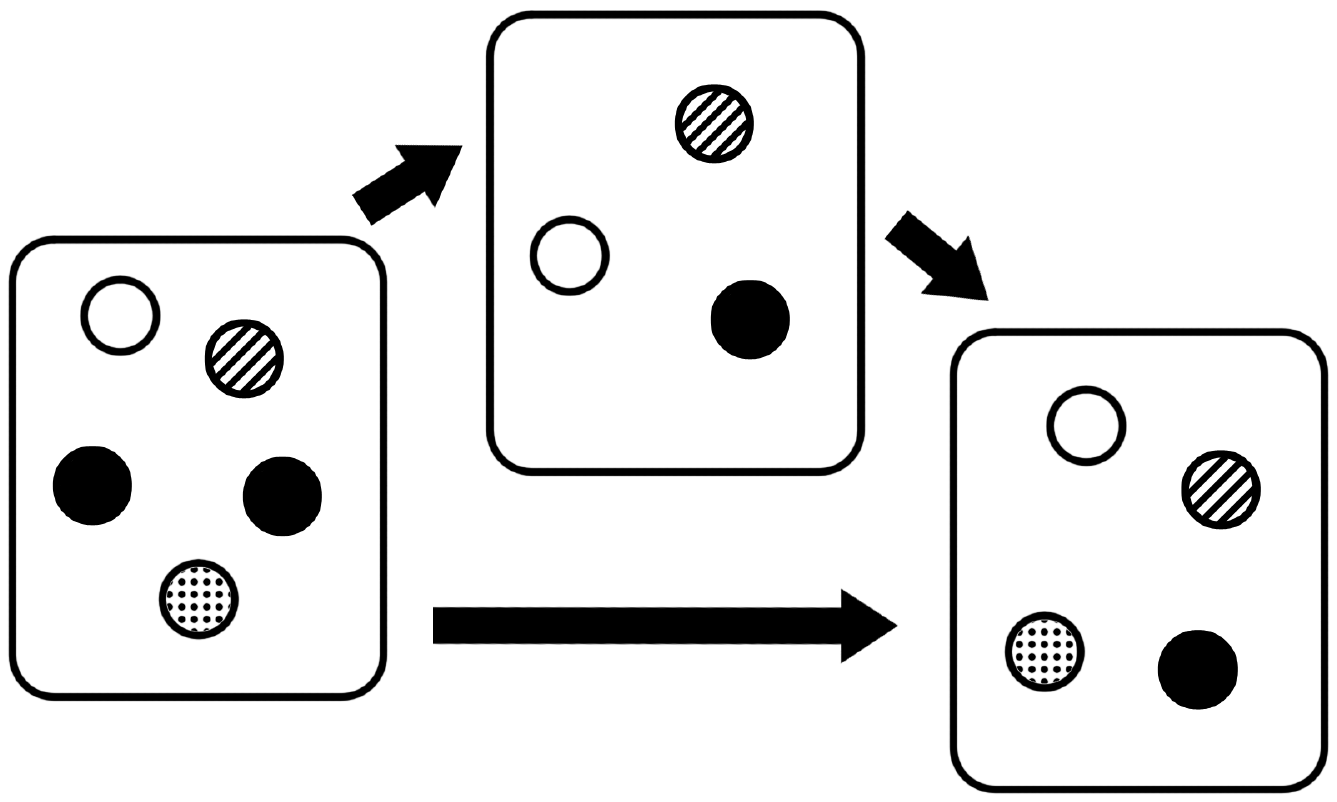}
	}
	\subfigure[]{
		\includegraphics[width=0.32\textwidth]{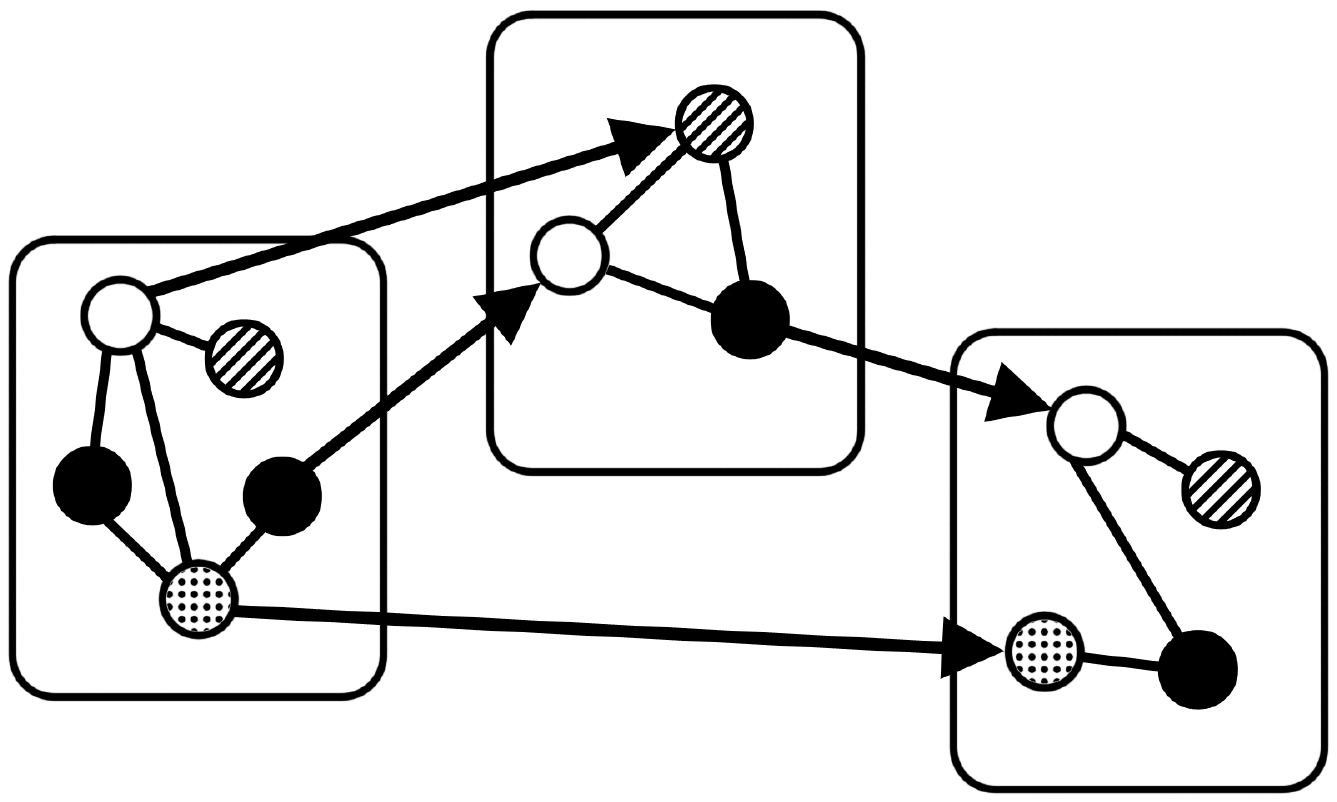}
	}
	\subfigure[]{
		\includegraphics[width=0.32\textwidth]{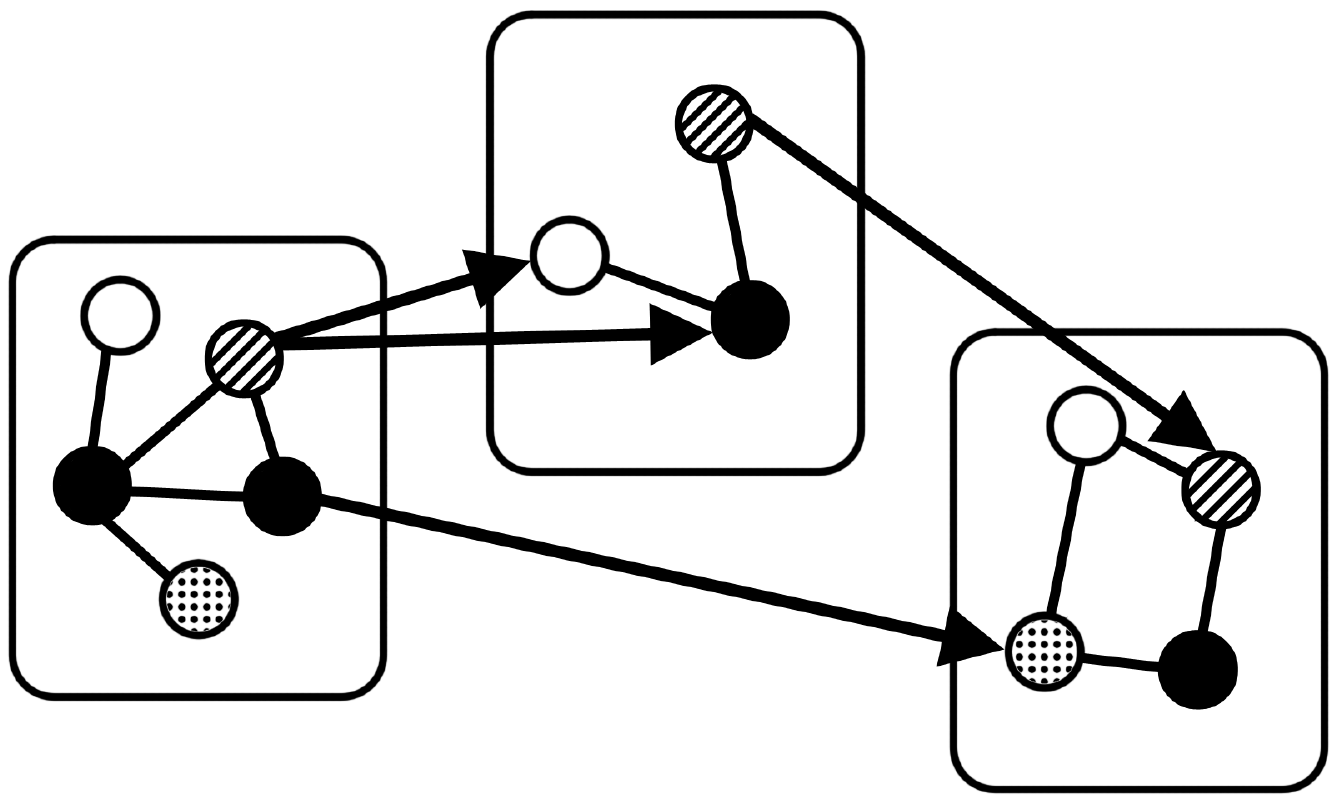}
	}
	\caption{Graphical model structural learning. Given $n$
          i.i.d. samples and 
          the partial DAG orderings as in fgure (a), our goal is to recover
          the graphical model structure of underlying distribution
          such as in figure (b) or figure (c).} 
	\label{Fig:Learning}
\end{figure}

In this section, we investigate fitting our class of BDMRF models in
\eqref{EqnBDMRF} to data, under high-dimensional settings where the
number of variables or nodes in the graph may potentially be larger
than the number of samples. Specifically, suppose we observe $n$
i.i.d. samples $\big\{X^{(j)}\big\}_{j=1}^n$ from an underlying BDMRF
distribution with unknown parameters 
$\theta^*$: 
\begin{align}
	\log \P[X] = \sum_{i=1}^{m} \bigg\{ \sum_{C \in \CV} \theta^*_{C}(\XPa) \prod_{s \in C} B_{s}(X_s)  + \sum_{s \in V_i} M_{s}(X_{s}) - A_{i}\big(\theta^*(\XPa)\big) \bigg\}.
	\label{EqnTrueDist}
\end{align}
Then, our objective is two-fold: (a) \emph{parameter learning}, or to
estimate the unknown parameters $\theta^{*}$, and (ii) \emph{structure
  learning}, or to estimate the unknown edge-set $E$ of the underlying
mixed graph.  These tasks are frequently also referred to as graphical
model estimation and selection respectively.

These two problems are especially challenging for our class of BDMRF
models. As with classical MRFs, the normalization constant for BDMRF
distributions is not available in closed form. Further, our class of
distributions need not even belong to an exponential family as noted
in earlier sections. Thus, estimating our class of models by directly
maximizing the full likelihood is typically intractable. Secondly, the
set of variables in our class of models is not only high-dimensional,
but also heterogeneous, belonging to varied data-types. Lastly, the
edge structure of the mixed graph underlying our class of BDMRF
distributions have both directed and undirected edges - directed
edges between variables from a parent block to variables in a child
block, and undirected edges connecting variables within each
block. Learning the structure of even solely directed edges of
directed graphical models is known to be an NP-hard
problem~\cite{Lauritzen}, unless a partial ordering of the variables
is known \cite{shojaie2010penalized}. Accordingly, in this paper, we
assume that the partial ordering of the blocked DAG is known a-priori,
based on domain knowledge.  This assumption is especially relevant in 
areas such as high-throughput genomics, as we discussed in
Section~\ref{Sec:BDMRFExamples}. Assuming these partial orderings over
the DAG blocks, our learning task is reduced to recovering the
undirected skeleton of the mixed graph in Figure \ref{Fig:Learning}.

By construction, our class of BDMRF distributions in \eqref{EqnBDMRF}
is completely specified by mixed CRF distributions
\eqref{EqnBlockMixedCRF} over the DAG blocks. Accordingly, we can
reduce the problem of estimating the overall BDMRF distribution to
that of estimating the corresponding mixed CRFs, as outlined in
Algorithm~\ref{alg:blockrec}.  In order to learn any of the mixed CRFs
$\P_{\exp}[\XV|\XPa]$, only the sample sub-vectors restricted to $\XV$
and $\XPa$ are required. Therefore, the overall BDMRF estimation
problem can be reduced to the set of sub-problems of estimating the
mixed CRFs comprising the BDMRF: 
\begin{align}
	\log \P_{\exp}[\XV|\XPa] \propto \sum_{C \in \C_{V_i}} \theta^*_{C}(X_{\text{PA}(C)}) \prod_{s \in C} B_{s}(X_s)  + \sum_{s \in V_i} M_{s}(X_{s}).
	\label{EqnTrueBlockDist}
\end{align}


In the following, we focus on pairwise graphical models with linear covariate functions meaning that  
the parameter functions corresponding to node-wise cliques $C := \{s\}$,
have the following form:  

%
%

\begin{align}\label{EqnTrueParam1}
	\theta^*_{C}(\XPa) := \theta^*_s + \sum_{t \in \textnormal{PA}(s)} \theta^*_{st} \,B_t(X_{t}),
\end{align}
while the parameter functions corresponding to pair-wise cliques $C := \{s,s'\}$, will be simply a constant parameterized by $\theta^*_{ss'}$: 
\begin{align}\label{EqnTrueParam2}
	\theta^*_{C}(\XPa) := \theta^*_{ss'}.
\end{align}
Thus, overall the mixed CRF distribution in \eqref{EqnTrueBlockDist},
when restricted to be pairwise, takes the form: 
\begin{align}\label{EqnTrueBlockPairwise}
	\log \P_{\exp}[\XV|\XPa] \propto & \sum_{s \in V_i} \bigg\{ \theta^*_s + \sum_{t \in \textnormal{PA}(s)} \theta^*_{st} \,B_t(X_{t}) \bigg\} B_s(X_s)   \nonumber\\
	& + \sum_{s' \in V_i \text{ s.t. } (s,s') \in E} \theta^*_{ss'} \, B_s(X_s) \, B_{s'}(X_{s'}) + \sum_{s \in V_i} M_s(X_s).
\end{align}

As noted above, the task of estimating a BDMRF distribution from data
reduces to the task of estimating these mixed
CRFs~\eqref{EqnTrueBlockPairwise}; since the graph factors according
to these mixed CRFs, we can estimate each one
independently. We propose to do so following the node-wise
neighborhood estimation approach of
\citep{Meinshausen06,RWL10,YRAL13,YRAL13CRF}, which allows us to
side-step the task of computing the log-partition function of the
mixed CRFs. These neighborhood selection approaches seek to learn the
network structure through an $\ell_{1}$-norm penalty that sparsely
estimates the set of edge parameters, the non-zeros of which
correspond to the selected node-neighbors.  Estimating the block mixed
CRF~\eqref{EqnTrueBlockPairwise} in turn reduces to estimating the
univariate node-conditional 
distributions of variables $X_s$ for $s \in V_i$ given all other nodes in $V_i$ and $\textnormal{PA}(i)$:
have the form: 
\begin{align}\label{EqnTrueBlockDistNode}
	\log \P_{\exp}[X_s | X_{V_i\setminus s}, & \XPa]  \propto  \bigg\{ \theta^*_s + \sum_{t \in \textnormal{PA}(s)} \theta^*_{st} \,B_t(X_{t})  \nonumber\\
	&+ \sum_{s' \in V_i \text{ s.t. } (s,s') \in E}\theta^*_{ss'} \, B_{s'}(X_{s'})  \bigg\}  B_s(X_s) + M_{s}(X_{s}).
\end{align}

Therefore, for each node $X_s$ in $V_i$, the node-conditional
distribution is specified by a parameter vector
$\boldsymbol{\theta^*}(s)$ with three components  
$\boldsymbol{\theta^*}(s) := (\theta^*_s,\boldsymbol{\theta^*}_{V_i},\boldsymbol{\theta^*}_{\textnormal{PA}(i)})$.
 Here, $\theta^*_s \in \real$ is the nodewise weight in the nodewise parameter function in \eqref{EqnTrueParam1}. $\boldsymbol{\theta^{*}}_{V_i} := \{\theta^*_{ss'}\}_{s' \in {V_i \backslash s}} \in \real^{p_i -1}$, where $p_i = |V_i|$, is the vector of intra-block edge-weights in the nodewise parameter function in \eqref{EqnTrueParam1}.
 And finally, $\boldsymbol{\theta^*}_{\textnormal{PA}(i)} := \{ \theta^*_{st} \}_{t \in \textnormal{PA}(i)} \in \real^{q_i - 1}$, where $q_i := |\textnormal{PA}(i)|$, is the vector of inter-block edge weights in the pairwise parameter function in \eqref{EqnTrueParam2}.

Let $\mathcal{N}^{*}_{V_i}(s) = \{s' \in {V_i \backslash s}: \theta^*_{ss'} \neq 0\}$ denote the set of true intra-block neighbors of $s \in V_i$. Let $d_{V_i} = \max_{s \in V_i}|\mathcal{N}^{*}_{V_i}(s)|$ denote the number of these intra-block neighbors. Similarly, let $\mathcal{N}^{*}_{\textnormal{PA}(i)}(s) = \{t \in \textnormal{PA}(i): \theta^*_{st} \neq 0\}$ denote the set of true inter-block neighbors of $s \in V_i$. Let $d_{\textnormal{PA}(i)} = \max_{s \in V_i}|\mathcal{N}^{*}_{\textnormal{PA}(i)}(s)|$ denote the number of these inter-block neighbors.

Then, given $n$ i.i.d. samples from our pairwise BDMRF of the form in \eqref{EqnTrueBlockPairwise}  with unknown parameters, we estimate the node-conditional distributions in \eqref{EqnTrueBlockDistNode} by solving for the $\ell_1$-regularized conditional MLEs:
\begin{align}\label{EqnNeighborSelection}
	\min_{ \boldsymbol{\theta}(s) \in \real^{1 + p_i + q_i}} -\frac{1}{n}\sum_{j=1}^{n} & \log \P\Big(X_s^{(j)}  | X_{V_i\setminus s}^{(j)}, \XPa^{(j)} ;\theta \Big) +  \regV \|\boldsymbol{\theta}_{V_i}\|_1  + \regPa \|\boldsymbol{\theta}_{\textnormal{PA}(i)}\|_1 \, ,
\end{align}
where $\regV$, $\regPa$ are
regularization constants, with $\regV$ determining the degree
of sparsity in the connections between $X_s$ and $X_{V_i\setminus s}$,
and $\regPa$ determining the degree of sparsity
in the connections between $X_s$ and the nodes in
$\textnormal{PA}(i)$.  Note that if the sufficient statistics are
linear, $B_{s}(X_{s}) = X_{s}$, then \eqref{EqnNeighborSelection} can
in turn be written in via exponential family natural parameters,
$\boldsymbol{\beta}$,  in the form of an $\ell_{1}$ regularized
Generalized Linear Model as first presented in \citet{YRAL12,YRAL13CRF}:
\begin{align*}
	\min_{ \boldsymbol{\beta}(s) \in \real^{1 + p_i + q_i}}
        & -\frac{1}{n}\sum_{j=1}^{n}   \ell \left( X_s^{(j)} ; \beta_s
        + \sum_{t \in \textnormal{PA}(i) } \beta_{st}
        \,B_t\big(X_{t}^{(j)}\big)  + \sum_{s' \in V_i } \beta_{ss'}
        \, B_{s'}\big(X_{s'}^{(j)}\big) \right)  \\
	&  +  \regV \|\boldsymbol{\beta}_{V_i}\|_1 + \regPa \|\boldsymbol{\beta}_{\textnormal{PA}(i)}\|_1 \, .
\end{align*}

Given the solution $\widehat{\boldsymbol{\theta}}(s)$ of the $M$-estimation problem \eqref{EqnNeighborSelection} above, we can obtain an estimate of the true intra-block neighbors $\mathcal{N}^{*}_{V_i}(s)$ of a node $s \in V_i$ by:
\[ \widehat{\mathcal{N}}_{V_i}(s) = \{s' \in {V_i \backslash s}: \widehat{\boldsymbol{\theta}}_{ss'} \neq 0\},\]
and obtain an estimate of the true inter-block neighbors $\mathcal{N}^{*}_{PA(i)}(s)$ of $s \in V_i$ by:
\[ \widehat{\mathcal{N}}_{PA(i)}(s) = \{t \in \textnormal{PA}(i): \widehat{\boldsymbol{\theta}}_{st} \neq 0\}.\]

As we show below, these graph-structure estimates $\{\widehat{\mathcal{N}}_{V_i}(s), \widehat{\mathcal{N}}_{PA(i)}(s)\}$ from our $M$-estimator come with strong statistical guarantees. 

\begin{theorem}\label{Thm:Main}
Consider a pairwise mixed CRF distribution as specified in \eqref{EqnTrueBlockPairwise}. Suppose we solve the $M$-estimation problem in \eqref{EqnNeighborSelection} with the regularization parameters set as
\begin{align*}
	\regV \geq M_1 \sqrt{\frac{\log p_i}{n}} \quad \regPa \geq M_1 \sqrt{\frac{\log q_i}{n}}, \quad \textrm{and} \quad \max\{\regV,\regPa\} \leq M_2,
\end{align*}
where $M_1$ and $M_2$ are some constants that depend on the types of exponential family in \eqref{EqnTrueBlockDistNode}. Further suppose that the minimum intra-block and inter-block edge-weights satisfy:
\begin{align*}
	\min_{t \in V_i\setminus s  \cup \textnormal{PA}(i)} |\theta_{st}^{*}| \ge \frac{10}{\Cmin}\max\Big\{\regV \sqrt{d_{V_i}} \, , \, \regPa\sqrt{d_{\textnormal{PA}(i)}}\Big\},
\end{align*}
where $\Cmin$ is the minimum eigenvalue of the Hessian of the loss
function at $\boldsymbol{\theta}^*(s)$, where we recall the notation that $p_i$ and $q_i$ are the number of intra-block and inter-block neighbors respectively of node $s \in V_i$. Then, for some universal positive constants $L$, $c_1$, $c_2$, and $c_3$, if 
\[ n \geq L \big(d_{V_i} + d_{\textnormal{PA}(i)}\big)^{2} \big(\log p_i + \log q_i\big) \Big(\max\big\{\log n,\log (p_i+q_i)\big\}\Big)^{2},\] then with probability at least $1- c_1 \max\{n,p_i+q_i\}^{-2} - \exp(-c_2 n) - \exp(-c_3 n)$, the solution $\widehat{\boldsymbol{\theta}}$ of the M-estimation problem in \eqref{EqnNeighborSelection} satisfies the following:
	\begin{description}[leftmargin=1cm,noitemsep]
		\item[(a)] \textnormal{(Parameter Error)} For each node $s \in V$, the solution $\widehat{\boldsymbol{\theta}}$ is unique with error bound:
			\begin{align*}
	\|\widehat{\boldsymbol{\theta}}_{\textnormal{PA}}-\boldsymbol{\theta}^*_{\textnormal{PA}}\|_2 + \|\widehat{\boldsymbol{\theta}}_{V_i}-\boldsymbol{\theta}^*_{V_i}\|_2
				\leq \,  \frac{5}{\Cmin} \max\Big\{ \regV\sqrt{d_{V_i}} \, , \, \regPa\sqrt{d_{\textnormal{PA}(i)}} \Big\}
			\end{align*}
		\item[(b)] \textnormal{(Structure Recovery in $V_i$)} The solution recovers the true intra-block neighborhoods exactly, so that $\widehat{\mathcal{N}}_{V_i}(s) = \mathcal{N}^{*}_{V_i}(s)$, for all $s \in V_i$.
		\item[(c)] \textnormal{(Structure Recovery between $V_i$ and $\textnormal{PA}(i)$)} The solution recovers the true inter-block neighborhoods exactly, so that $\widehat{\mathcal{N}}_{\textnormal{PA}(i)}(s) = \mathcal{N}^{*}_{\textnormal{PA}(i)}(s)$, for all $s \in V_i$.
	\end{description}
\end{theorem}

Note that since the neighborhood of each node is recovered with
high-probability, the entire graph structure estimate in turn is
recovered exactly with high-probability, by a simple application of a
union bound.  These statistical guarantees extend the graph estimation
and selection results from \citet{YRAL12,YRAL13CRF} to our BDMRF
setting.

\section{Numerical Examples}\label{Sec:SimData}

\begin{figure}[t!]
	\subfigure[Gau-Ising Mixed MRF]{\includegraphics[width=.32\linewidth]{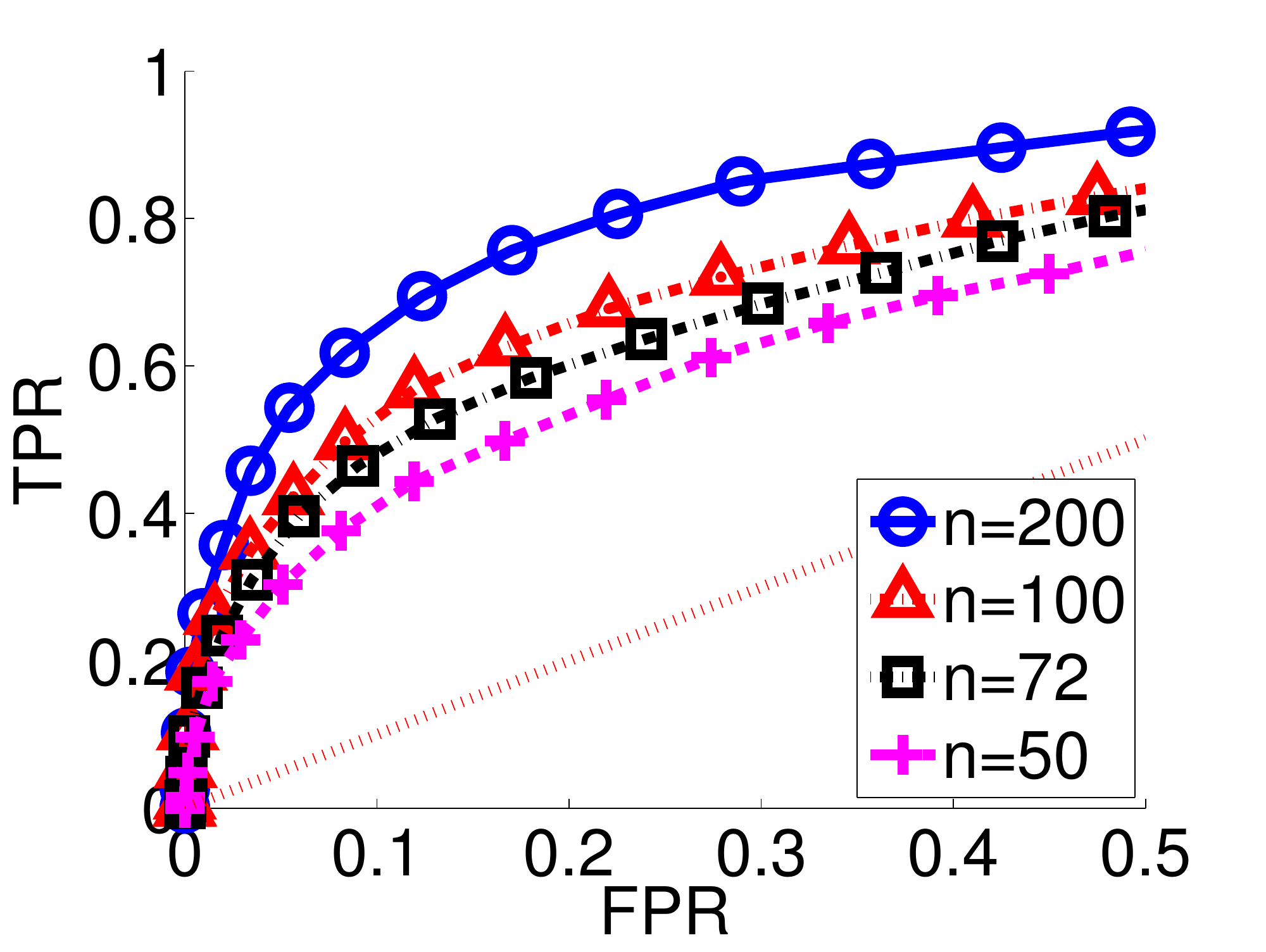}}\hfill
	\subfigure[Gau MRF-Ising CRF]{\includegraphics[width=.32\linewidth]{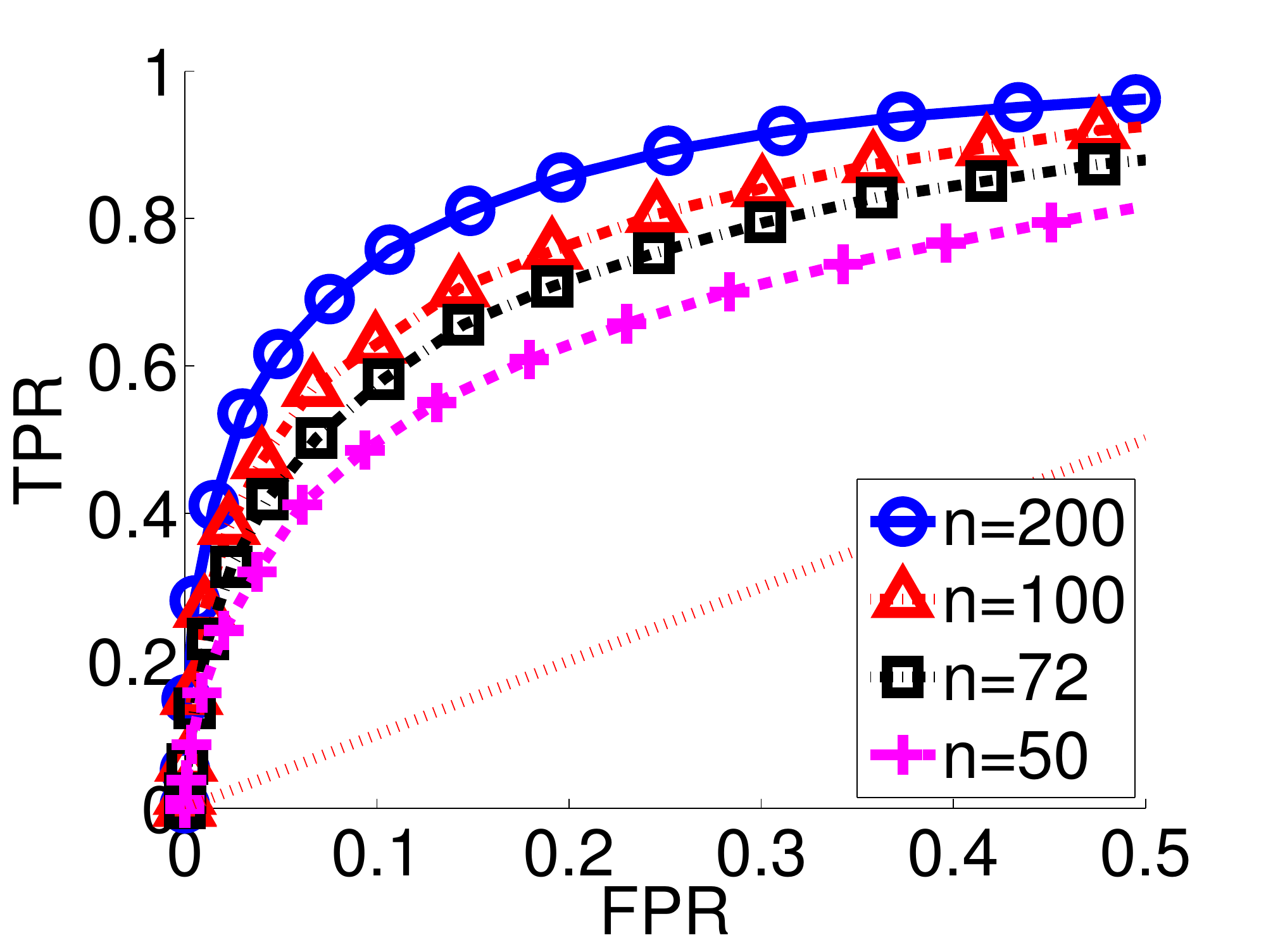}}\hfill
	\subfigure[Gau CRF-Ising MRF]{\includegraphics[width=.32\linewidth]{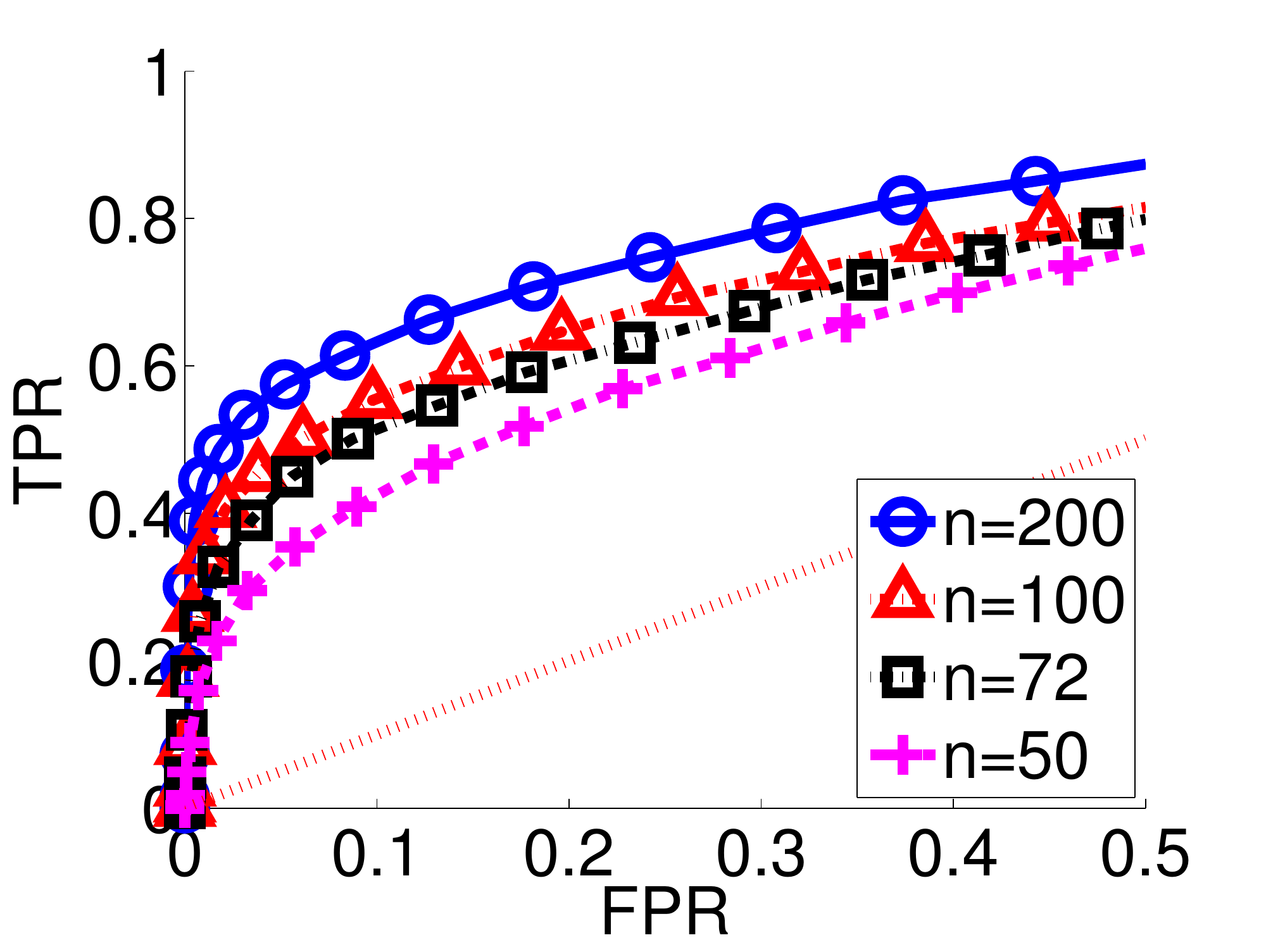}}
	\subfigure[Poi-Ising Mixed MRF]{\includegraphics[width=.32\linewidth]{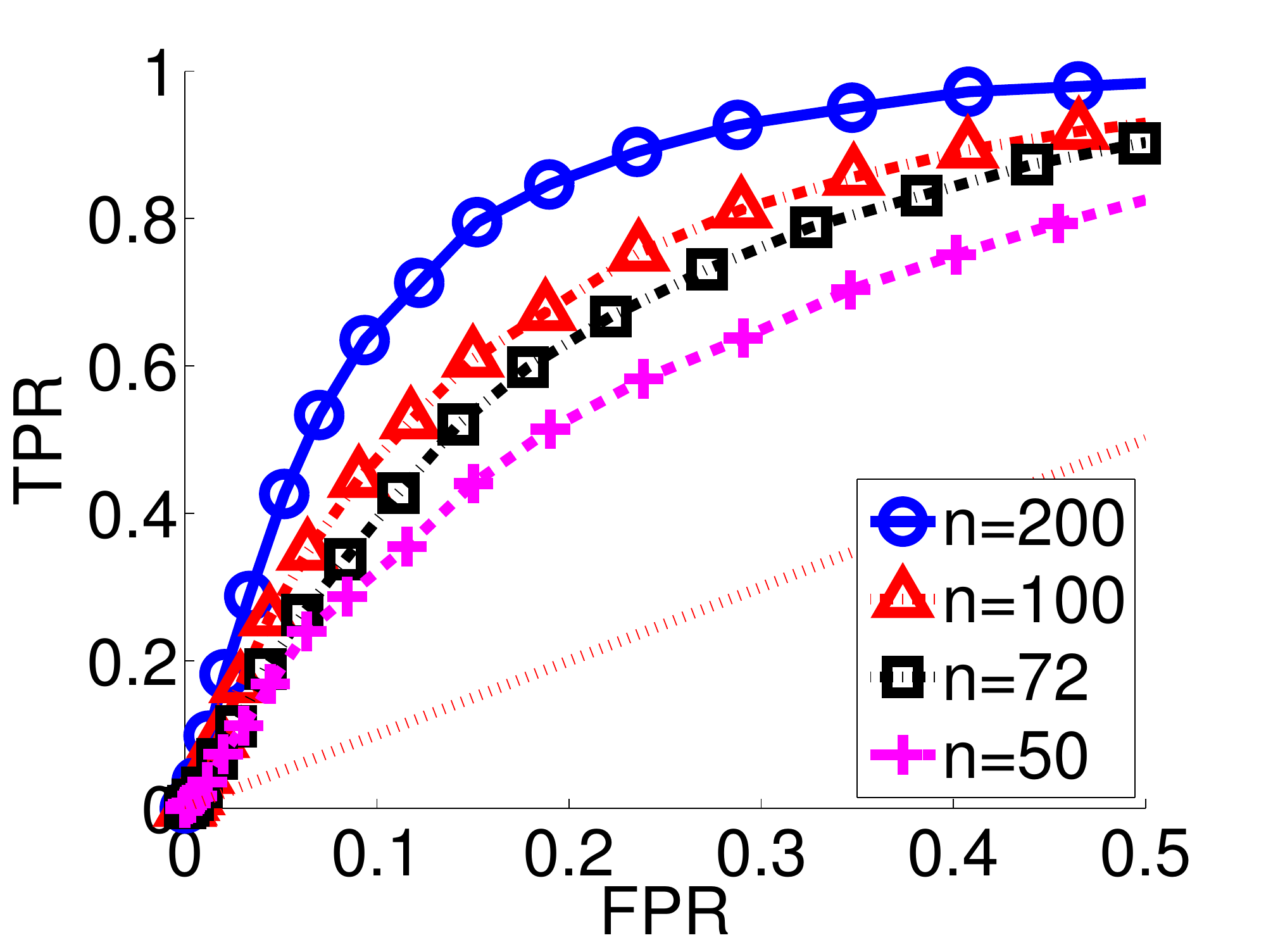}}\hfill
	\subfigure[Poi MRF-Ising CRF]{\includegraphics[width=.32\linewidth]{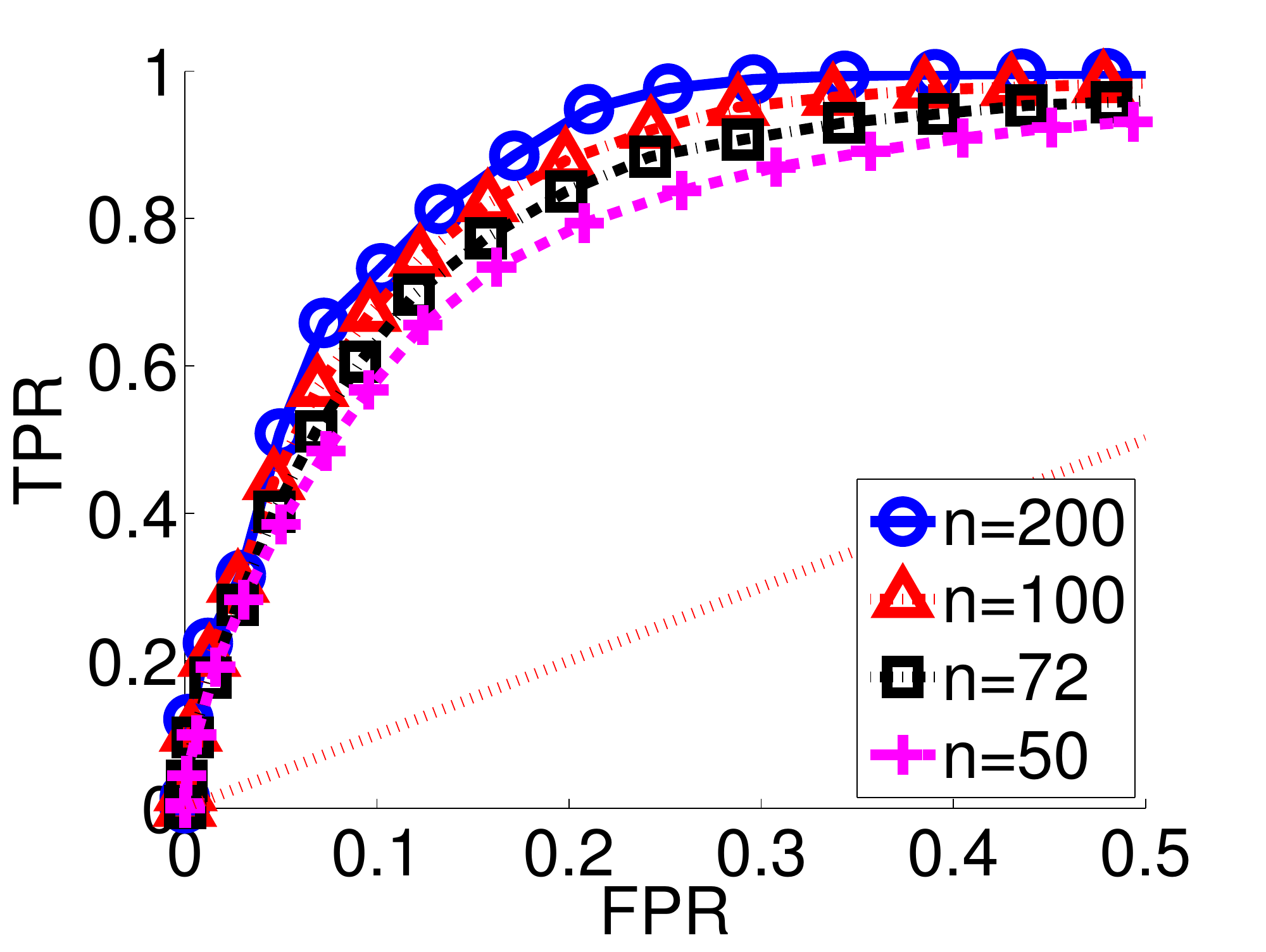}}\hfill
	\subfigure[Poi CRF-Ising MRF]{\includegraphics[width=.32\linewidth]{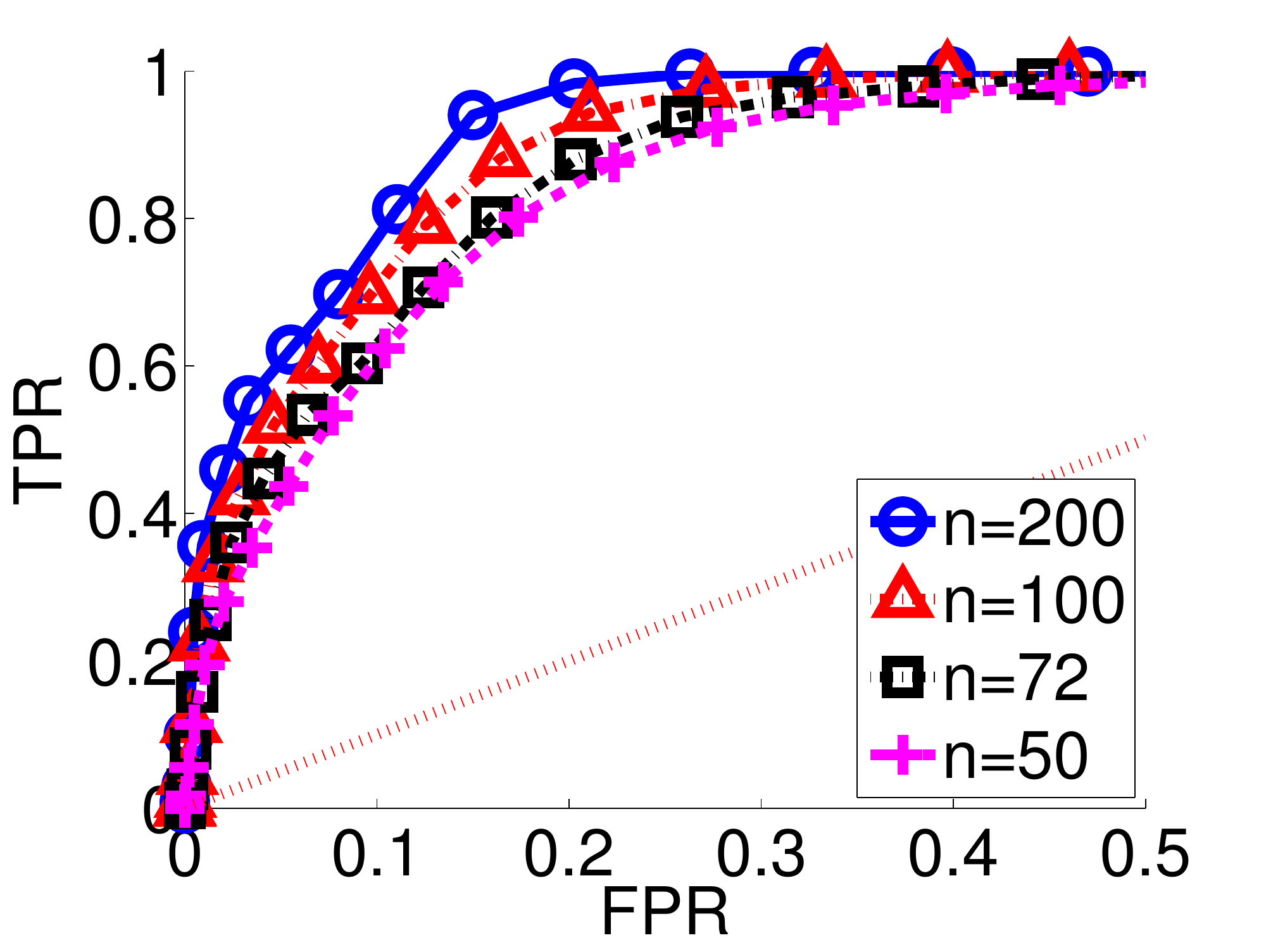}}
	\subfigure[Gau-TPGM Mixed MRF]{\includegraphics[width=.32\linewidth]{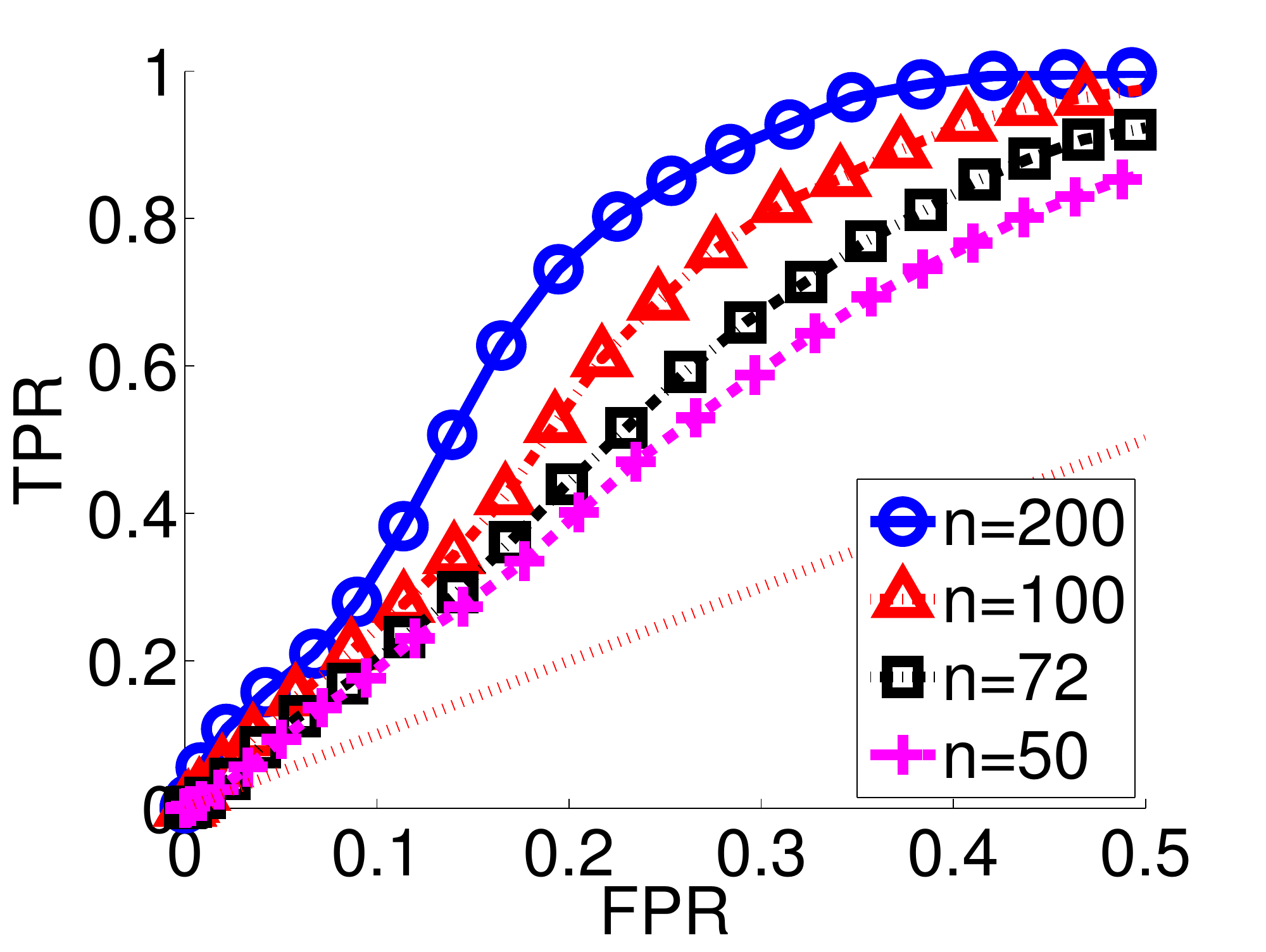}\label{Fig:TPGM-Gau}}\hfill
	\subfigure[Gau MRF-Poi CRF]{\includegraphics[width=.32\linewidth]{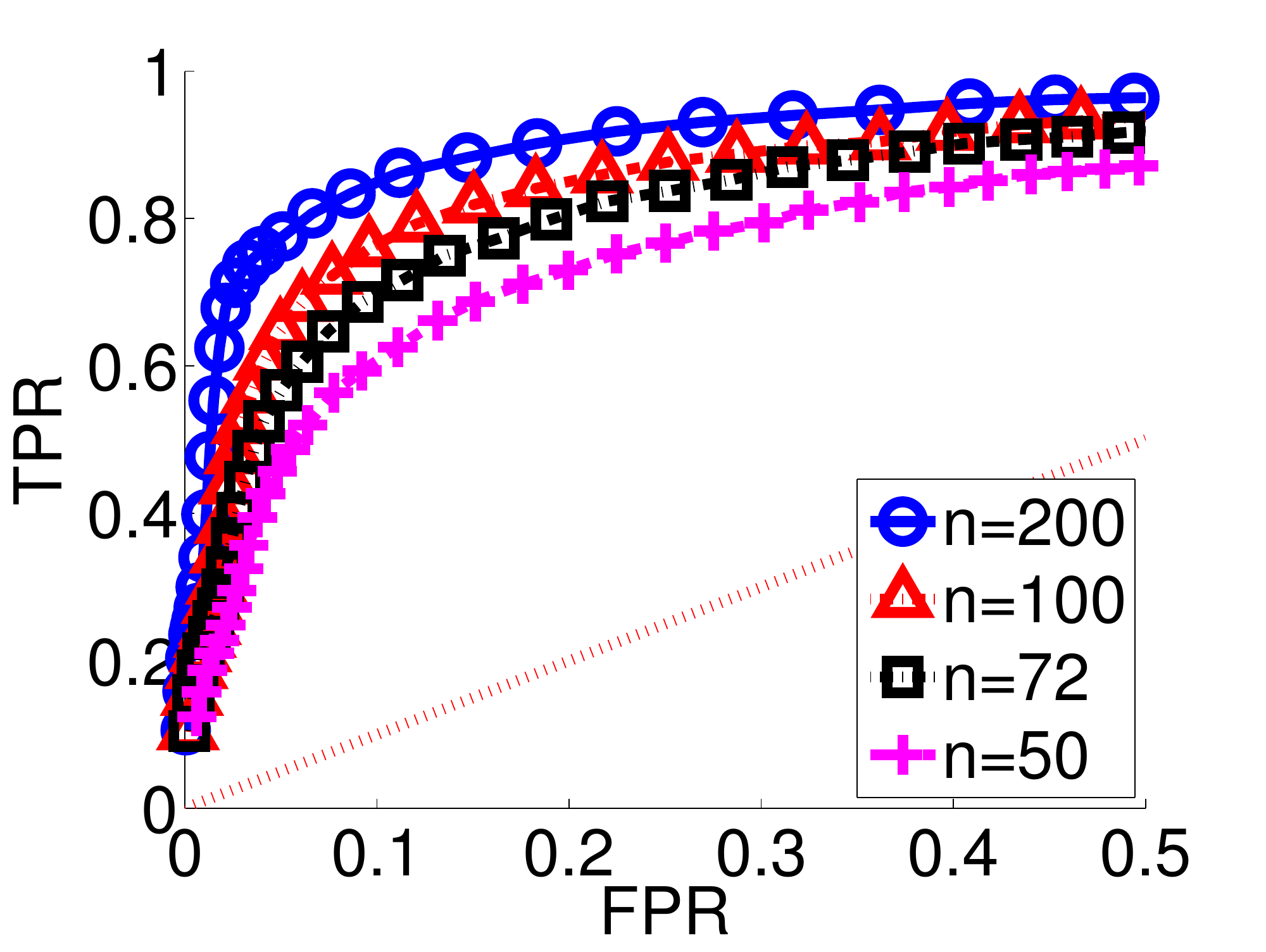}\label{Fig:Pos-Gau1}}\hfill
	\subfigure[Gau CRF-Poi MRF]{\includegraphics[width=.32\linewidth]{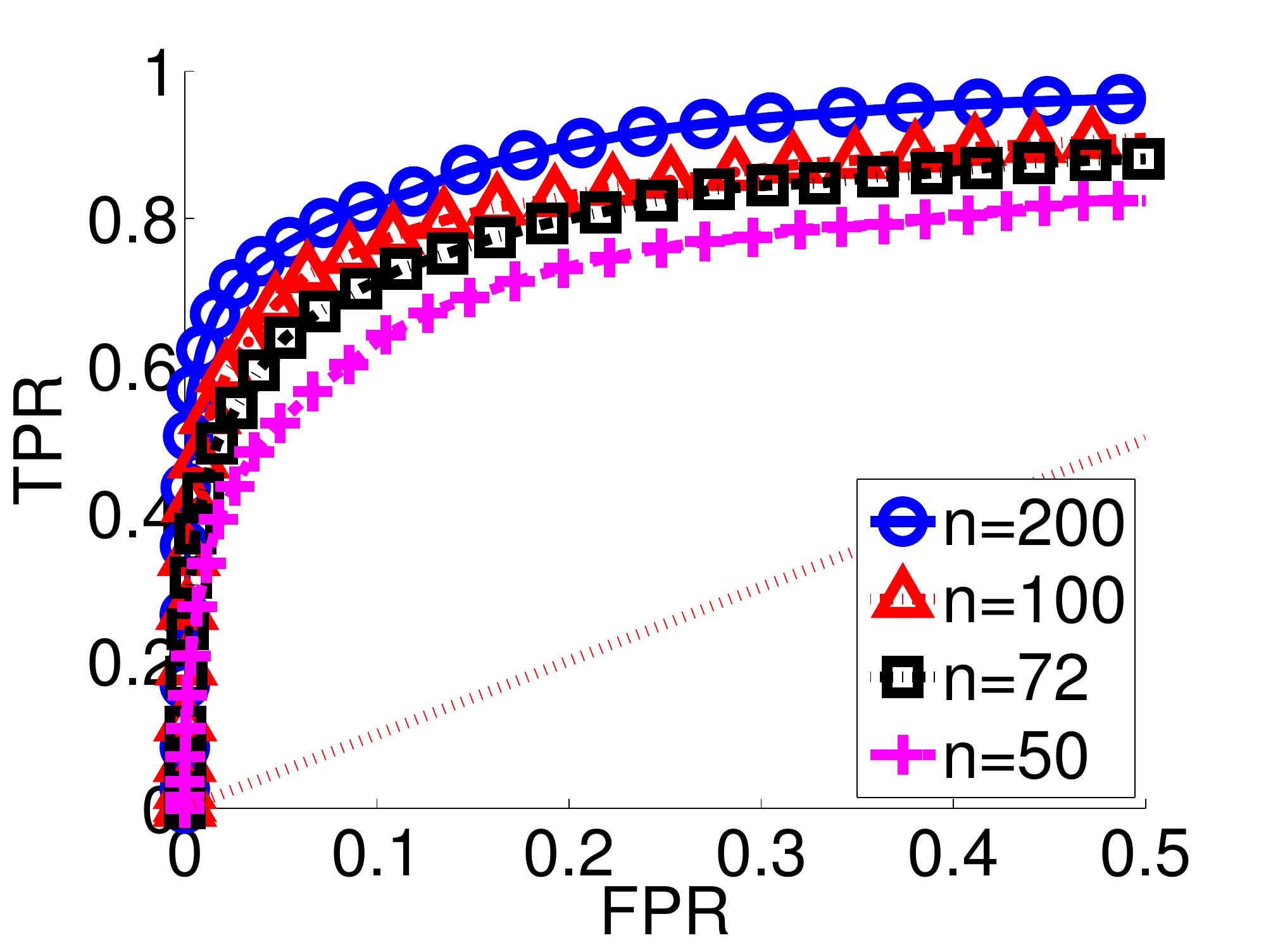}\label{Fig:Pos-Gau2}}
	\caption{Receiver Operator Characteristic (ROC) curves for graph selection
	  of various examples of homogeneous mixed MRFs and our EBDMRFs whose form is
	  specified by the MRF and CRF component as labeled(Gau indicates Guassian and Poi does Poisson).  Data was
	  generated from our models via Gibbs samplers with 
	  $p = 72$ nodes on a lattice structure; curves are generated by
	  varying the regularization parameter and shown for several samples sizes,
	  $n$. } 
	\label{fig:sim2way}
\end{figure}

We first evaluate the performance of our $M$-estimator in the previous
section for estimating our class of BDMRF models from data through a
series of simulated numerical examples.  All simulated examples used a
two-dimensional lattice graph structure.  Care was taken to ensure that all
parameters meet normalizability conditions discussed in
Section~\ref{SecEBDMRFNormalizability}; specific parameter values for
all our simulated examples are given in Appendix~\ref{AppSec:SimSettings}. We
generated samples from the models via a Gibbs 
sampler based on the node conditional distributions specified in
\eqref{EqnTrueBlockDistNode}. Our $M$-estimator, based on
$\ell_{1}$-penalized node-conditional maximum likelihood as described
in Section~\ref{Sec:Learning}, was implemented via projected gradient
descent~\cite{beck2009fast}. 

\subsection{EBDMRFs}
Our first set of simulations are for our class of Elementary Block
Directed Markov Random Fields (EBDMRFs), with two blocks of
variables, as described in Section~\ref{Sec:ElemChain}. We constructed models
with $p = 72$ variables, with two blocks of $p_{x} =  36$ and $p_{y} =
36$ variables each. For various instances of our model, we generated
datasets of four different sample sizes, $n = 50, 72, 100, 200$. We
computed Receiver operator characteristic (ROC) curves for graph
structure recovery averaged over 10 replicates by varying the
regularization parameter as shown in Figure~\ref{fig:sim2way}.  The
top two panels display graph structure recovery for three instances of our
EBDMRF class of distributions $\P[X,Y]$ over the two blocks of
variables $(X,Y)$, by varying the block-DAG over the two blocks. The
left panels use a mixed MRF~\cite{YBRAL14}, where all the variables
are in a single block, the middle panels use the decomposition $\Pe[Y
  | X]\, \Pe[X]$, while the right panels use the decomposition $\P[X,Y]
= \Pe[X | Y] \,\Pe[Y]$, employing products of mixed CRFs and mixed
MRFs. As the figure shows, our $M$-estimators (specified by the
corresponding model) are able to recover the underlying network
structure, even in high-dimensional sampling regimes.

The bottom panel in particular highlights the modeling advantages of
our class of models when compared with the mixed
MRFs~\cite{YBRAL14}. As we have noted earlier in
Section~\ref{SecEBDMRFExamples}, the mixed MRF framework does
not permit dependencies between Gaussian and Poisson type variables
for reasons of normalizability. Our class of EBDMRF distributions, on
the other hand, allow for such dependencies. Figure \ref{Fig:Pos-Gau1} and Figure \ref{Fig:Pos-Gau2} show graph structure recovery for EBDMRFs through Gaussian MRF-Poisson CRF, and Gaussian CRF-Poisson MRF, 
respectively; these represent instances where the corresponding mixed
MRF does not permit inter-data-type 
dependencies. Figure \ref{Fig:TPGM-Gau} provides an additional way to model
inter-data-type dependencies using the class of mixed MRFs, by
using a truncated variant of the Poisson MRF instead as introduced in
~\cite{YRAL13PGM}. Overall, as these simulations suggest, our class of
EBDMRFs permit a wide class of normalizable distributions.

\subsection{BDMRFs}
\begin{figure}[t!]
	\centering
	\subfigure[{$\P[X|Y,Z]\, \P[Y|Z]\, \P[Z]$}]{\includegraphics[width=.32\linewidth]{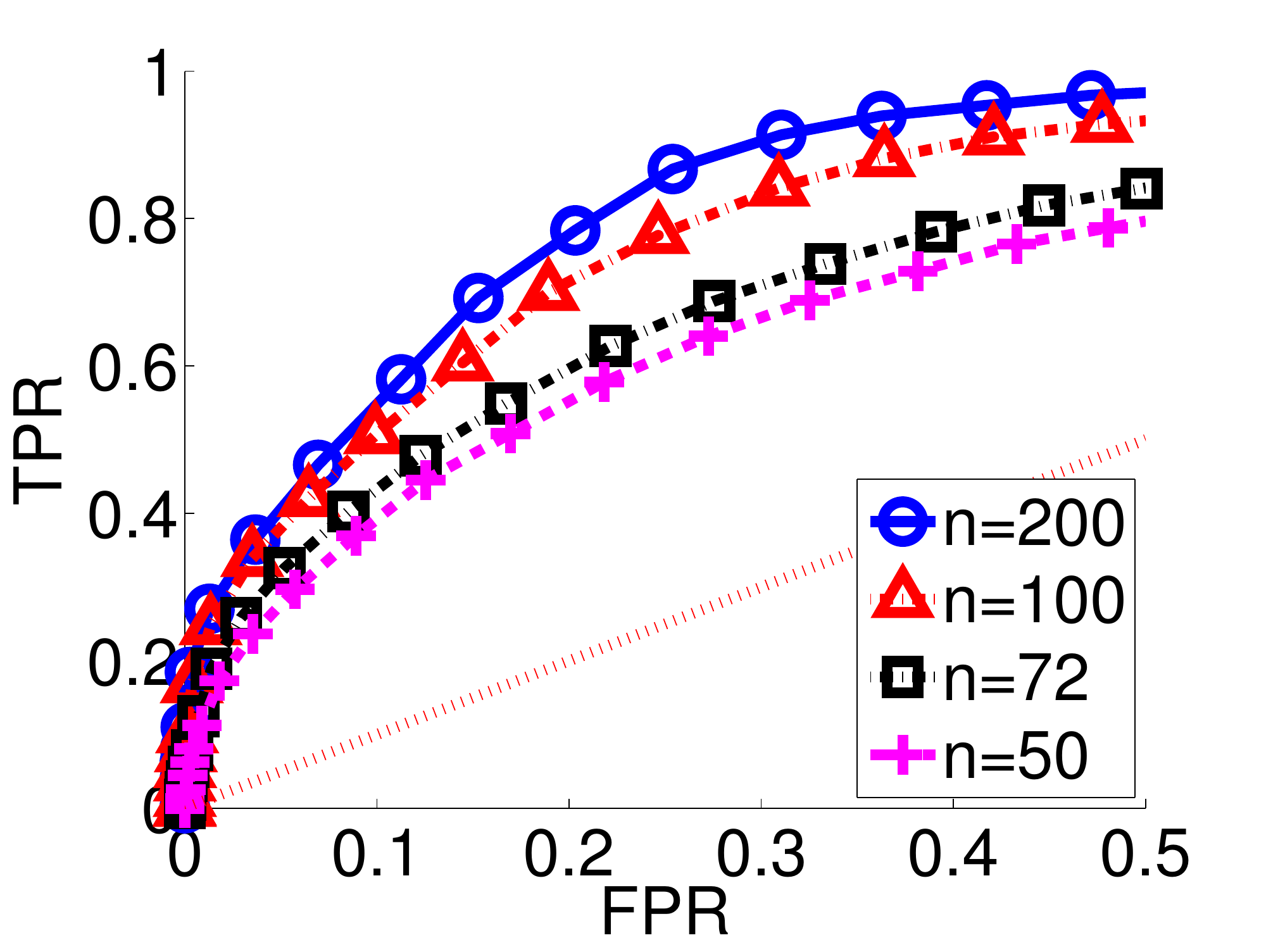}}\hfill
	\subfigure[{$\P[Y|X,Z]\, \P[X|Z]\, \P[Z]$}]{\includegraphics[width=.32\linewidth]{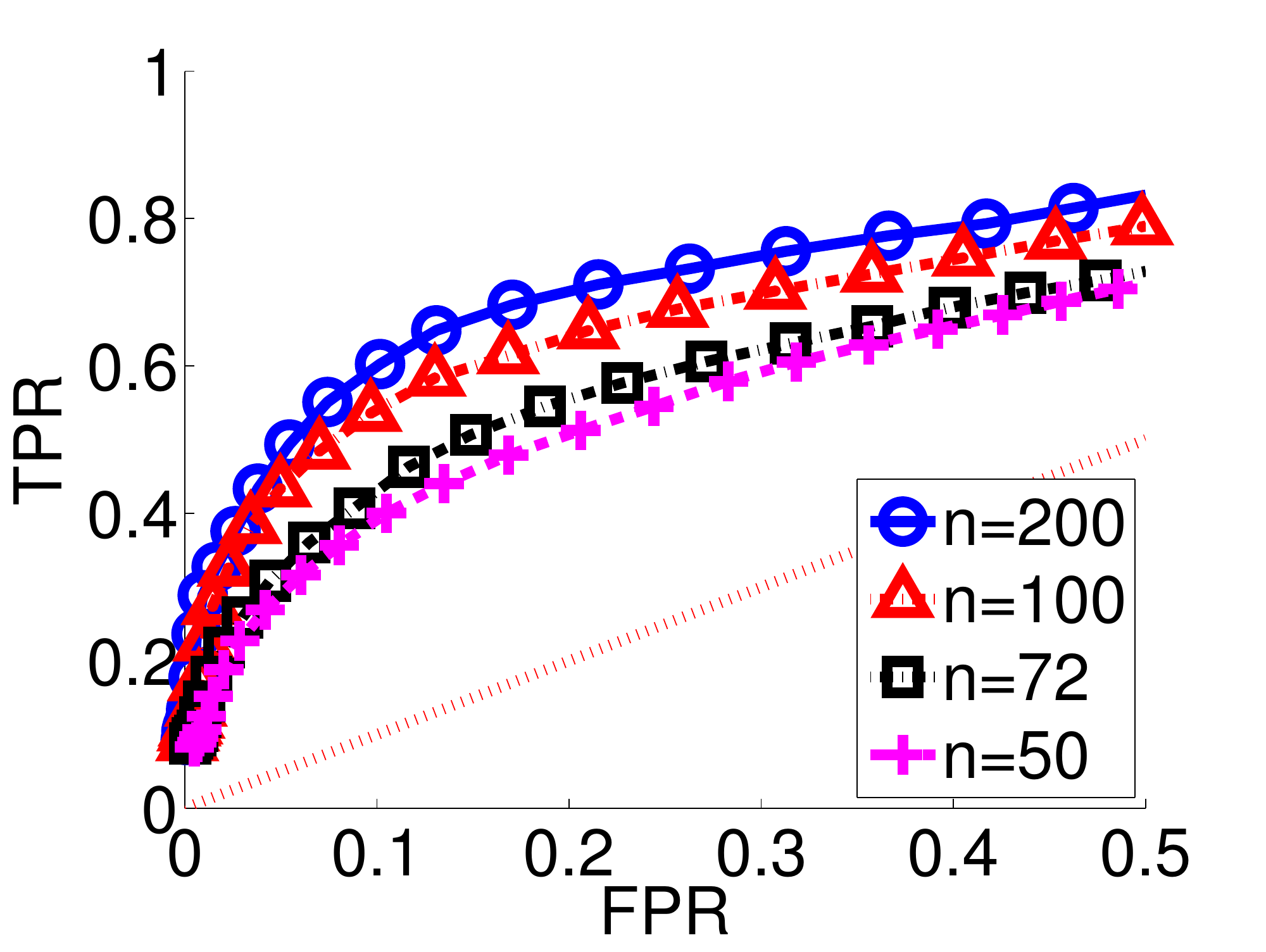}}\hfill
	\subfigure[{$\P[Y|X]\, \P[X|Z]\, \P[Z]$}]{\includegraphics[width=.32\linewidth]{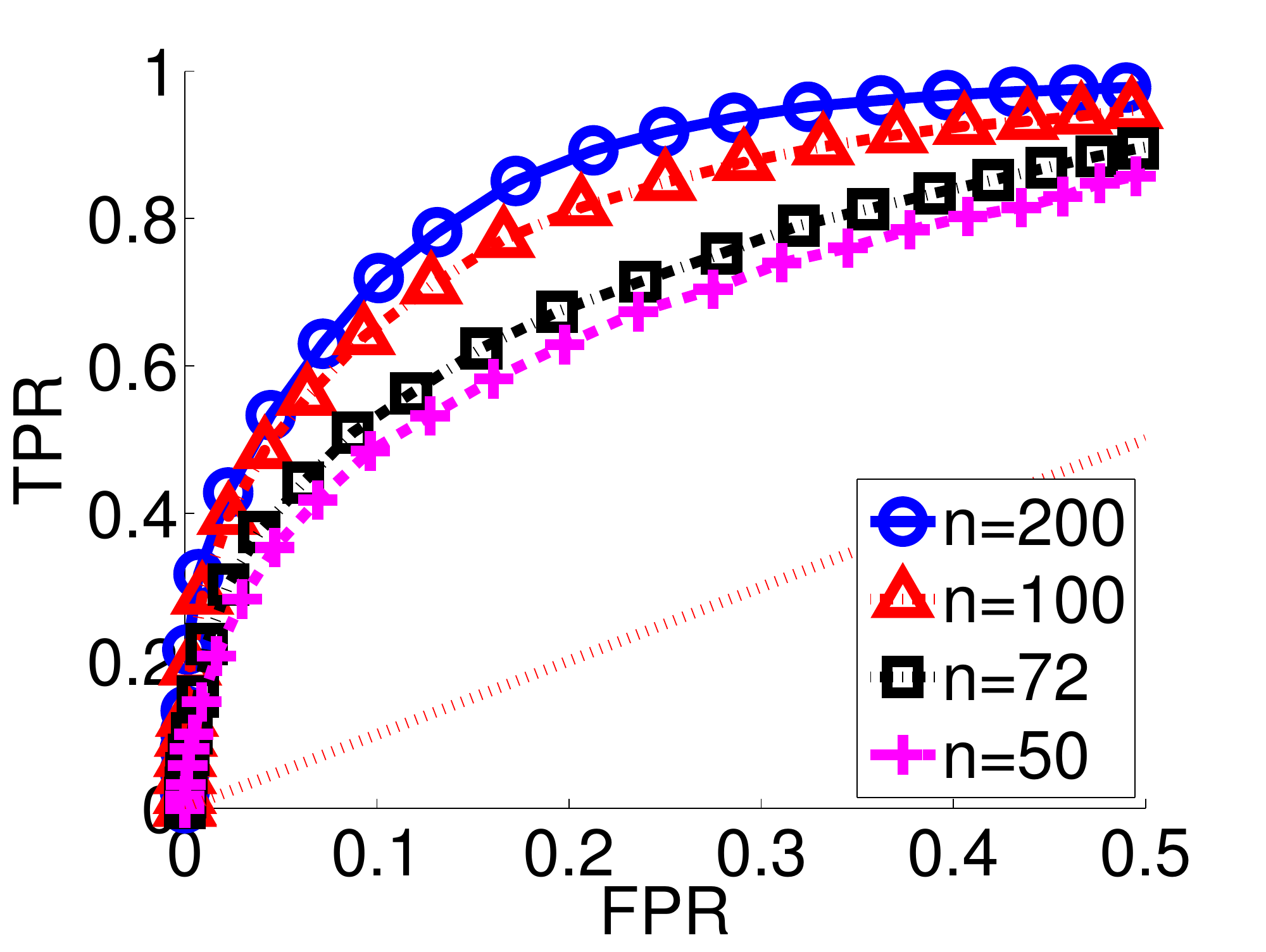}}
	\subfigure[{DAG structure of (a)}]{\includegraphics[width=.32\linewidth]{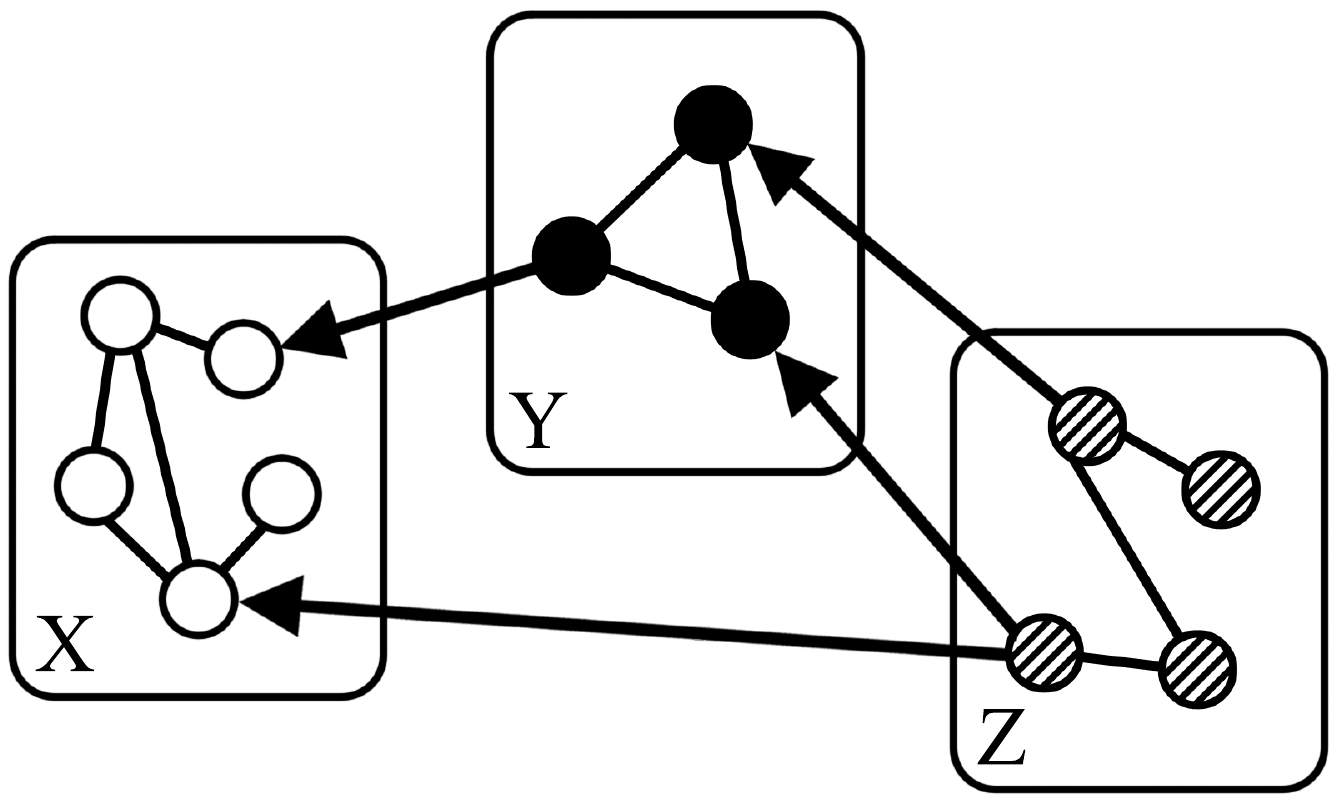}}\hfill
	\subfigure[{DAG structure of (b)}]{\includegraphics[width=.32\linewidth]{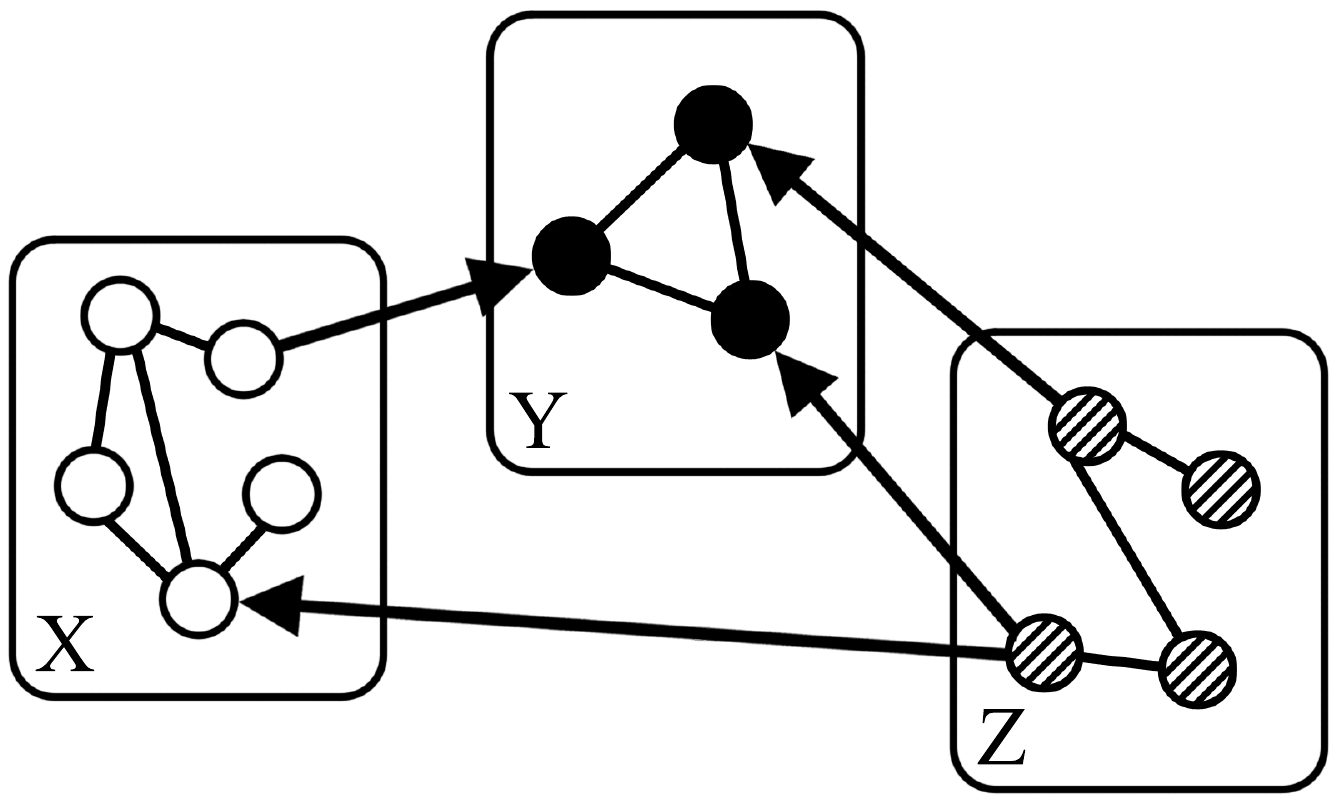}}\hfill
	\subfigure[{DAG structure of (c)}]{\includegraphics[width=.32\linewidth]{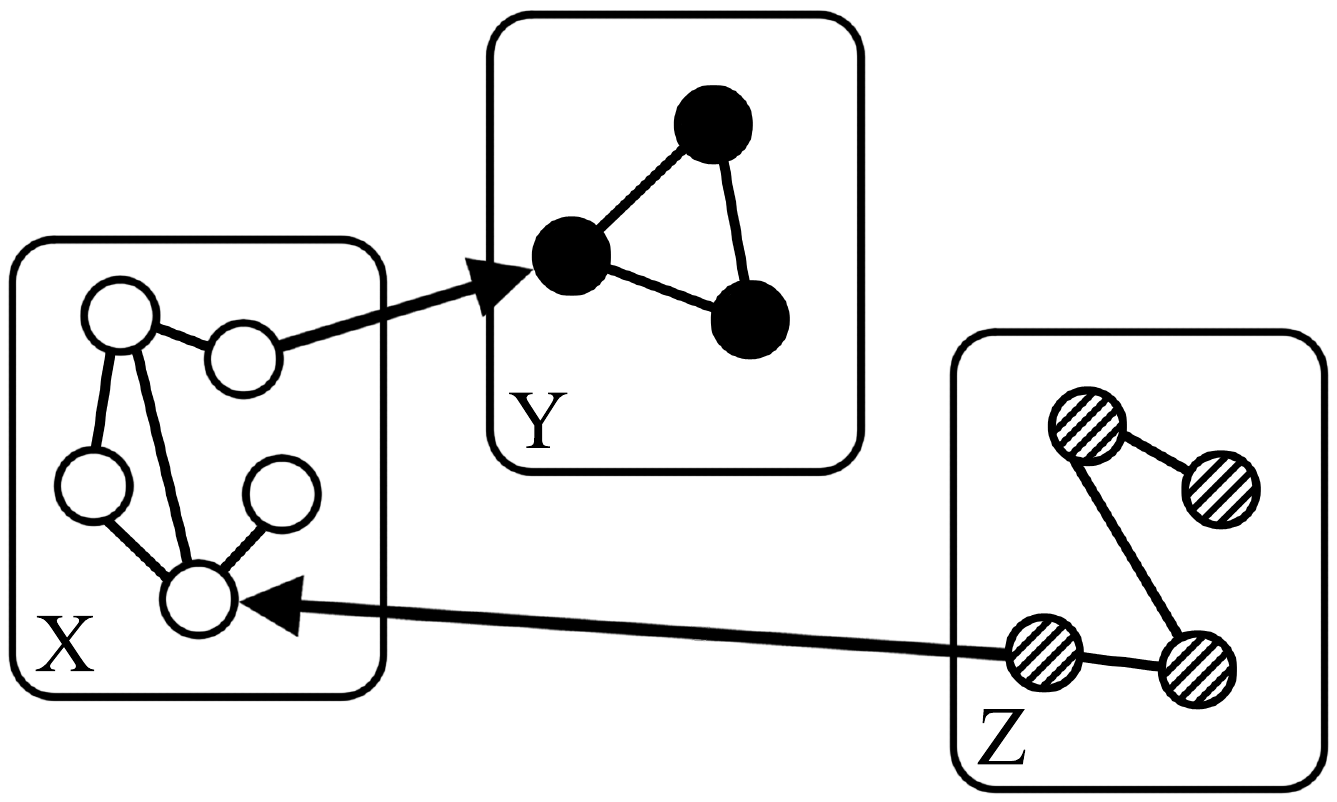}}
	\caption{Receiver Operator Characteristic (ROC) curves (top) for graphical
	  structural learning for various configurations of our BDMRFs for
	  three homogeneous blocks of variables: binary (Ising, $X$), continuous 
	  (Gaussian, $Y$) and counts (Poisson, $Z$).  Model formulations with
	  an illustrative picture are given in the bottom panel. }
	\label{fig:sim3way}
\end{figure}

Our next set of simulations study graph selection recovery for our
general class of BDMRF distributions with three homogeneous blocks of
variables where we know the partial ordering over the variable
blocks. Data was again generated via Gibbs sampler with a lattice
graph structure connecting nearest neighbors on a two dimensional
grid and with $p=75$  variables split into three blocks of $p_X = 25$,
$p_Y = 25$, and $p_Z = 25$ variables. In Figure~\ref{fig:sim3way}, we
plot ROC curves for graph structure recovery of three instances of our
model with three blocks of variables: binary (for which we use Ising CRFs
and MRFs), counts (for which we use Poisson CRFs and MRFs) and continuous
(for which we use Gaussian CRFs and MRFs). The specific models with
illustrated directionality are given in the bottom panel. These
results show that our class of $M$-estimators are able to reasonably
recover the underlying network even under extremely sample-limited
settings. We note that as before, we cannot model dependencies between
these three sets of variables in the classical mixed MRF framework for
reasons of normalizability.  Thus again, these 
illustrate the breadth and flexibility of our class of models.

\subsection{BDMRFs: When the DAG Partial Ordering is Unknown}
\begin{figure}
	\subfigure[{\scriptsize Mixed MRF}]{\includegraphics[width=.33\linewidth]{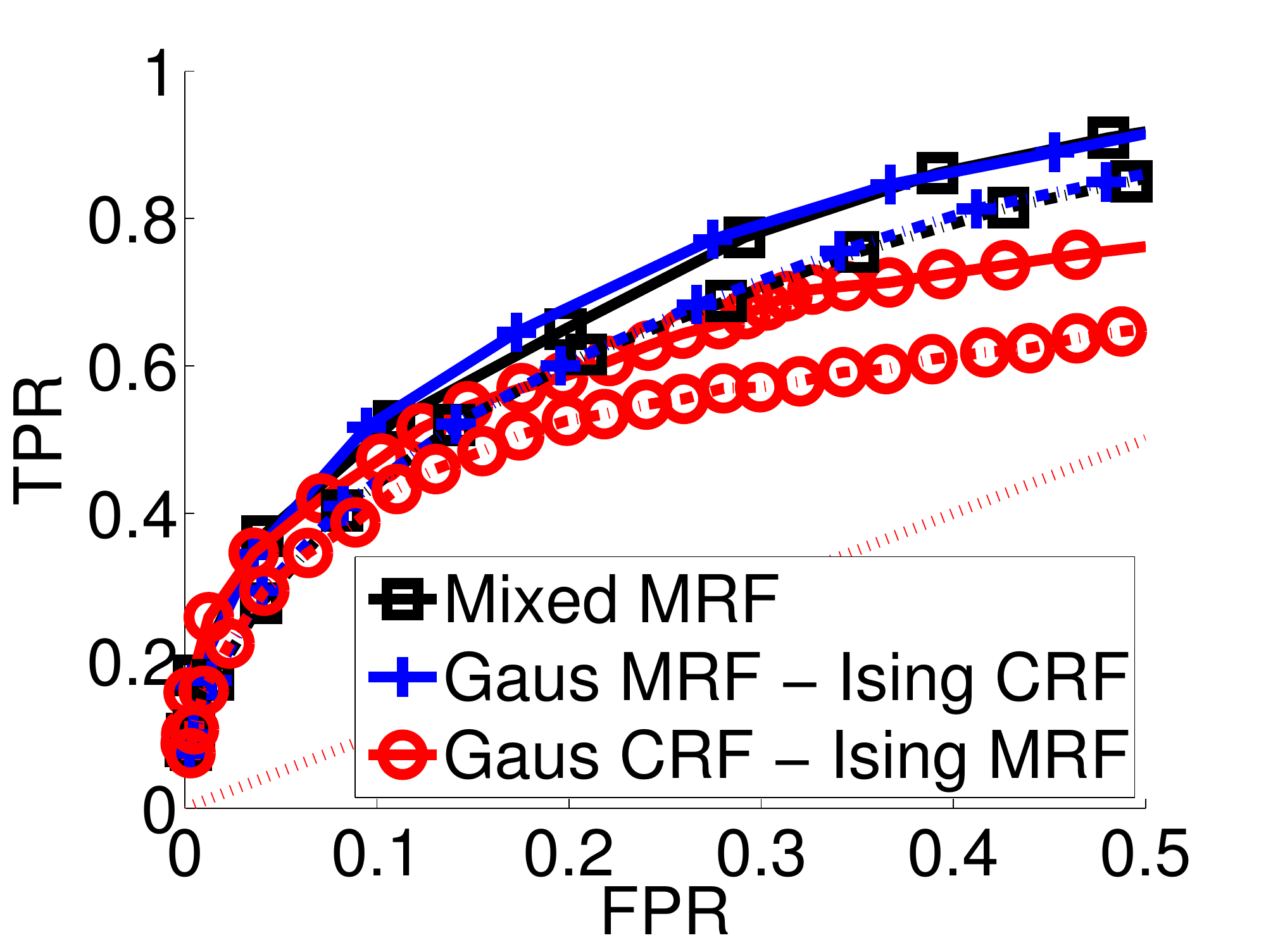}}\hfill
	\subfigure[{\scriptsize Gaussian MRF-Ising CRF}]{\includegraphics[width=.33\linewidth]{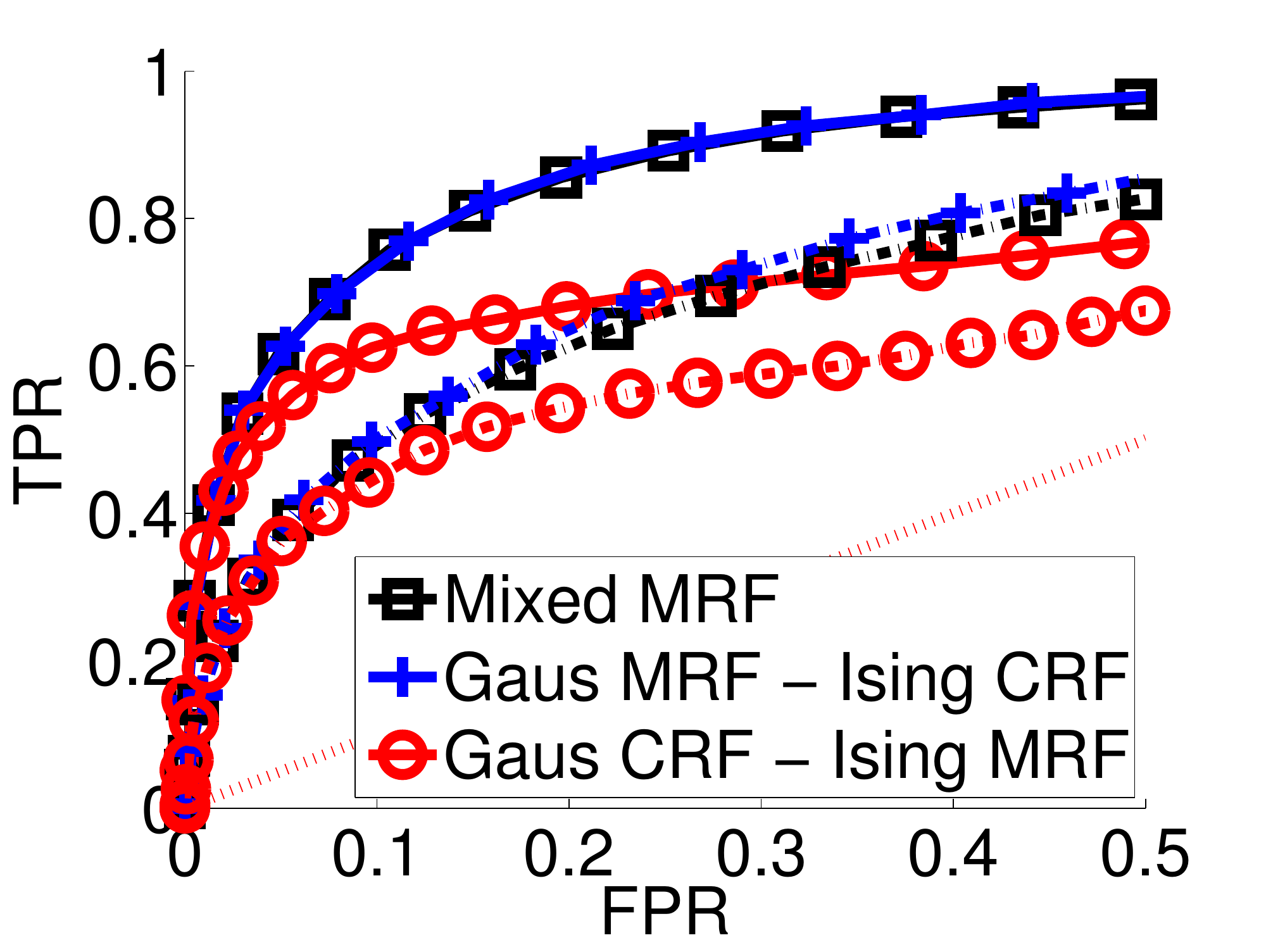}}\hfill
	\subfigure[{\scriptsize Gaussian CRF-Ising MRF}]{\includegraphics[width=.33\linewidth]{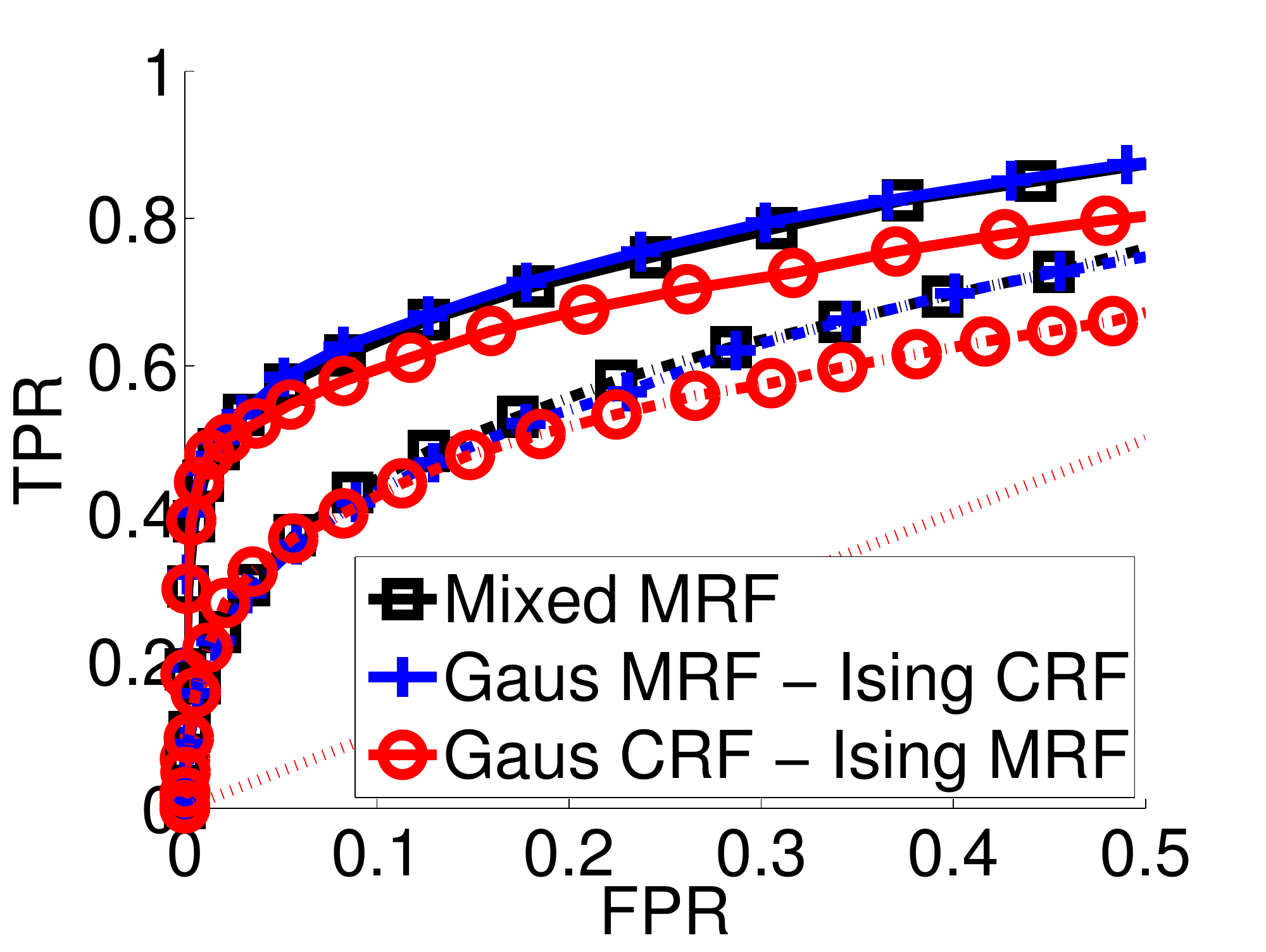}}
	\caption{Receiver Operator Characteristic (ROC) curves for graphical
	  structural learning with two sets of homogeneous variables, binary
	  (Ising) and continuous (Gaussian) when the true model (a) mixed MRF (b) Gaussian MRF-Ising CRF or (c) Gaussian CRF-Ising MRF, is unknown. Solid lines show the results for $n=200$ and dotted lines do for $n=50$.}
	\label{fig:sim_unknown}
\end{figure}

In the previous simulations, and in the main section of the paper, we
studied graph structure recovery assuming the partial ordering between
the variable blocks of our BDMRF is known. Graph structure recovery when this
partial ordering of the blocked DAG is unknown is a much more
challenging problem. While a complete theoretical and empirical
treatment of this question is beyond the scope of this paper, in this
section, we provide brief numerical results demonstrating that even
under this setting, recovering some dependencies is achievable via
node-neighborhood estimation. In 
Figure~\ref{fig:sim_unknown}, we plot ROC curves for graph structure
recovery for BDMRFs with two blocks of variables, binary
(corresponding to Ising CRF and MRF components), and continuous 
(corresponding Gaussian CRF and MRF components), for the same three
BDMRF instances and simulations settings as in the top panel of
Figure~\ref{fig:sim2way}.  Here however, for each of the model
instances, we do not just present results for the $M$-estimator
corresponding to the specific model instance, but for three different
$M$-estimators corresponding to the three different model
instances. As seen in the figure, the $M$-estimator corresponding to
the true model is not always the best estimator under extremely
sample-limited-settings (though we might expect it to be so as 
$n$ increases). As the plots in this brief investigation show, the
$M$-estimator corresponding to the mixed MRF model may be the most
robust to model misspecification, which in turn indicates that we
could use this estimator to learn dependencies even when the true
partial ordering between the variable blocks may be unknown or
misspecified. A full study extending these preliminary results  is
left as a promising area of future work.

\section{Case Study: High-Throughput Genomics}\label{Sec:RealData}

\begin{figure}[!!t]
\centering
\includegraphics[width=5in]{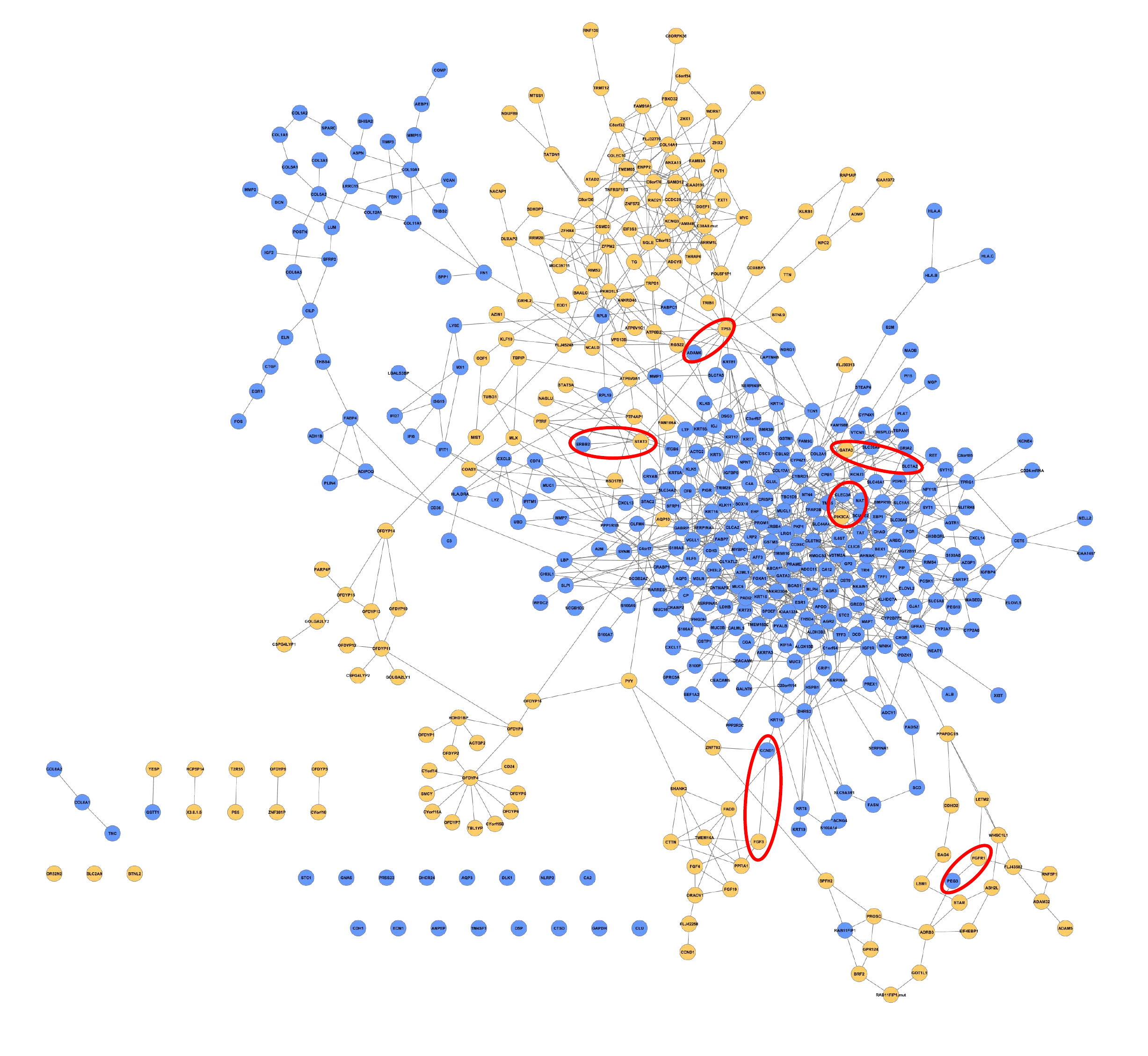}
\caption{Invasive breast carcinoma network of gene expression
  (blue; counts via RNA-sequencing; $Y$) and genomic aberrations
  (yellow; binary via mutations and copy number aberrations; $X$) 
  estimated via our BDMRF distribution.  We learn the graph structure
  assuming the 
  components $\Pe[Y | X]$ and $\Pe[X]$ are specified by a sub-linear
  Poisson CRF and an Ising MRF 
  respectively.  Our method discovered several biologically important
  connections between different types of biomarkers (circled in red).}
\label{fig:genomics}
\end{figure}

We apply our class of BDMRF models to high-throughput cancer genomics
data to find connections both within and between different types of
biomarkers.  Specifically, we study how genomic aberrations, which
include both SNPs and copy number variations, affect gene expression
levels. Note that the standard method for finding connections between
mutations and expression levels, commonly called expression
Quantitative Trait Loci (eQTL) analysis \citep{lee2010adaptive}, uses
regression models to find connections between the two different types of
biomarkers but cannot model the connections within the set of
mutations and gene expression levels.

For our analysis, we use publicly available data on invasive breast
carcinoma available from the Cancer Genome Atlas
(TCGA) \citep{tcga_breast_2012}.  Level III RNA-sequencing data for 806
patients was downloaded and pre-processed using techniques described
in \citet{allen2013local}, so that the expression levels can be
well-modeled with the Poisson distribution.  For the aberration data,
we used Level II non-silent somatic mutations and Level III copy
number variation data for 951 patients.  The later was segmented using
standard techniques \citep{CNTools} and merged with the mutation data
at the gene level to form a binary matrix, indicating whether a
mutation or copy number aberration is present or absent in the coding
region of the gene. This leaves us with $n= 697$ patients that are
common to both the gene expression, $Y$, and aberration, $X$, data
sets.  We 
filtered the biomarkers further to consider the top 2\% of genes whose
expression levels had the highest variance across patients, $p_{X} =
329$, and only the aberrations present in at least 10\% of patients,
$p_{Y} = 169$, yielding $p = 498$ total biomarkers.

We fit our BDMRF model to this data to learn the structure of the
breast cancer genomic network. As we know that aberrations affect
expression levels but not the converse, we set the partial ordering of
the mixed graph underlying our BDMRF to have the aberration variable
block ordered before the expression level variable block.
Specifically, we used the factorization $\P[X,Y] = \Pe[Y | X] \Pe[X]$,
where $\Pe[Y | X]$ is a pairwise Poisson conditional random field as
RNA-sequencing data is count-valued and $\Pe[X]$ is a pair-wise Ising
model as aberration data is binary. Recall  that for the Poisson MRF
or CRF, only negative conditional dependencies are permitted, which is
unrealistic for genomics data. Thus, instead of the usual Poisson CRF,
we fit a Sub-Linear Poisson CRF, an extension of the Sub-Linear
Poisson MRF model described in \citet{YRAL13PGM}, which uses a
sub-linear function as the exponential family sufficient statistic to
relax the resulting normalizability conditions, permitting both
positive and negative conditional dependencies; see \citet{YRAL13PGM}
for further details. Node-wise neighborhood selection as described in
Section~\ref{Sec:Learning} was employed to learn the edge structure of
the network.  Stability 
selection, as described in \citet{liu_STARS_2010}, was used with
parameter level $\beta = 0.01$ to determine the optimal level of
regularization.

Our estimated network is presented in Figure~\ref{fig:genomics} where
blue nodes denote gene expression biomarkers and yellow nodes denote
aberration biomarkers.  
Our network identified several key connections
between biomarkers of 
different types; these are denoted via red circles in
Figure~\ref{fig:genomics}. These include both links that have been
previously indicated in the literature, as well as some novel
discoveries.  
First, several connections that we identify are well-known breast
cancer biomarkers: the GATA3 mutation is linked to SLC39A6 expression;
the ratio of these gene's expressions levels are used to define breast
cancer sub-types \citep{kapp2006discovery} and both of these biomarkers
have been previously implicated in breast cancer \citep{voduc2008gata}.
The FGFR1 mutation is linked to PEG3 
expression; the former  regulates growth factors that are known to be amplified in breast cancer \citep{hynes2010potential} while the latter modulates the related process of cancer  progression \citep{su1999peg}.  The STAT3 mutation is linked to ERBB2 expression; these are known to be amplified in HERB2 sub-types \citep{chung2014stat3}.

Our estimated network also discovers several novel connections to be
investigated and further validated in future work; these include: the
TP53 mutation is linked to ADAM6 expression; TP53 is a well known
tumor suppressor gene \citep{levine1991p53} and ADAM6 is a long
non-coding RNA over-expressed in breast cancer \citep{seals2003adams}.
The FGF3 mutation is linked to CCND1 expression; FGF3 regulates
estrogen expanding breast cancer stem cells
\citep{fillmore2010estrogen} and CDN1 leads to over-expression of
hormone receptors in breast cancer \citep{arnold2005cyclin}.  The
PIK3CA mutation is linked to CLEC3A expression and NAT1 expression;
PIK3CA is a known oncogene 
\citep{cizkova2012pik3ca}, CLEC3A affects tumor metastasis
\citep{tsunezumi2009matrilysin}, and NAT1 is a potential marker for the
estrogen receptor positive sub-type \citep{kim2008promoter}.

Overall, this genomics example demonstrates the direct applicability
of our class of BDMRF models for learning relevant connections both
within and between cancer biomarkers of different 
types.

\section{Discussion}

In this paper, we have constructed a general class of mixed graphical
models through block directed exponential family mixed conditional
random fields. Our work generalizes much of the 
existing literature on parametric Markov Random Field (MRF) and
Conditional Random Field (CRF)  
distributions, all of which are special cases of our framework. Our
so-called class of Block Directed Markov Random Fields (BDMRFs) are an
extremely flexible class of models, that are represented by mixed
graphs, with both directed edges (through a DAG over blocks of
variables) and undirected edges (within these blocks), between mixed
types of variables. We have also shown that these are normalizable
under fairly relaxed conditions.

To our knowledge, these models represent the first class of
multivariate densities over mixed variables that directly permits and
parameterizes such a rich set of dependencies. Our work thus has broad
implications for multivariate analysis in general, especially that
involving integrative analysis of mixed variables. At the very least,
we have greatly expanded the class of off-the-shelf graphical models
beyond canonical instances such as Ising and Gaussian graphical
models, and have also provided estimators for our class of models that
come with strong statistical guarantees even  under high-dimensional
settings.

While the main contributions of this paper are theoretical in nature,
our novel class of BDMRF distributions has broad applicability in a
variety of fields that yield mixed, big data. Beyond the
high-throughput genomics example discussed in this paper, these
include imaging genetics, national security, climate studies, spatial
statistics, Internet data, marketing and advertising, and economics,
among many others.

There are several avenues for further research building on the work in
this paper. These include providing a goodness of fit test for understanding
whether our models are applicable to any particular data
set. \citet{YRAL13} discuss a heuristic for this, but further
work needs to be done to develop a proper post selection inferential
procedure.  While we have shown that our BDMRF distributions impose
fairly relaxed conditions upon their parameters for normalizability,
these conditions 
can still be limiting in practice, and more work is needed on general
strategies for relaxing these 
normalizability conditions; for instance, this could be achieved by
potentially extending  
the work in \citep{YRAL13PGM} on relaxing normalizability
conditions for Poisson MRFs. Additionally, our work has concentrated
on constructing models via single-parameter univariate exponential 
family distributions, which can be further extended to the case of
multi-parameter exponential families. Finally, in this paper, we
assumed that the partial ordering or DAG over the variable blocks is
known a priori based on domain knowledge.  For mixed genomics data,
this assumption is realistic as we know how various biomarkers affect
each other in the biological system. For 
other application areas, the directionality of links between blocks of
variables may be unknown or the block designations of the variables
may be unknown.  In these situations, the question remains: can we
still learn the dependence structure between mixed variables and/or
estimate both the edges are their directionality?  Our final set of
simulations in Figure~\ref{fig:sim3way} suggests that we may be able to
learn dependencies even if the true underlying BDMRF is unknown;
further theoretical and empirical investigations that build on these
are needed. Also, while there are several popular algorithms for
learning DAGs~\citep{kalisch2007estimating}, little is known about 
estimating MRFs with mixed directed and undirected edges, an open area
of future research.   

Overall, we have proposed a novel class of distributions for mixed
variables, that has broad theoretical and practical implications
across many statistical and scientific domains.

\bibliographystyle{abbrvnat}
\bibliography{mixedGM,network,sml}

\newpage
\section{Appendix}
\subsection{Proof of Theorem~\ref{ThmExpMixedCRF}}\label{SecProofMixedCRF}

	We follow the proof policy of \citep{YRAL12,YRAL13CRF}: Define $Q(\by | X)$ as
    \begin{align*}
    	Q(\by | X) := \log(\P[\by | X ] / \P[{\bf 0} | X]),
    \end{align*}
    for any $\by = (Y_1,\hdots,Y_p) \in \mathcal{Y}^{p}$ given $X$ where ${\bf 0}$ indicates a zero vector (The number of zeros vary appropriately in the context below). For any $\by$, also denote $\bys := (Y_1,\hdots,Y_{s-1},0,Y_{s+1},\hdots,Y_p)$. 
    
    Now, consider the following general form for $Q(\by|X)$:
    \begin{align}\label{EqnQFactorForm}
    	Q(\by | X) = & \sum_{t_1 \in V} Y_{t_1} G_{t_1}\big(Y_{t_1},X\big) + \hdots + \\
            & \sum_{t_1,\hdots,t_k \in V}Y_{t_1}\hdots Y_{t_k} G_{t_1,\hdots,t_k}\big(Y_{t_1},\hdots,Y_{t_k},X\big), \nonumber
    \end{align}
    where $k$ is the maximum clique size of $G_Y$. It can then be seen that
    \begin{align}\label{EqnTmpA}
    	& \exp\Big(Q(\by|X) - Q(\bys|X)\Big) = \P[\by|X]/ \P[\bys|X]\nonumber\\
    	 						=\ & \frac{\P[Y_s|Y_1,\hdots,Y_{s-1},Y_{s+1},\hdots,Y_p,X]}{\P[0|Y_1,\hdots,Y_{s-1},Y_{s+1},\hdots,Y_p,X]},
    \end{align}
    where the first equality follows from the definition of $Q$, and the second equality follows from some algebra. Now, consider simplifications of both sides of \eqref{EqnTmpA}. Given the form of $Q(\by|X)$ in \eqref{EqnQFactorForm}, we have
    \begin{align}\label{EqnTmpB}
    	& Q(\by|X) - Q(\byi|X) = \\
        & Y_{1}\bigg(G_1\big(Y_1,X\big) + \sum_{t=2}^{p} Y_t G_{1t}\big(Y_1, Y_t,X\big) + \hdots + \nonumber\\
        & \sum_{t_2,\hdots,t_k \in \{2,\hdots,p\}}Y_{t_2}\hdots Y_{t_k} G_{1,t_2,\hdots,t_k}\big(Y_1,\hdots,Y_{t_k},X\big)\bigg).\nonumber
    \end{align}
    Also, given the exponential family form of the node-conditional distribution specified in \eqref{EqnExpMixedCRF_Cond}, we obtain
    \begin{align}\label{EqnTmpC}
        &\log \frac{\P[Y_s|Y_1,\hdots,Y_{s-1},Y_{s+1},\hdots,Y_p,X]}{\P[0|Y_1,\hdots,Y_{s-1},Y_{s+1},\hdots,Y_p,X]} = \\
        &E_s(Y_{V \backslash s},X)(B_s(Y_s) - B_s(0)) + (C_s(Y_s) - C_s(0)).\nonumber
    \end{align}
    Setting $Y_{t} = 0$ for all $t \neq s$ in \eqref{EqnTmpA}, and using the expressions for the left and right hand sides in \eqref{EqnTmpB}
    and \eqref{EqnTmpC}, we obtain,
    \begin{align*}
    	&Y_{s} G_s\big(Y_s,X\big) \\
        = \ &E_s({\bf 0},X)(B_s(Y_s) - B_s(0)) + (C_s(Y_s) - C_s(0)).
    \end{align*}
    Setting $Y_{r} = 0$ for all $r \not\in \{s,t\}$,
    \begin{align*}
    	& Y_{s} G_s\big(Y_s,X\big) + Y_{s} Y_{t} G_{st}\big(Y_s,Y_t,X\big) \\
        =\ & E_s({\bf 0},Y_t,{\bf 0},X)(B_s(Y_s) - B_s(0)) + (C_s(Y_s) - C_s(0)).
    \end{align*}
    Combining these two equations yields
    \begin{align}\label{EqnGst1}
        &Y_{s} Y_{t} G_{st}\big(Y_s,Y_t,X\big) \nonumber\\
        =\ &\big(E_s({\bf 0},Y_t,{\bf 0},X)-E_s({\bf 0},X)\big)(B_s(Y_s) - B_s(0)).
    \end{align}
    Similarly, from the same reasoning for node $t$, we have
    \begin{align*}
    	& Y_{t} G_t\big(Y_t,X\big) + Y_{s} Y_{t} G_{st}\big(Y_s,Y_t,X\big)  \\
        = \ & E_t({\bf 0},Y_s,{\bf 0},X)(B_t(Y_t) - B_t(0)) + (C_t(Y_t) - C_t(0)),
    \end{align*}
    and at the same time,
    \begin{align}\label{EqnGst2}
        &Y_{s} Y_{t} G_{st}\big(Y_s,Y_t,X\big) \nonumber\\
        =\ &\big(E_t({\bf 0},Y_s,{\bf 0},X)-E_t({\bf 0},X)\big)(B_t(Y_t) - B_t(0)).
    \end{align}
    From the assumption of the statement of \eqref{CondMixedCRFMarkov}, $E_s(\cdot)$ depends $X$ only through $N_{YX}(s)$, and similarly $E_t(\cdot)$ does only through $N_{YX}(t)$. Furthermore, by the equality of \eqref{EqnGst1} and \eqref{EqnGst2}, we obtain
    \begin{align}\label{EqnEtmp}
        &E_t({\bf 0},Y_s,{\bf 0},X)-E_t({\bf 0},X) = \frac{E_s({\bf 0},Y_t,{\bf 0},X)-E_s({\bf 0},X)}{B_t(Y_t) - B_t(0)}(B_s(Y_s) - B_s(0)).
    \end{align} 
    %
    %
    Since \eqref{EqnEtmp} should hold for all possible combinations of $Y_s$, $Y_t$ and $N_{YX}(s) \cap N_{YX}(t)$, for any fixed $Y_t \neq 0$,
    \begin{align}\label{EqnE}
        &E_t({\bf 0},Y_s,{\bf 0},X)-E_t({\bf 0},X) = \theta_{st}\big(N_{YX}(s) \cap N_{YX}(t)\big)(B_s(Y_s) - B_s(0))
    \end{align}
    where $\theta_{st}(\cdot)$ is a function on $N_{YX}(s) \cap N_{YX}(t)$. Plugging \eqref{EqnE} back into \eqref{EqnGst2},
    \begin{align*}
    	&Y_{s} Y_{t} G_{st}\big(Y_s,Y_t,X\big) \\
        =\ & \theta_{st}\big(N_{YX}(s) \cap N_{YX}(t)\big) (B_s(Y_s) - B_s(0)) (B_t(Y_t) - B_t(0)).
    \end{align*}
    More generally, by considering non-zero triplets, and setting $Y_{r} = 0$ for all $r \not\in \{s,t,u\}$, we obtain,
    \begin{align*}
    	&Y_{s} G_s\big(Y_s,X\big) + Y_{s} Y_{t} G_{st}\big(Y_s,Y_t,X\big)\\
        & + Y_{s} Y_{u} G_{su}\big(Y_s,Y_u,X) + Y_{s} Y_{t} Y_{u} G_{stu}\big(Y_s,Y_t,Y_u,X\big) \\
         = \, & E_s({\bf 0},Y_t,{\bf 0},Y_u,{\bf 0},X)(B_s(Y_s) - B_s(0)) \\
            & + (C_s(Y_s) - C_s(0)),
    \end{align*}
    so that by a similar reasoning we can obtain
    {\small
    \begin{align*}
    	&Y_{s} Y_{t} Y_{u} G_{stu}\big(Y_s,Y_t,Y_u,X\big) = \\
        &\theta_{stu}\big(N_{YX}(s) \cap N_{YX}(t) \cap N_{YX}(u)\big)(B_s(Y_{s}) - B_s(0))(B_t(Y_{t}) - B_t(0))(B_u(Y_{u}) - B_u(0)).
    \end{align*}
    }More generally, we can show that
    \begin{align*}
    	&Y_{t_1}\hdots Y_{t_k} G_{t_1,\hdots,t_k}\big(Y_{t_1},\hdots,Y_{t_k},N_Y(X)\big) = \\
        &\theta_{t_1,\hdots,t_k}\big(\cap_{t = t_1,\hdots,t_k} N_{YX}(t)\big)(B_{t_1}(Y_{t_1}) - B_{t_1}(0))\hdots(B_{t_k}(Y_{t_k}) - B_{t_k}(0)).
    \end{align*}
    Thus, the $k$-th order factors in the joint distribution as specified in \eqref{EqnQFactorForm} are tensor products of $(B_s(Y_s) - B_s(0))$, thus proving the statement of the theorem.

\subsection{Proof of Theorem~\ref{ThmMarkovProp} and \ref{ThmRecurMarkovProp}}\label{SecProofMakrovProp}
The joint distribution \eqref{EqnElemBDMRF} factors according to $E_{X}$ and $E_{Y}$ by construction. Hence, (a) and (b) hold by Hammersley-Clifford theorem. Now, suppose that $\{d\} \cup C$ is not complete with respect to $G$. Then, since $C$ is complete and every node in $d$ is fully connected to $C$, we can conclude that again by Hammersley-Clifford theorem, $\theta_{C,d}(X)$ is decomposable $\sum_{j=1}^i \theta_{C,d_j}(X)$ where $d_1,\hdots,d_i$ are the disjoint sets such that $\cup_{j=1,\hdots,i} {d_j} = d$ and $d_j \cup C$ is complete with respect to $G$ for all $j$, which means $\theta_{C,d}(X) = 0$. 

The proof of Theorem \ref{ThmRecurMarkovProp} can be trivially derived by extending this proof.

\subsection{Proof of Theorem~\ref{ThmEBDMRFNormalizability}}\label{SecProofEBDMRFNormalizability}
We first introduce shorthand notations for unnormalized probabilities:  
\begin{align*}
	&U_X = \exp\Big( \sum_{C \in \C_{X}} \theta_{C} \prod_{t \in C} B_{t}(X_t) + \sum_{t \in V_X} M_{t}(X_{t}) \Big) \, , \nonumber\\
	&U_{XY} = \exp\Big( \sum_{C \in \C_{Y}} \theta_{C}\big( X_{N_{YX}(C)}  \big) \prod_{s \in C} B_{s}(Y_s) + \sum_{s \in V_Y} M_{s}(Y_{s}) - A_{Y|X}\big(\theta\big) \Big) \, .
\end{align*}
Suppose that $\Pe[X]$ and $\Pe[Y|X]$ are well-defined and normalizable, which means that $\sum_{X} U_X < \infty$, and $\sum_{Y} U_{XY} < \infty$ for any $X$. Note that $\sum$ operator can be substituted by integral for continuous case. 
Then, the log partition function of \eqref{EqnElemBDMRF} can be shown to be finite as follows:
\begin{align*}
	&\sum_{X,Y} \Big[ U_X U_{XY} \Big] = \sum_{X} \Big[ U_X \big( \sum_{Y} U_{XY} \big) \Big] \nonumber\\
	\leq \, & \sum_{X} \Big[ U_X \max_{X} \big( \sum_{Y} U_{XY} \big) \Big] =  \Big[ \max_{X}  \sum_{Y} U_{XY} \Big] \, \Big[ \sum_{X} U_X \Big] < \infty.
\end{align*}

\subsection{Proof of Theorem~\ref{Prop:MRF-RecMRF}}\label{Proof:Prop:MRF-RecMRF}
Since the unnormalized probability mass $\exp\{ F(X;\theta) + F(X,Y;\theta) \}$ is strictly positive, $\exp\{ F(X;\theta)\} >0$ and $\exp\{F(X,Y;\theta) \} >0$ for any $(X,Y) \in \mathcal{X}^q \times \mathcal{Y}^p$.
	 
Suppose that $A(\theta) = \log \sum_{X} \Big[\exp \{ F(X;\theta) \} \sum_Y \exp\{ F(X,Y;\theta) \} \Big] < \infty$ as in the statement. Then, for any $X \in \mathcal{X}^q$, the partial summation of $A(\theta)$, $\exp \{ F(X;\theta) \} \sum_Y \exp\{ F(X,Y;\theta) \}$ is also summable;
\[
		\exp \{ F(X;\theta) \} \sum_Y \exp\{ F(X,Y;\theta) \} = c_1,
\] which implies that $\sum_Y \exp\{ F(X,Y;\theta) \} = c_1 / \exp \{ F(X;\theta) \} < \infty$, and $A_{Y|X}\big(\vartheta(X)\big)$ is summable. 
	
Now, we consider the second term $A_X(\theta)$ in \eqref{EqnEBDMRFPairwise}. Since $\exp\{F(X,Y;\theta) \} >0$ for any $X \in \mathcal{X}^q$,  $\min_{X} \sum_Y \exp\{ F(X,Y;\theta) \}$ is strictly positive; let us define $c_{\min} := \min_{X} \sum_Y \exp\{ F(X,Y;\theta) \} >0$. We then have:
\begin{align*}
	& c_{\min} \sum_X \exp \{ F(X;\theta) \} = \sum_X \Big[c_{\min} \exp \{ F(X;\theta) \}\Big]  \\
	\leq \, & \sum_{X} \Big[\exp \{ F(X;\theta) \} \sum_Y \exp\{ F(X,Y;\theta) \} \Big].  
\end{align*}
Therefore, 
\[
	\sum_X \exp \{ F(X;\theta) \} \leq \frac{1}{c_{\min}} \sum_{X} \Big[\exp \{ F(X;\theta) \} \sum_Y \exp\{ F(X,Y;\theta) \} \Big].
\]
Since $c_{\min} >0$ and $\Big[\exp \{ F(X;\theta) \} \sum_Y \exp\{ F(X,Y;\theta) \} \Big] < \infty$ by assumption, we can conclude $A_X(\theta) < \infty$, which completes the proof.

\subsection{Normalizability of Mixed CRFs}\label{Sec:NormalizabilityMixedCRFs}


The normalizability conditions for the mixed MRF $\Pe[X]$ is well
developed in \citet{YBRAL14}. 
This can be seamlessly extended for the mixed CRF
\eqref{EqnExpMixedCRF} case as well. 
\begin{proposition}\label{Prop_Finite}
	Suppose that some random variables in $Y$ are finite. Let $Z$
        be a set of such random variables: $\mathcal{Z}_{s}$ is
        finite, with $\max\{z: z \in \mathcal{Z}_{s}\} < \infty$ and
        $\min\{z: z \in \mathcal{Z}_{s}\} > - \infty$. We also define
        $W$ to denote the remaining (infinite) random variables in $Y$
        so that $Y 
        = (Z,W)$. Suppose that the following conditional distribution
        given $Z$ is well-defined (i.e. normalizable) for all possible
        $Z \in \prod_{s \in V_Z} \mathcal{Z}_{r}$: 
	\begin{align*}
		\log \Pe[ W | Z, X]  \propto \sum_{ C \in \C_W} \theta_C(Z, X) \prod_{s \in C} B_{s}(W_s) + \sum_{s \in V_W} M_s(W_s) . 
	\end{align*}
	 Then, $\Pe[ W | Z, X]$ is normalizable if and only if mixed CRF $\P_{\exp}[Y|X]$ \eqref{EqnExpMixedCRF} is normalizable.
\end{proposition}

\begin{proposition}\label{Prop_Infinite}
	Suppose that all domains $\mathcal{Y}_s$ for $Y$  are
        infinite, and that the maximum size of cliques is limited to two
        with linear sufficient statistics $B_s(Y_s)=Y_s$ for $\forall
        s \in V_Y$. Then the mixed CRF in \eqref{EqnExpMixedCRF} is
        normalizable only if at least one of the following conditions
        is satisfied for all $s,t \in V_Y$ with non-zero
        $\theta_{st}(X)$: 
		\begin{description}[leftmargin=0.5cm] 
			\item[(a)] Both $\mathcal{Y}_s$ and
                          $\mathcal{Y}_t$ are infinite only from one
                          side, so that $\sup\{y: y \in
                          \mathcal{Y}_s\} < \infty$ or $\inf\{y: y \in
                          \mathcal{Y}_s\} > - \infty$, and similarly
                          for $\mathcal{Y}_t$.  
			\item[(b)] For all $\alpha, \beta > 0$
                          such that $- M_{s}(Y_s) = O(Y_s^\alpha)$ and
                          $ - M_{t}(Y_t) = O(Y_t^\beta)$, it holds
                          that $(\alpha-1)(\beta-1) \geq 1$. 
		\end{description}
\end{proposition}

Even though Proposition \ref{Prop_Finite} and \ref{Prop_Infinite} are
simply the generalization of analogous results on normalizability
conditions of mixed MRF~\citep{YBRAL14}, they are useful to study to better
understand why our BDMRF formulations are less restrictive 
than mixed MRFs. In Proposition~\ref{Prop_Infinite},
consider that conditions (a) and (b) only depend on the response
variables $Y$ instead of both $X$ and $Y$.  This can dramatically
relax normalizability restrictions in many settings.  For example,
consider the Gaussian-Poisson EBDMRF studied in
Section~\ref{SecEBDMRFExamples} for Poisson-associated variables $X$
and Gaussian-associated variables $Y$.  For the mixed MRF, $\Pe[X,Y]$,
we would consider both $X$ and $Y$ as response variables and neither
(a) nor  (b) in Proposition~\ref{Prop_Infinite} hold. (For (a),
 $Y_s$ is \emph{not} infinite only from one side, and for (b), $-M_{s}(Y_s) =
O(Y_s^2)$ and $-M_{t}(X_t) = O(X_t^{1.5})$, $(2-1)(1.5-1)$ is
\emph{not} larger than $1$). Hence, the Gaussian-Poisson mixed MRF is
not normalizable.

On the other hand, both of our EBDMRF formulations are normalizable
because of the relaxed restrictions with our mixed CRFs.
Consider taking the variables $X$ and $Y$ to be two homogeneous
blocks.  Then the conditions in Proposition~\ref{Prop_Infinite} need
only to be satisfied for $X$ {\em or} $Y$.  For the Gaussian-Poisson
case, the statement in (a) holds for the Poisson CRF, $\Pe[X|Y]$ as
all Poisson variables $X$ are unbounded in only one
direction. Similarly, statement (b) holds for the Gaussian CRF,
$\Pe[Y|X]$ since $-M_{s}(Y_s) = O(Y_s^2)$ for all Gaussian variable
$Y_s$.  (Further conditions on the normalizability of these Poisson and
Gaussian CRFs are given in Section~\ref{SecEBDMRFExamples}).  Hence,
{\em both} the Poisson and Gaussian CRFs are normalizable meaning that
both of our EBDMRF formulations, $\Pe[X|Y]\Pe[Y]$ and
$\Pe[Y|X]\Pe[X]$, are normalizable and permit dependencies in
Gaussian-Poisson mixed graphical models.

\subsection{Simulation Settings}\label{AppSec:SimSettings}

\begin{figure}[t]
	  \centering
	  \includegraphics[width=.4\linewidth]{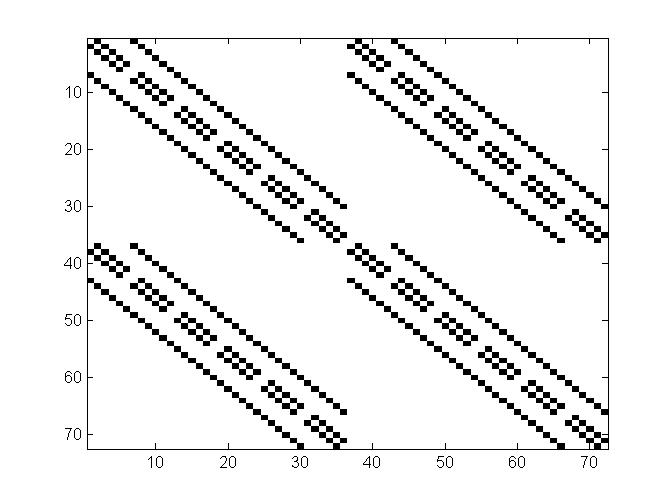}
	  \caption{Lattice Graph Structure}
	  \label{fig:lattice}
\end{figure}

\begin{table}[!t]
	\caption{Parameters used in data generation for various examples of homogeneous mixed MRFs and our EBDMRFs in Figure \ref{fig:sim2way} whose form is specified by the MRF and CRF component as labeled.} \label{tab:2waytable}
	\centering
	\begin{tabular}{c|c c c}
	Model & $\theta_{xx}$ & $\theta_{yy}$ & $\theta_{xy}$\tabularnewline
	\hline 
	Gaussian (X) - Ising (Y) Mixed MRF & 0.2 & 0.3 & 0.1\tabularnewline
	Gaussian MRF - Ising CRF & 0.5 & 0.3 & 0.1\tabularnewline
	Gaussian CRF - Ising MRF & 0.5 & 0.3 & 0.1\tabularnewline
	Poisson (X) - Ising (Y) Mixed MRF & -0.8 & 0.4 & 0.6\tabularnewline
	Poisson MRF - Ising CRF & -0.8 & 0.4 & 0.6\tabularnewline
	Poisson CRF - Ising MRF & -0.8 & 0.4 & 0.6\tabularnewline
	Gaussian (X) CRF - TPGM (Y) MRF & 0.1 & 0.1 & 0.2\tabularnewline
	Gaussian (X) MRF - Poisson (Y) CRF & 0.6 & -1 & 1\tabularnewline
	Gaussian (X) CRF - Poisson (Y) MRF & 0.6 & -1 & 1\tabularnewline
	\end{tabular}
\end{table}

\begin{table}[!t]
	\caption{Parameters used in data generation for various configurations of our BDMRFs in Figure \ref{fig:sim3way}, for three homogeneous blocks of variables: binary (Ising, X), continuous (Gaussian, Y) and counts (Poisson, Z).}\label{tab:3waytable}
	\centering
	\begin{tabular}{ c| c c c c c c }
	Model & $\theta_{xx}$ & $\theta_{yy}$ & $\theta_{zz}$ & $\theta_{xy}$ & $\theta_{xz}$ & $\theta_{yz}$\tabularnewline
	\hline 
	 $\Pe[X|Y,Z]\, \Pe[Y|Z]\, \Pe[Z]$ & 0.1 & 0.3 & -0.8 & 0.1 & 0.1 & 0.1\tabularnewline
	 $\Pe[Y|X,Z]\, \Pe[X|Z]\, \Pe[Z]$ & 0.1 & 0.3 & -0.8 & 0.1 & 0.1 & 0.1\tabularnewline
	 $\Pe[Y|X]\, \Pe[X|Z]\, \Pe[Z]$ & 0.1 & 0.3 & -0.8 & 0.1 & - & 0.1\tabularnewline
	\end{tabular}
\end{table}

We generate samples from our BDMRF models using Gibbs sampling for the
lattice network structure shown in Figure~\ref{fig:lattice}.
Parameters for each model were specified as constant edge weights for
each block corresponding to the values given in Table~\ref{tab:2waytable}
and Table~\ref{tab:3waytable}.

\subsection{Genomics Case Study Comparisons}

\begin{figure}
\subfigure[EBDMRF\label{fig:circ_chained}]
  {\includegraphics[trim = 3cm 15cm 5cm 0cm, clip, width=.5\textwidth]{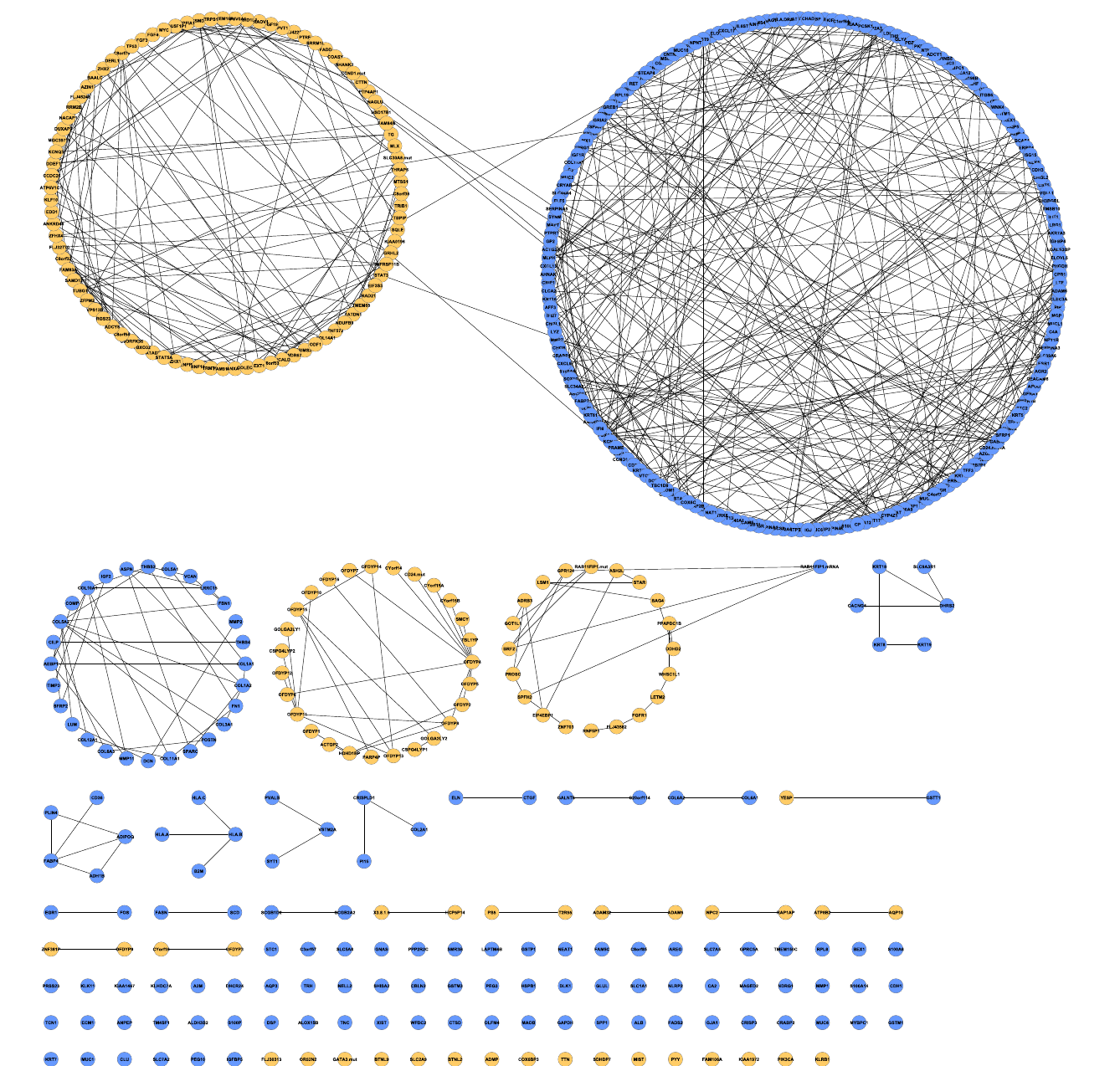}}\hfill
\subfigure[Mixed MRF\label{fig:circ_mixedMRF}]
  {\includegraphics[trim = 0cm 5cm 0cm 0cm, clip,width=.5\textwidth]{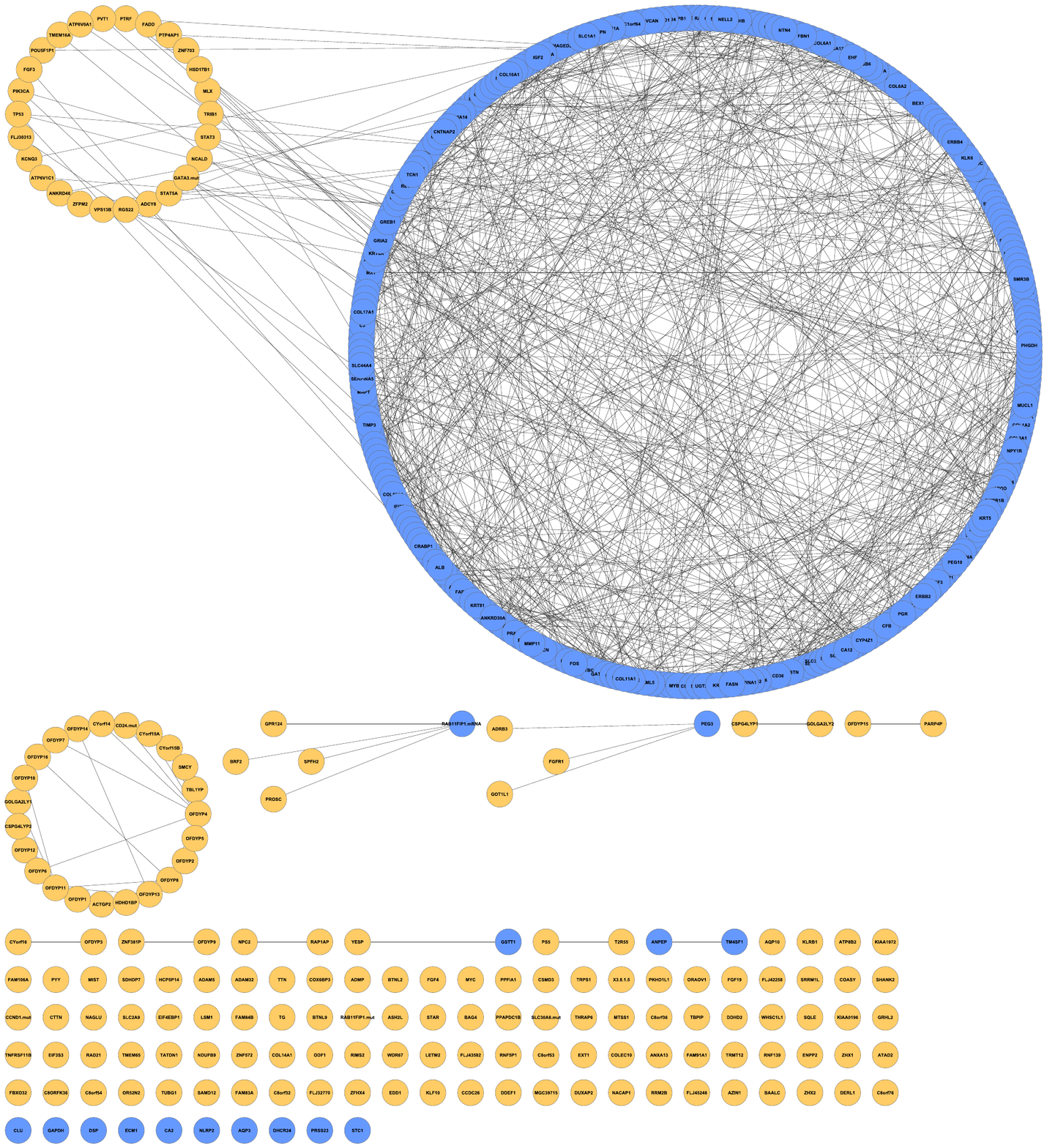}}\hfill
\caption{Comparison of breast cancer genomics networks estimated via
  our EBDMRF vs. the mixed MRF.}
\label{fig:circle}
\end{figure}

In circle plots in 
Figure~\ref{fig:circle}, we 
compare our EBDMRF network estimate to that of the mixed MRF network
studied 
in \cite{YBRAL14} and reviewed in Section~\ref{Sec:Background}.  We
see that the mixed MRF 
finds many connections within gene expression biomarkers, but few
within aberration biomarkers. 
In contrast, our BDMRF model finds both within connections as well as
a few key between connections.  Biologically, this result is not
surprising as gene expression levels are known to be highly dependent,
and thus when conditioning on all biomarkers as in the mixed MRF,
connections within gene expressions dominate.  With our network,
conditioning based on the directions of known biological influence
(i.e. aberrations influence gene expression), allows us to estimate a
more biologically meaningful network.

\end{document}